\documentclass[12pt]{article}
\usepackage{geometry}
\usepackage{graphicx}
\usepackage{amsmath, amsfonts, amsthm, amssymb}
\usepackage{url}
\usepackage[dvipsnames]{xcolor} 
\usepackage{tikz}
\usetikzlibrary{arrows.meta,positioning,calc}
\usepackage{aliascnt}
\usepackage[hidelinks]{hyperref}
\usepackage[nameinlink]{cleveref}
\usepackage{microtype}
\usepackage{comment}
\usepackage[section]{placeins}
\usepackage{thmtools}
\usepackage{thm-restate}

\theoremstyle{definition}

\theoremstyle{definition}

\newtheorem{theorem}{Theorem}[section]

\theoremstyle{definition}

\newtheorem{lemma}[theorem]{Lemma}

\newtheorem{corollary}[theorem]{Corollary}

\newtheorem{definition}[theorem]{Definition}

\newtheorem{proposition}[theorem]{Proposition}

\newtheorem{problem}[theorem]{Problem}

\newtheorem{example}[theorem]{Example}

\newtheorem{remark}[theorem]{Remark}

\crefname{theorem}{theorem}{theorems}
\Crefname{theorem}{Theorem}{Theorems}
\crefname{lemma}{lemma}{lemmas}
\Crefname{lemma}{Lemma}{Lemmas}
\crefname{corollary}{corollary}{corollaries}
\Crefname{corollary}{Corollary}{Corollaries}
\crefname{definition}{definition}{definitions}
\Crefname{definition}{Definition}{Definitions}
\crefname{proposition}{proposition}{propositions}
\Crefname{proposition}{Proposition}{Propositions}
\crefname{remark}{remark}{remarks}
\Crefname{remark}{Remark}{Remarks}
\crefname{example}{example}{examples}
\Crefname{example}{Example}{Examples}
\crefname{problem}{problem}{problems}
\Crefname{problem}{Problem}{Problems}
\crefname{appendix}{appendix}{appendices}
\Crefname{appendix}{Appendix}{Appendices}
\crefname{maintheorem}{theorem}{theorems}
\Crefname{maintheorem}{Theorem}{Theorems}

\usepackage{enumitem}

\usepackage[normalem]{ulem}
\newif\ifshowrevisions
\showrevisionstrue            
\definecolor{revref}{rgb}{0,0,1}   
\definecolor{revnew}{rgb}{1,0,0}   
\ifshowrevisions

  \newcommand{\refdel}[1]{{\color{revref}\sout{#1}}}   
  \newcommand{\newdel}[1]{{\color{revnew}\sout{#1}}}
  \newenvironment{refblockrev}{\color{revref}}{}
  \newenvironment{newblockrev}{\color{revnew}}{}
\else

  \newcommand{\refdel}[1]{}
  \newcommand{\newdel}[1]{}

\fi

\DeclareMathOperator{\ReLU}{ReLU}
\newcommand{\act}{\varsigma} 

\newcommand{\Lam}{\Lambda}                 
\newcommand{\dimfunzero}{\dim^{0}_{\mathrm{fun}}}

\newcommand{\clockset}[1]{\mathrm{Clk}_{#1}}        
\newcommand{\clocksign}{\epsilon}                   
\newcommand{\adjdiff}[1]{\Delta_{#1}}               
\newcommand{\slopeembed}[1]{\iota_{#1}}             
\newcommand{\rankonevar}{\mathcal{R}}               

\usepackage{framed} 
\usepackage{subcaption} 
\DeclareCaptionFont{revnew}{\ifshowrevisions\color{revnew}\fi} 
\usepackage{booktabs} 
\usepackage{arydshln} 
\usepackage{mathtools} 

\newcommand{\demph}[1]{\emph{#1}}
\title{\bf Canonical activation patterns and ReLU pattern varieties}

\author{Bella Finkel, Yaoying Fu, and Chenxi Wu \medskip\\(with an appendix by Patrick Frem, Eduardo Madero, and Weiya Wang)}
\date{}

\begin{document}

\maketitle

\begin{abstract}
A ReLU pattern variety records the polynomial relations among the coefficients of the affine maps associated with a prescribed set of activation patterns of a ReLU network. 
We study how the combinatorics of the canonical set of activation patterns of an unbiased ReLU network is reflected in the geometry of its pattern variety. In particular, we introduce the notion of pattern stability and characterize pattern stability and combinatorial stability for shallow networks. 
When the associated wall graph is connected and contains an edge of every neuron label, we resolve a problem of Alexandr and Mont\'ufar by determining the defining ideal of the associated pattern variety explicitly. We establish a hierarchy of progressively weaker equivalence relations relating labeled combinatorial equivalence, equality and permutation equivalence of canonical pattern sets, and equivalences of the associated varieties. Finally, we relate functional dimension to the differential of the pattern map and show that it is bounded above by the dimension of the canonical pattern variety at every depth. For nondegenerate shallow networks, we give an exact formula for the functional dimension; for every shallow architecture, it equals the dimension of the canonical pattern variety on an open dense set of full Lebesgue measure.
\end{abstract}

\section*{Introduction}

A feedforward ReLU network partitions its input space into polyhedral regions on which it computes affine maps. Without biases, these regions are cones and the maps are linear. The number of regions provides one measure of expressive complexity and has motivated both worst-case bounds and analyses of typical networks \cite{montufar2014regions,hanin2019complexity,hanin2019deep}. To describe the network more thoroughly, however, one must also understand which activation patterns occur, how their regions fit together, and which algebraic relations hold among the corresponding maps.

The representation of piecewise-linear functions has a
long history, including the work of Chua and Deng on canonical piecewise-linear representations \cite{ChuaDeng}. In the context of ReLU networks, Alexandr and Mont\'ufar \cite{alexandr2025constraining} introduced ReLU pattern varieties to study polynomial relations among the maps associated with prescribed activation patterns. In the unbiased setting considered here, we fix a finite pattern set $\mathcal A$ and let $M_A(\theta)$ be the matrix obtained by inserting the activation masks of $A$ into the product of weight matrices. The polynomial map $\varphi^{\mathcal A}$ records these matrices simultaneously, and its image has Zariski closure $U^{\mathcal A}$. The patterns remain fixed while the weights vary over the entire parameter space, including parameters that do not realize those patterns. Thus $U^{\mathcal A}$ records polynomial relations among the prescribed maps, without imposing the inequalities required for all of their patterns to occur.

A parameter $\theta$ determines a canonical choice of
patterns. Its canonical polyhedral complex
$\mathcal C(\theta)$ is obtained by successive subdivision
by hidden-neuron preactivations, and $\mathcal{A}_{\theta}$ consists
of the patterns realized on its full-dimensional cells.
Canonical complexes and their combinatorial structure
have been studied by Grigsby and Lindsey
\cite{grigsby2022transversality} and Masden \cite{Masden2025}. Taking $\mathcal A=\mathcal{A}_{\theta}$ associates a pattern variety with the canonical patterns of an individual parameter. This construction leads to four questions:
\begin{enumerate}[label=(\roman*),leftmargin=2em]
\item When is $\mathcal{A}_{\theta}$ locally constant in parameter
space, and how does this relate to combinatorial stability of $\mathcal C(\theta)$?
\item Which polynomial equations generate the defining
ideal of $U^{\mathcal{A}_{\theta}}$?
\item How are labeled
combinatorial equivalence of canonical complexes,
equality of canonical pattern sets up to neuron
permutations, and the different notions of equivalence
of their ReLU pattern varieties related? Which implications can be reversed, and under what hypotheses?

\item How does the unbiased functional dimension at
$\theta$ relate to $\dim U^{\mathcal{A}_{\theta}}$, and under which
hypotheses are they equal?
\end{enumerate}
Our main results answer these questions for networks
without biases and with the convention that the output layer map is not ReLU-activated.

We first consider stability. A parameter is
\emph{pattern-stable} if its canonical pattern set is
constant throughout some neighborhood
(\Cref{def:pattern_stable_neighborhood}). Combinatorial stability requires the canonical polyhedral complex and its ternary labelings to retain their structure under sufficiently small parameter perturbations.
For shallow networks, these two conditions coincide and admit an exact characterization in terms of the first-layer rows.

\begin{restatable}[Stability of shallow networks]
{maintheorem}{shallowstability}
\label{thm:shallow_pattern_combinatorial_stable}
Let $\theta=(W_1,W_2)\in\mathcal P_{n_0,n_1,n_2}$ be a parameter of an unbiased shallow ReLU network where $n_0,n_1\geq1$, and let $r=\min(n_0,n_1)$.
The following are equivalent:
\begin{enumerate}[label=(\roman*)]
    \item $\theta$ is combinatorially stable;
    \item $\theta$ is pattern-stable;
    \item every set of $r$ rows of $W_1$ is linearly independent.
\end{enumerate}
Let $\mathcal G_{n_0,n_1,n_2}\subseteq
\mathcal P_{n_0,n_1,n_2}$ denote the set of parameters satisfying
the conditions above. The collection $\mathcal G_{n_0,n_1,n_2}$ is
open and dense in $\mathcal P_{n_0,n_1,n_2}$, and its complement
is the zero set of a nonzero polynomial. In particular, the
complement has Lebesgue measure zero.
\end{restatable}

\Cref{thm:shallow_pattern_combinatorial_stable} identifies general position of the first-layer rows as the exact condition for both forms of shallow stability.
At greater depth, the two notions differ: we exhibit
a pattern-stable parameter of architecture $(3,4,1,1)$
that is not combinatorially stable. Pattern stability
nevertheless remains generic. For every architecture,
the pattern-stable parameters form an open dense subset
of full Lebesgue measure (\Cref{thm:pattern_stable_generic}).

Grigsby, Lindsey, Meyerhoff, and Wu conjectured that combinatorial stability holds on a set of full measure for networks with biases \cite[Section~5.2 and Remark~4.16]{GLMWFunctionalDimension}. \Cref{thm:shallow_pattern_combinatorial_stable} proves the unbiased analogue of this conjecture for shallow networks, with an explicit description of the stable locus. For deep networks, \Cref{thm:pattern_stable_generic} proves, in every architecture, that pattern-stable parameters form an open dense set of full Lebesgue measure. Since combinatorial stability implies pattern stability, this is a necessary condition for the unbiased analogue at every depth. We expect our approach to extend to deep networks with biases--the remaining step toward the conjecture is to pass from generic pattern stability to generic combinatorial stability.

We next determine the defining ideal of a class of
shallow ReLU pattern varieties. Alexandr and Mont\'ufar
\cite{alexandr2025constraining} obtained rank constraints on linear combinations of the pattern-indexed maps. Their Problem~35 asks which finite choices of these constraints suffice to produce the ideal's minimal generators. We obtain an explicit finite generating family under a combinatorial condition on the pattern set. The two-block case for shallow networks is currently being studied by \cite{alexandr2026defining}.

For a nonempty set $\mathcal A\subseteq\{0,1\}^{n_1}$,
the \emph{wall graph} $\Gamma_{\mathcal A}$ has vertex
set $\mathcal A$ and an edge labeled $k$ connecting two patterns whenever they differ only in coordinate $k$. Write $N=|\mathcal A|$, let $Y_A$ denote the
$n_2\times n_0$ matrix of coordinate variables indexed
by $A$, and let $J^{\mathcal A}$ be the defining ideal
of $U^{\mathcal A}$.

\begin{restatable}[Generators for shallow networks]
{maintheorem}{shallowgenerators}
\label{thm:shallow_generators}
Let $\mathcal A\subseteq\{0,1\}^{n_1}$ be nonempty.
Suppose $\Gamma_{\mathcal A}$ is connected and has an edge
labeled $k$ for every $k=1,\ldots,n_1$. 
For each $k$, choose a pair
$A^{k,0},A^{k,1}\in\mathcal A$
with $A^{k,1}-A^{k,0}=e_k$, and set
\[
\adjdiff{k}=Y_{A^{k,1}}-Y_{A^{k,0}}.
\]
If a base pattern $A_*\in\mathcal A$ is fixed, then the following hold.
\begin{enumerate}[label=(\roman*)]
\item The ideal $J^{\mathcal A}$ is generated by:
\begin{enumerate}[label=(G\arabic*)]
\item the entries of $Y_{B'}-Y_B-\adjdiff{k}$ for every edge $B'-B=e_k$; 
\item the entries of
$Y_{A_*}-\sum_{k=1}^{n_1}(A_*)_k\adjdiff{k}$;
\item the $2\times2$ minors of $\adjdiff{k}$, for
$k=1,\ldots,n_1$.
\end{enumerate}

\item Let $\rankonevar\subseteq\mathbb C^{n_2\times n_0}$
be the variety of matrices of rank at most one. The linear map
\[
\begin{aligned}
\slopeembed{\mathcal A}:
(\mathbb C^{n_2\times n_0})^{n_1}
&\longrightarrow (\mathbb C^{n_2\times n_0})^N,\\
(V_1,\ldots,V_{n_1})
&\longmapsto
\left(\sum_{k=1}^{n_1}A_kV_k\right)_{A\in\mathcal A}
\end{aligned}
\]
is injective and restricts to an isomorphism
$\rankonevar^{n_1}\cong U^{\mathcal A}$.
Moreover, $J^{\mathcal A}$ is prime and
\[
\dim U^{\mathcal A}=n_1(n_0+n_2-1).
\]
\end{enumerate}
\end{restatable}

The wall-graph hypotheses hold for $\mathcal A=\mathcal{A}_{\theta}$
whenever $\theta$ is shallow and its first-layer rows
are nonzero and pairwise linearly independent
(\Cref{lem:cube}). Consequently,
\Cref{thm:shallow_generators} gives a resolution
of \cite[~Problem~35]{alexandr2025constraining} for this class of canonical pattern sets. The theorem determines the entire defining ideal and identifies the variety as a product of determinantal varieties.
Every generator in the theorem arises from the rank
constraints of \cite[eq.~(5)]{alexandr2025constraining}. 

The product description in \Cref{thm:shallow_generators} also shows why the third
question matters. For a fixed shallow architecture,
all pattern sets satisfying the wall-graph hypotheses
give abstractly isomorphic varieties, even though the
pattern sets themselves need not be equivalent.
We therefore compare the labeled canonical complex,
its canonical pattern set, and the associated variety
under several distinct notions of equivalence.
We write $\theta\sim\theta'$ when their canonical
pattern sets agree after permutations of the hidden
neurons within each layer
(\Cref{def:patternequivalentparameter}).
A block permutation permutes the matrix blocks indexed
by activation patterns, without changing entries
within a block. For shallow networks,
$\mathcal P^{\mathrm{nd}}_{n_0,n_1,n_2}$ denotes the
parameters whose first-layer rows are nonzero and
pairwise linearly independent. Combinatorial
equivalence in the following theorem preserves
the ternary labels.

\begin{restatable}[Hierarchy of equivalences]
{maintheorem}{equivalencehierarchy}
\label{thm:hierarchy_parameters}
For $\theta,\theta'\in\mathcal P_{n_0,\ldots,n_m}$, consider the following statements:
\begin{enumerate}[label=(\alph*)]
\item $\theta$ and $\theta'$ are combinatorially equivalent
\item $\mathcal{A}_{\theta}=\mathcal{A}_{\theta'}$
\item $\theta\sim\theta'$
\item $|\mathcal{A}_{\theta}|=|\mathcal{A}_{\theta'}|\eqqcolon N$ and
$U^{\mathcal{A}_{\theta'}}=\Pi(U^{\mathcal{A}_{\theta}})$ for a block permutation
$\Pi$ of $\mathbb C^{n_m\times Nn_0}$
\item $|\mathcal{A}_{\theta}|=|\mathcal{A}_{\theta'}|\eqqcolon N$ and
$U^{\mathcal{A}_{\theta'}}=T(U^{\mathcal{A}_{\theta}})$ for a linear automorphism
$T$ of $\mathbb C^{n_m\times Nn_0}$
\item $U^{\mathcal{A}_{\theta'}}\cong U^{\mathcal{A}_{\theta}}$ as affine varieties
over $\mathbb C$.
\end{enumerate}
Then
\[
\mathrm{(a)}\Longrightarrow\mathrm{(b)}
\Longrightarrow\mathrm{(c)}
\Longrightarrow\mathrm{(d)}
\Longrightarrow\mathrm{(e)}
\Longrightarrow\mathrm{(f)}.
\]
Every reverse implication fails in general. More precisely:
\begin{enumerate}[label=(\roman*)]
\item For shallow networks, (a) and (b) are equivalent,
including at degenerate parameters. The converse
(b)$\Rightarrow$(a) fails for architecture $(3,4,1,1)$.
\item The converse (c)$\Rightarrow$(b) fails on
$\mathcal P^{\mathrm{nd}}_{2,3,n_2}$ for every $n_2\geq1$.
\item On $\mathcal P^{\mathrm{nd}}_{n_0,n_1,n_2}$ with
$n_0,n_2\geq2$, (c) and (d) are equivalent. The converse
(d)$\Rightarrow$(c) fails, for example, on
$\mathcal P^{\mathrm{nd}}_{2,5,1}$.
\item The converse (e)$\Rightarrow$(d) fails on
$\mathcal P^{\mathrm{nd}}_{2,3,n_2}$ for every $n_2\geq2$.
\item The converse (f)$\Rightarrow$(e) fails on
$\mathcal P^{\mathrm{nd}}_{3,3,n_2}$ for every $n_2\geq1$:
abstractly isomorphic varieties may have different
numbers of pattern blocks.
\end{enumerate}
\end{restatable}

\Cref{thm:hierarchy_parameters} gives the implications
between these equivalence relations and identifies
the hypotheses under which particular converses hold.
For shallow networks, the canonical pattern set
determines the labeled combinatorial type.
For nondegenerate shallow networks with input and
output dimensions at least two, equality of pattern
varieties up to block permutation is equivalent to
pattern equivalence. Arbitrary ambient linear
equivalence and abstract affine isomorphism can
identify varieties whose pattern sets are not
equivalent.

\begin{figure}[htbp]
\centering
\resizebox{\textwidth}{!}{%
\begin{tikzpicture}[
 box/.style={draw,rounded corners,align=center,font=\small,text width=6.0cm,minimum height=0.75cm},
 imp/.style={-{Stealth[length=2.2mm]},thick},
 rev/.style={-{Stealth[length=2mm]},dashed},
 lab/.style={font=\footnotesize,align=left,text width=6.9cm}]
\node[font=\bfseries] at (0,0.8) {Pairs of parameters};
\node[box] (a) at (0,0) {(a) Ternary-label-preserving\\combinatorial equivalence};
\node[box] (b) at (0,-1.65) {(b) $\mathcal{A}_{\theta}=\mathcal{A}_{\theta'}$};
\node[box] (c) at (0,-3.3) {(c) $\theta\sim\theta'$};
\node[box] (d) at (0,-4.95) {(d) Equality after a block permutation};
\node[box] (e) at (0,-6.6) {(e) Ambient linear equivalence};
\node[box] (f) at (0,-8.25) {(f) Abstract affine isomorphism};
\foreach \u/\v in {a/b,b/c,c/d,d/e,e/f} {
 \draw[imp] ([xshift=-0.55cm]\u.south)--([xshift=-0.55cm]\v.north);
 \draw[rev] ([xshift=0.55cm]\v.north)--([xshift=0.55cm]\u.south);
}
\node[lab,anchor=west] at (3.35,-0.825) {Reverse holds for shallow networks;\\fails for architecture $(3,4,1,1)$.};
\node[lab,anchor=west] at (3.35,-2.475) {Reverse fails on $\mathcal P^{\mathrm{nd}}_{2,3,n_2}$.};
\node[lab,anchor=west] at (3.35,-4.125) {Reverse holds on $\mathcal P^{\mathrm{nd}}_{n_0,n_1,n_2}$ if $n_0,n_2\geq2$;\\fails on $\mathcal P^{\mathrm{nd}}_{2,5,1}$.};
\node[lab,anchor=west] at (3.35,-5.775) {Reverse fails on $\mathcal P^{\mathrm{nd}}_{2,3,n_2}$ for $n_2\geq2$.};
\node[lab,anchor=west] at (3.35,-7.425) {Reverse fails on $\mathcal P^{\mathrm{nd}}_{3,3,n_2}$:\\six versus eight pattern blocks.};
\end{tikzpicture}%
}

\vspace{1em}
\resizebox{\textwidth}{!}{%
\begin{tikzpicture}[
 box/.style={
   draw,
   rounded corners,
   align=center,
   font=\small,
   text width=4.5cm,
   minimum height=0.8cm
 },
 imp/.style={-{Stealth[length=2.2mm]},thick},
 rev/.style={-{Stealth[length=2mm]},dashed},
 lab/.style={font=\footnotesize,align=center}
]

\node[font=\bfseries] at (3.7,0.85) {A single parameter};

\node[box] (cs) at (0,0)
  {Combinatorially stable};
\node[box] (ps) at (0,-1.8)
  {Pattern-stable};
\node[box] (pi) at (0,-3.6)
  {Pattern-inclusive};
\node[box] (le) at (7.4,-1.8)
  {Interior point of its\\pattern-equivalence class};

\draw[imp]
  ([xshift=0.5cm]cs.south) --
  ([xshift=0.5cm]ps.north);
\draw[imp]
  ([xshift=0.5cm]ps.south) --
  ([xshift=0.5cm]pi.north);

\draw[rev]
  ([xshift=-0.5cm]ps.north) --
  ([xshift=-0.5cm]cs.south);
\draw[rev]
  ([xshift=-0.5cm]pi.north) --
  ([xshift=-0.5cm]ps.south);

\draw[imp]
  ([yshift=0.18cm]ps.east) --
  ([yshift=0.18cm]le.west);
\draw[rev]
  ([yshift=-0.18cm]le.west) --
  ([yshift=-0.18cm]ps.east);

\node[
  lab,
  anchor=east,
  text width=4.2cm
] at (-0.85,-0.9)
  {Reverse: shallow, yes;\\deep, no.};

\node[
  lab,
  anchor=east,
  text width=4.2cm
] at (-0.85,-2.7)
  {Reverse fails for architecture $(2,2,1)$.};

\node[
  lab,
  text width=4.2cm
] at (3.7,-2.8)
  {Reverse: shallow, yes;\\deep, open.};

\end{tikzpicture}%
}
\caption{The hierarchy of \Cref{thm:hierarchy_parameters}, the local implications \eqref{eq:hierarchy_local}, and \Cref{prop:local_pattern_equivalence}. Solid arrows are implications valid for every architecture; dashed arrows are converses, whose status is stated beside them. Conditions (d) and (e) require a common number of blocks. Pattern stability is generic at every depth by \Cref{prop:pattern_stable_generic}.}
\label{fig:hierarchy}
\end{figure}
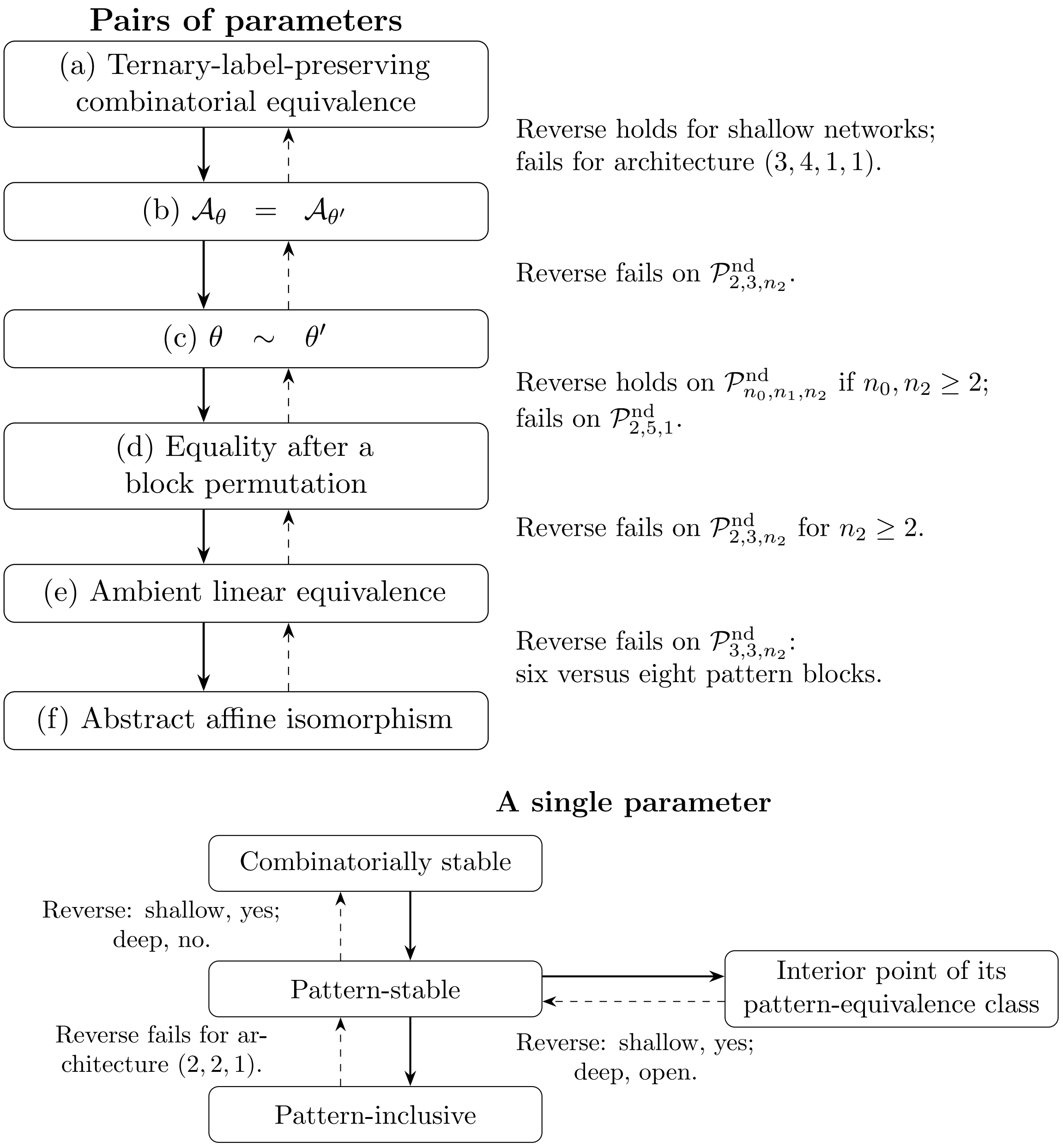
Figure~\ref{fig:hierarchy} summarizes this hierarchy
and the corresponding local stability implications.
Here, pattern inclusiveness means that every pattern
in $\mathcal{A}_{\theta}$ persists throughout a common neighborhood
of $\theta$. The two panels distinguish equivalence
relations between pairs of parameters from local
properties of a single parameter.

For two-dimensional input, we further classify the
canonical pattern sets of nondegenerate shallow
networks up to equivalence by dihedral orbits
of antipodal sign sequences. Their enumeration
agrees with a classical count \cite{PalmerRobinson1984}. This gives a concrete
description of the pattern-equivalence classes
appearing in \Cref{thm:hierarchy_parameters}. We also show that the pattern variety of every shallow network with one-dimensional output is a linear space.

Our final main result gives a functional interpretation of $\dim U^{\mathcal{A}_\theta}$, the dimension of the canonical ReLU pattern variety. We compare this algebraic dimension with the number of independent first-order variations of the realized function detectable by finitely many evaluations. Following the evaluation-map approach of \cite{GLMWFunctionalDimension}, we measure these variations by the ranks of evaluation Jacobians.

For this comparison, call an input $z$ \emph{$\theta$-generic} if every hidden preactivation at $z$ is nonzero, and write
\[
\mathcal{A}_\theta^{\mathrm{str}}
\coloneqq
\{\alpha_\theta(z):z\text{ is }\theta\text{-generic}\}.
\]
For a finite ordered set $Z=(z_1,\ldots,z_k)$ of such inputs, let
\[
E_Z(\theta')
\coloneqq
[\,f_{\theta'}(z_1)\mid\cdots\mid f_{\theta'}(z_k)\,],
\qquad
\dimfunzero(\theta)
\coloneqq
\max_Z\operatorname{rank}JE_Z(\theta),
\]
where derivatives are taken with respect to the unbiased
weights. The empty sample set is allowed and has rank zero. This defines an unbiased functional dimension on strict activation regions. In each differential
$d\varphi^{\mathcal A}(\theta)$ below, the pattern set
$\mathcal A$ is held fixed as the weights vary.

\begin{restatable}[Functional dimension and the canonical pattern variety]
{maintheorem}{functionaldimension}
\label{thm:functional_dimension}
Let $\theta\in\mathcal P_{n_0,\ldots,n_m}$.
\begin{enumerate}[label=(\roman*)]
\item At every depth,
\[
\dimfunzero(\theta)
=\operatorname{rank}d\varphi^{\mathcal{A}_\theta^{\mathrm{str}}}(\theta)
\leq\operatorname{rank}d\varphi^{\mathcal{A}_\theta}(\theta)
\leq\dim U^{\mathcal{A}_\theta}.
\]
For $\mathcal{A}_\theta^{\mathrm{str}}=\varnothing$, the first differential
is the zero map.

\item If the interior of every maximal-dimensional cell of
$\mathcal C(\theta)$ contains a $\theta$-generic point, then
\[
\dimfunzero(\theta)
=\operatorname{rank}d\varphi^{\mathcal{A}_\theta}(\theta).
\]
For shallow networks, this hypothesis holds if and only if
$W_1$ has no zero row.

\item Suppose $\theta\in\mathcal P^{\mathrm{nd}}_{n_0,n_1,n_2}$.
Write $w_k^\top$ for row $k$ of $W_1$ and
$c_k\in\mathbb R^{n_2}$ for column $k$ of $W_2$. Then
\[
\dimfunzero(\theta)
=n_1(n_0+n_2-1)
-(n_0-1)\cdot\big|\{k:c_k=0\}\big|.
\]
Consequently, $\dimfunzero(\theta)=\dim U^{\mathcal{A}_\theta}$
if and only if $n_0=1$ or $W_2$ has no zero column.

\item For every shallow architecture $(n_0,n_1,n_2)$,
there is an open dense subset of
$\mathcal P_{n_0,n_1,n_2}$ of full Lebesgue measure on which
\[
\dimfunzero(\theta)
=\operatorname{rank}d\varphi^{\mathcal{A}_\theta}(\theta)
=\dim U^{\mathcal{A}_\theta}.
\]
\end{enumerate}
\end{restatable}

\Cref{thm:functional_dimension} gives the algebraic dimension a direct interpretation: generically for shallow networks, it counts the independent first-order variations of the realized function detectable by finite evaluations. The two inequalities distinguish separate obstructions to this interpretation at an arbitrary parameter. Strict sampling may omit canonical patterns, while the polynomial parametrization may have less than its generic rank. We give explicit examples of both phenomena.

The sampling convention is substantive. Our definition
differs from that of \cite{GLMWFunctionalDimension} in its restriction to nonzero hidden preactivations, as well as its unbiased parameter space and linear output layer. Under the hypothesis of
\Cref{thm:functional_dimension}~(ii), allowing all parametrically smooth sample points \cite[Definition~3.9]{GLMWFunctionalDimension} gives the same value for our architecture (\Cref{prop:functional_smooth_sampling}).
This applies, in particular, to every shallow parameter
with no zero first-layer row.

Our results connect with several other approaches to ReLU networks. Tropical geometry describes their piecewise-linear structure through polyhedral constructions \cite{zhang2018tropical,brandenburg2024real}, while Fu's toric formulation gives intersection-number compatibility conditions for rational-weight shallow networks \cite{fu2025toric}. Our bending relations express the corresponding equalities of matrix jumps across activation hyperplanes, with arbitrary real weights. The study of pattern varieties also forms part of the broader algebraic-geometric study of parametrized function spaces \cite{marchetti2025neuroalgebraic}. Our functional-dimension comparison connects this algebraic viewpoint with the evaluation-map approach of
\cite{GLMWFunctionalDimension}: for generic shallow parameters, the dimension of the canonical pattern variety equals the number of independent first-order function variations detected by our admissible evaluations.

Work on invariant embeddings, reconstruction, and hidden
symmetries addresses the complementary question of recovering parameters from realized functions
\cite{stock2023embedding,rolnick2020reverse,grigsby2023hidden}. Pattern equivalence concerns the canonical activation sets and does not require equality of the realized functions, whereas functional dimension measures the rank of their infinitesimal variation under parameter perturbations.

The constructions motivating our questions also make sense for networks with biases. We expect biased analogues of our stability and generator results to follow by related methods, but their precise formulation requires a separate treatment of affine activation arrangements and their degeneracies. Extending the functional-dimension comparison would also require including bias derivatives in the evaluation Jacobians and the constant terms of the affine pieces in the pattern parametrization. We leave these extensions for future work.

\Cref{sec:prelim} introduces the objects and conventions. \Cref{sec:pattern_stability} studies perturbations and shallow stability, and \Cref{sec:generators} determines the shallow generators. \Cref{sec:permutations} introduces pattern equivalence and the clock classification. \Cref{sec:hierarchy} establishes the hierarchy of equivalences and generic pattern stability, and states the remaining local-stability questions. \Cref{sec:functional_dimension} compares functional dimension with the differential rank and dimension of the canonical pattern variety, and gives examples of strict inequalities. The appendix introduces clocks and provides the relevant proofs.

\subsection*{Acknowledgements} We are grateful to Yulia Alexandr, Elisenda Grigsby, Moritz Grillo, Kathryn Lindsey, and Marissa Masden for insightful discussions. 

{\footnotesize BF was supported by the National Science Foundation Graduate Research Fellowship Program under Grant No. 2137424.  Any opinions, findings, and conclusions or recommendations expressed in this material are those of the authors and do not necessarily reflect the views of the National Science Foundation. Support was also provided by the Graduate School and the Office of the Vice Chancellor for Research at the University of Wisconsin-Madison with funding from the Wisconsin Alumni Research Foundation, as well as by National Science Foundation Grant No. 2510307. YF was partially supported by the Institute for the Foundations of Machine Learning (IFML). CW was supported partly by a Simons Collaboration Grant. PF, EM and WW were supported by undergraduate research fellowships to conduct research under YF's mentorship with the Boston College Math and Machine Learning Lab.
}

\subsection*{Statement of AI Use}

Publicly available AI models from OpenAI and Anthropic, such as ChatGPT 6 and Claude 5.5, were used for proofreading and revising the text.

\setcounter{tocdepth}{2}
\tableofcontents

\section{Preliminaries}\label{sec:prelim}

When possible, we follow the notation and terminology of
\cite{alexandr2025constraining,fu2025toric}. We also recall the notion of combinatorial stability introduced in \cite{GLMWFunctionalDimension}.

\subsection{Feedforward ReLU neural networks}\label{subsec:feedfoward ReLU}

The \demph{rectified linear unit (ReLU)} is the map
$\ReLU:\mathbb{R}\to\mathbb{R}$ defined by
\[
\ReLU(t):=\max\{0,t\}.
\]
For any positive integer $n$, we denote by
$\act:\mathbb{R}^n\to\mathbb{R}^n$ the component-wise
ReLU activation function
\[
\act(x_1,\ldots,x_n)
=
(\ReLU(x_1),\ldots,\ReLU(x_n)).
\]

Fix an architecture $(n_0,n_1,\ldots,n_m)$ with $m\geq2$, where
$n_0$ is the input dimension, $n_m$ is the output
dimension, and $n_1,\ldots,n_{m-1}$ are the widths
of the hidden layers.
An \demph{unbiased feedforward ReLU neural network}
with this architecture is a composition
\[
f_\theta
=
W_m\circ\act\circ W_{m-1}\circ\cdots
\circ\act\circ W_1,
\]
where
\[
W_\ell:\mathbb{R}^{n_{\ell-1}}
\longrightarrow\mathbb{R}^{n_\ell}
\]
is a linear map represented by a matrix
$W_\ell\in\mathbb{R}^{n_\ell\times n_{\ell-1}}$.
The \demph{parameter space} associated with the
architecture is
\[
\mathcal{P}_{n_0,\ldots,n_m}
:=
\prod_{\ell=1}^{m}
\mathbb{R}^{n_\ell\times n_{\ell-1}},
\]
and a parameter
\[
\theta=(W_1,\ldots,W_m)
\in\mathcal{P}_{n_0,\ldots,n_m}
\]
records all the weight matrices of the network.
There is no ReLU activation after the output layer.
A network with exactly one hidden layer (equivalently,
$m=2$) is called \demph{shallow}, while a network
with at least two hidden layers is called \demph{deep}.

For each input $x\in\mathbb{R}^{n_0}$, define
recursively
\[
z_\theta^{(0)}(x):=x,
\]
and, for $\ell=1,\ldots,m-1$,
\[
h_\theta^{(\ell)}(x)
:=
W_\ell z_\theta^{(\ell-1)}(x),
\qquad
z_\theta^{(\ell)}(x)
:=
\act\bigl(h_\theta^{(\ell)}(x)\bigr).
\]
Here $h_\theta^{(\ell)}(x)$ and
$z_\theta^{(\ell)}(x)$ denote the vectors of
preactivation and postactivation values,
respectively, in the $\ell$-th hidden layer.
The output function is therefore
\[
f_\theta(x)
=
W_m z_\theta^{(m-1)}(x).
\]

Since the network has no biases, $f_\theta$ and
all hidden-layer preactivation and postactivation
functions are positively homogeneous of degree one:
for every $\lambda\geq0$,
\[
f_\theta(\lambda x)=\lambda f_\theta(x).
\]
Moreover, $f_\theta$ is continuous and finitely
piecewise linear.

\subsection{Canonical polyhedral complexes and
activation patterns}\label{subsec:canonical polyhedral complex}
Let $\theta\in\mathcal{P}_{n_0,\ldots,n_m}$ be a parameter of an unbiased feedforward ReLU network. Following Grigsby and Lindsey
\cite[Definition~6.5]{grigsby2022transversality}, we define the
\demph{canonical polyhedral complex} $\mathcal{C}(\theta)$ by
successively subdividing $\mathbb{R}^{n_0}$ according to the
preactivations of the hidden neurons.

The first hidden layer induces a central hyperplane arrangement. For each $i\in\{1,\ldots,n_1\}$ such that the $i$-th row of $W_1$ is nonzero, the corresponding hyperplane is
\[
H_i=\{x\in\mathbb{R}^{n_0}:(W_1x)_i=0\}.
\]
The closures of the chambers of this arrangement, together with
all their faces, form the first-layer decomposition.

For $\ell=1,\ldots,m-1$, let $C$ be a cell of the decomposition obtained after the first $\ell-1$ hidden layers (with $C=\mathbb{R}^{n_0}$ for the zeroth hidden layer). The preactivations $h_{\theta,1}^{(\ell)},\ldots,h_{\theta,n_\ell}^{(\ell)}$ are linear on $C$. Their zero sets subdivide $C$, with any preactivation that vanishes identically on $C$ contributing no cut. The resulting closed cells and all their faces form the decomposition after layer $\ell$. The decomposition obtained after the last hidden layer is
$\mathcal{C}(\theta)$.

By construction, every hidden-layer preactivation and postactivation
is linear on each cell of $\mathcal{C}(\theta)$, and hence so is
$f_\theta$. On the relative interior of each cell, every hidden
preactivation is either strictly positive, strictly negative, or
identically zero.
Since the network is unbiased, the cells of $\mathcal{C}(\theta)$ are polyhedral cones. Its full-dimensional cones cover $\mathbb{R}^{n_0}$.  

We now define the activation patterns that will be used throughout this paper.
\begin{definition}[Activation pattern]
\label{def:activation_pattern}
Let $\theta\in\mathcal P_{n_0,\ldots,n_m}$ and $x\in\mathbb R^{n_0}$. For each hidden layer $\ell=1,\ldots,m-1$, define
\[
\alpha_\theta^{(\ell)}(x)
=
\bigl(
\alpha_{\theta,1}^{(\ell)}(x),
\ldots,
\alpha_{\theta,n_\ell}^{(\ell)}(x)
\bigr)
\in\{0,1\}^{n_\ell},
\]
where
\[
\alpha_{\theta,i}^{(\ell)}(x)
=
\begin{cases}
1, & h_{\theta,i}^{(\ell)}(x)>0,\\
0, & h_{\theta,i}^{(\ell)}(x)\leq0.
\end{cases}
\]
That is, a neuron is active if its preactivation is strictly positive and inactive otherwise.

The \demph{activation pattern} of $x$ associated with $\theta$ is
\[
\alpha_\theta(x)
=
\bigl(
\alpha_\theta^{(1)}(x),
\ldots,
\alpha_\theta^{(m-1)}(x)
\bigr)
\in
\prod_{\ell=1}^{m-1}\{0,1\}^{n_\ell}.
\]
\end{definition}

Thus, $\alpha_\theta$ is constant on the relative interior of each cell. For a full-dimensional cell $C$, write $A(C)$ for this common value on $\operatorname{int}C$.

\begin{definition}[Canonical set of activation patterns]\label{def:canonical_set_of_activation_patterns}
Let $\theta\in\mathcal{P}_{n_0,\ldots,n_m}$.
The \demph{canonical set of activation patterns}
associated with $\theta$ is
\[
\mathcal{A}_\theta
:=
\left\{
A(C):
C\in\mathcal{C}(\theta),\ \dim C=n_0
\right\}.
\]
Thus $\mathcal{A}_\theta$ consists of the activation patterns
on the interiors of the maximal-dimensional cells
of $\mathcal{C}(\theta)$. We call its elements the
\demph{realizable activation patterns} of $\theta$,
and say that $\theta$ \emph{realizes} $A$ if $A\in \mathcal{A}_\theta$. 
\end{definition}

The set $\mathcal{A}_\theta$ is finite. For a fixed activation
pattern $A$, its \demph{activation region} is
\[
R_A(\theta):=\{x\in\mathbb{R}^{n_0}:\alpha_\theta(x)=A\}.
\]
Once the activation statuses of the preceding
neurons are fixed, each preactivation is a linear
function of the input. Hence $R_A(\theta)$ is
defined by finitely many strict or weak linear
inequalities and is therefore convex.

Distinct maximal-dimensional cells of $\mathcal{C}(\theta)$ have distinct activation patterns. Since these cells are obtained by successive subdivisions of $\mathbb{R}^{n_0}$, two distinct cells must first be separated at some hidden layer. At that layer, some preactivation is positive on the interior of one cell and negative on the interior of the other, so their activation patterns differ. 
Thus $C\mapsto A(C)$ gives a one-to-one correspondence between the
maximal-dimensional cells of $\mathcal{C}(\theta)$
and the elements of $\mathcal{A}_\theta$.
Patterns that occur exclusively on
lower-dimensional cells are excluded from $\mathcal{A}_\theta$. By definition, a neuron whose
preactivation vanishes identically on a
maximal-dimensional cell is inactive there.

For any activation pattern
\[
A=(A^{(1)},\ldots,A^{(m-1)}),
\qquad A^{(\ell)}\in\{0,1\}^{n_\ell},
\]
let
\[
D_\ell(A):=\operatorname{diag}(A^{(\ell)}),
\]
and define
\[
M_A(\theta)
:=
W_mD_{m-1}(A)W_{m-1}\cdots D_1(A)W_1.
\]
Then
\[
f_\theta(x)=M_A(\theta)x
\qquad\text{for every }x\in R_A(\theta).
\]
By continuity, this equality also holds on the
whole closed cell $C$ whenever $A(C)=A$.

\subsection{ReLU pattern varieties and pattern ideals}\label{subsec:ReLU pattern varieties}

We recall the ReLU pattern varieties and their defining ideals introduced in \cite[Section 3.3]{alexandr2025constraining}, specialized to unbiased networks.

Fix an architecture $(n_0,\ldots,n_m)$ and a finite set
$\mathcal A=\{A_1,\ldots,A_N\}$ of distinct binary activation
patterns. Consider the polynomial parametrization
\[
\begin{aligned}
\varphi^{\mathcal A}:
\mathcal P_{n_0,\ldots,n_m}
&\longrightarrow \mathbb R^{n_m\times Nn_0},\\
\theta
&\longmapsto
[M_{A_1}(\theta)\mid\cdots\mid M_{A_N}(\theta)].
\end{aligned}
\]
Here $\mathcal A$ is fixed while $\theta$ varies over
the entire parameter space, including parameters that
do not realize every pattern in $\mathcal A$.
In \cite{alexandr2025constraining}, the parametrization also includes the constant terms of the affine pieces, which are absent in the unbiased setting.

\begin{definition}[ReLU pattern variety]\label{rpv}
The \demph{ReLU pattern variety} associated with
$\mathcal A$ is the Zariski closure
\[
U^{\mathcal A}
:=
\overline{
\varphi^{\mathcal A}(\mathcal P_{n_0,\ldots,n_m})
}^{\,\mathrm{Zar}}
\subseteq \mathbb C^{n_m\times Nn_0}.
\]
\end{definition}

In this paper, we are particularly interested in
$\mathcal A=\mathcal{A}_{\theta_0}$, where $\mathcal{A}_{\theta_0}$ is
the canonical set of activation patterns associated
with a fixed parameter $\theta_0$
(Definition~\ref{def:canonical_set_of_activation_patterns}). Thus $\theta_0$ determines the prescribed patterns, while the parameter $\theta$ in $\varphi^{\mathcal A}$ is free to vary.

\begin{remark}\label{rmk:domain_of_varphi}
When $\mathcal A=\mathcal{A}_{\theta_0}$, it is natural to restrict
$\varphi^{\mathcal{A}_{\theta_0}}$ to a neighborhood of $\theta_0$
on which all these patterns remain realizable, whenever such a neighborhood exists.
This motivates the notions of combinatorial stability, pattern-stability, and pattern-inclusiveness (Definitions~\ref{def:combinatorial_stable},
\ref{def:pattern_stable_neighborhood}, and
\ref{def:patterninclusivenbhd}).
Combinatorial stability and pattern-stability require
$\mathcal{A}_{\theta}=\mathcal{A}_{\theta_0}$ throughout a sufficiently small
neighborhood, while pattern-inclusiveness requires only
$\mathcal{A}_{\theta_0}\subseteq \mathcal{A}_{\theta}$.

Restricting the domain to any nonempty Euclidean-open
subset $\mathcal{N}_{\theta_0}$ of the parameter space gives the same Zariski closure:
\[
\overline{\varphi^{\mathcal{A}_{\theta_0}}(\mathcal{N}_{\theta_0})}^{\,\mathrm{Zar}}
=
U^{\mathcal{A}_{\theta_0}}.
\]
Indeed, the parameter space is isomorphic to
$\mathbb R^D$, where $D=\sum_{\ell=1}^m n_\ell n_{\ell-1}$,
and $\mathcal{N}_{\theta_0}$ is Zariski dense in $\mathbb C^D$. If a polynomial $F$ vanishes on
$\varphi^{\mathcal{A}_{\theta_0}}(\mathcal N_{\theta_0})$,
then $F\circ\varphi^{\mathcal{A}_{\theta_0}}$ vanishes on
$\mathcal N_{\theta_0}$ and, by Zariski density,
on the entire parameter space. Thus restricting the parameters to $\mathcal N_{\theta_0}$ gives the same polynomial relations among the matrices $M_A(\theta)$, $A\in \mathcal{A}_{\theta_0}$.
\end{remark}

For each $j=1,\ldots,N$, let $Y_j=(y_{ab}^{(j)})$
be an $n_m\times n_0$ matrix of indeterminates, and write
\[
\mathbf Y=[Y_1\mid\cdots\mid Y_N].
\]
We denote by $\mathbb C[\mathbf Y]$ the polynomial ring
in the entries of these matrices and write $Y_A:=Y_j$
when $A=A_j$. Reordering the patterns in $\mathcal A$
amounts to permuting the coordinate blocks.

\begin{definition}[ReLU pattern ideal]
The \demph{ReLU pattern ideal} associated with $\mathcal A$
is the defining ideal
\[
J^{\mathcal A}:=I(U^{\mathcal A})
\subseteq\mathbb C[\mathbf Y].
\]
Equivalently,
\[
J^{\mathcal A}
=
\left\{
F\in\mathbb C[\mathbf Y]:
F\bigl(\varphi^{\mathcal A}(\theta)\bigr)=0
\text{ for every }\theta\in\mathcal P_{n_0,\ldots,n_m}
\right\}.
\]
\end{definition}

Let $\mathbf W$ denote the collection of weight matrix
entries, regarded as indeterminates. Since $\mathbb R^D$
is Zariski dense in $\mathbb C^D$, we also have
\[
J^{\mathcal A}=\ker(\varphi^{\mathcal A})^*,
\]
where
\[
\begin{aligned}
(\varphi^{\mathcal A})^*:
\mathbb C[\mathbf Y]
&\longrightarrow \mathbb C[\mathbf W],\\
y_{ab}^{(j)}
&\longmapsto \bigl(M_{A_j}(\mathbf W)\bigr)_{ab},
\end{aligned}
\]
for $1\leq j\leq N$, $1\leq a\leq n_m$, and
$1\leq b\leq n_0$.
Because $\mathbb C[\mathbf W]$ is an integral domain,
$J^{\mathcal A}$ is prime and consequently
$U^{\mathcal A}$ is irreducible.
The ideal $J^{\mathcal A}$ records the polynomial
relations among the matrices associated with the
prescribed activation patterns. Determining its
generators is the implicitization problem for
$\varphi^{\mathcal A}$.

\subsection{Combinatorial stability}
\label{subsec:combinatorial_stability}

We adapt the notion of combinatorial stability from
\cite[Definition~5.3]{GLMWFunctionalDimension} to our
setting of unbiased networks with a linear output layer.
All neighborhoods below are taken in the weight space
$\mathcal P_{n_0,\ldots,n_m}$, so perturbations preserve
the absence of biases.

For each input $x\in\mathbb R^{n_0}$, define its
\demph{ternary activation labeling} by
\[
s_\theta(x)
=
\bigl(s_\theta^{(1)}(x),\ldots,
s_\theta^{(m-1)}(x)\bigr)
\in\prod_{\ell=1}^{m-1}\{-1,0,1\}^{n_\ell},
\]
where
\[
s_{\theta,i}^{(\ell)}(x)
=
\operatorname{sgn}\bigl(h_{\theta,i}^{(\ell)}(x)\bigr),
\qquad
1\leq\ell\leq m-1,\quad 1\leq i\leq n_\ell,
\]
and $\operatorname{sgn}(0)=0$.
The ternary labeling distinguishes negative and zero
preactivations, both of which correspond to inactive
neurons in the binary pattern.
It is constant on the relative interior of each cell
of $\mathcal C(\theta)$; we call this value the
\demph{ternary labeling of the cell}.

Two polyhedral complexes are
\demph{combinatorially equivalent} if there is a
dimension-preserving bijection $\psi$ between their
cells such that $F$ is a face of $G$ if and only if
$\psi(F)$ is a face of $\psi(G)$
\cite[Definition~5.1]{GLMWFunctionalDimension}.
For complexes with ternary labelings, such an
equivalence \demph{respects the ternary labelings}
if corresponding cells have the same labeling
\cite[Definition~5.2]{GLMWFunctionalDimension}.

\begin{definition}[Combinatorially stable parameter
    {\cite[Definition~5.3]{GLMWFunctionalDimension}}]
\label{def:combinatorial_stable}
A parameter $\theta\in\mathcal P_{n_0,\ldots,n_m}$
is \demph{combinatorially stable} if there exist
an open neighborhood $\mathcal U_\theta$ of $\theta$
in $\mathcal P_{n_0,\ldots,n_m}$ and a continuous map
\[
\Phi:\mathcal U_\theta\times\mathbb R^{n_0}
\longrightarrow\mathbb R^{n_0}
\]
such that, for every $\theta'\in\mathcal U_\theta$,
the map $\Phi_{\theta'}:=\Phi(\theta',\cdot)$ is a
homeomorphism of $\mathbb R^{n_0}$ that maps cells
of $\mathcal C(\theta)$ onto cells of
$\mathcal C(\theta')$ and induces a combinatorial
equivalence respecting their ternary labelings.
\end{definition}

Since the combinatorial equivalence preserves the ternary labelings, it also preserves the binary activation patterns: a ternary entry of $1$ gives a binary entry of $1$, while entries of $0$ and $-1$ both give a binary entry of $0$. Corresponding maximal-dimensional cells therefore have
the same activation pattern, and hence
\[
\mathcal{A}_{\theta'}=\mathcal{A}_\theta
\qquad\text{for every }\theta'\in\mathcal U_\theta.
\]
The pattern-stability and pattern-inclusiveness
properties introduced in the next section concern
the preservation of these sets of binary patterns.
Combinatorial stability additionally requires the
preservation of the face relations and ternary
labelings of all cells through the continuous family
of homeomorphisms above.

\section{Activation patterns under perturbations}\label{sec:pattern_stability}

\subsection{Pattern stability and its obstructions}\label{subsec:pattern_stability}
\begin{definition}[Pattern-stable neighborhood]\label{def:pattern_stable_neighborhood}
A parameter $\theta\in\mathcal{P}_{n_0,\ldots,n_m}$ is
\emph{pattern-stable} if there exists an open neighborhood $\mathcal N_\theta$
of $\theta$ such that, for every $\theta'\in \mathcal N_\theta$,
\[
\mathcal{A}_{\theta'}=\mathcal{A}_\theta.
\]
\end{definition}

\begin{proposition}\label{prop:pattern_stable_parameter}

A parameter $\theta \in \mathcal{P}_{n_0,\ldots,n_m}$ is not pattern-stable
in either of the following cases:
\begin{enumerate}[label=(\roman*)]
    \item there is a zero row in the first hidden-layer
    weight matrix $W_1$, or

    \item $n_0\geq 2$, and there exist two opposite
    neurons in the first hidden layer, i.e.,
    two distinct rows $w_i^\top,w_j^\top$ of
    $W_1$ such that
    \[
    w_j=-\lambda w_i
    \]
    for some $\lambda>0$ and $w_i\neq 0$.
\end{enumerate}
\end{proposition}

\begin{proof}
We show that in each case, there exist arbitrarily
small perturbations $\theta'$ of $\theta$ such that
\[
\mathcal{A}_{\theta'}\neq \mathcal{A}_\theta.
\]

\emph{Case (i): $W_1$ contains a zero row.}
Suppose that the $k$-th row of $W_1$ is zero, i.e.,
$w_k=0$. Then
\[
w_k^\top x=0
\qquad
\text{for every }x\in\mathbb{R}^{n_0}.
\]
Thus, the $k$-th neuron is inactive on every
maximal-dimensional cell of the canonical polyhedral complex.

Choose any nonzero vector $v\in\mathbb{R}^{n_0}$.
For any $\varepsilon>0$, perturb only the $k$-th row by
setting
\[
w_k'=\varepsilon v,
\]
while leaving all other parameters unchanged.
Then
\[
(w_k')^\top x=\varepsilon v^\top x>0
\]
on the nonempty open half-space
\[
U=\{x\in\mathbb{R}^{n_0}:v^\top x>0\}.
\]
Hence the $k$-th neuron is active on $U$ under
the perturbed parameter $\theta'$.

Since $U$ is nonempty and open, it meets the interior of a maximal-dimensional cell of $\mathcal C(\theta')$.
Thus $\mathcal{A}_{\theta'}$ contains a pattern $A$ with $A^{(1)}_k=1$, whereas every pattern $A$ in $\mathcal{A}_\theta$
has $A^{(1)}_k=0$. Consequently,
\[
\mathcal{A}_{\theta'}\neq \mathcal{A}_\theta.
\]
Since $\varepsilon>0$ can be arbitrarily small,
$\theta$ is not pattern-stable.

\emph{Case (ii): There are two opposite neurons in the first hidden layer.}
Suppose that $n_0\geq 2$ and there exist
distinct indices $i,j$ such that
\[
w_j=-\lambda w_i
\]
for some $\lambda>0$ and $w_i\neq 0$.

For every $x\in\mathbb{R}^{n_0}$,
\[
w_j^\top x=-\lambda w_i^\top x.
\]
Thus the $i$-th and $j$-th neurons cannot be
simultaneously active. In particular, no pattern
$A$ in $\mathcal{A}_\theta$ has $A^{(1)}_i=A^{(1)}_j=1$.

Since $n_0\geq 2$, we can choose a nonzero vector
$v\in\mathbb{R}^{n_0}$ satisfying
\[
v^\top w_i=0.
\]
For any $\varepsilon>0$, perturb only the $j$-th
row by setting
\[
w_j'=w_j+\varepsilon v
=-\lambda w_i+\varepsilon v,
\]
while leaving all other parameters unchanged.

Choose
\[
x=v+t w_i,
\]
where
\[
0<t<
\frac{\varepsilon\|v\|^2}
{\lambda\|w_i\|^2}.
\]
Since $v^\top w_i=0$, we have
\[
w_i^\top x=w_i^\top (v+tw_i)=t\|w_i\|^2>0.
\]
Moreover,
\begin{align*}
(w_j')^\top x
&=(-\lambda w_i+\varepsilon v)^\top
(v+t w_i)\\
&=-\lambda t\|w_i\|^2
+\varepsilon\|v\|^2\\
&>0.
\end{align*}
Therefore the $i$-th and $j$-th neurons are
simultaneously active at $x$ under $\theta'$.

Since both inequalities are strict, by continuity, there
exists a nonempty open neighborhood $U$ of $x$ such that
\[
w_i^\top y>0
\qquad\text{and}\qquad
(w_j')^\top y>0
\]
for every $y\in U$.

Recall that the canonical polyhedral complex $\mathcal C(\theta')$ consists of finitely many polyhedral cells,
and the union of the interiors of its maximal-dimensional
cells is dense in $\mathbb{R}^{n_0}$. Hence $U$ contains
a point in the interior of some maximal-dimensional
cell $C'$ of $\mathcal C(\theta')$.
Since $\operatorname{int}C'$ meets $U$, and the activation pattern
is constant on $\operatorname{int}C'$, both the $i$-th and $j$-th neurons
are active on $\operatorname{int}C'$.

Thus $\mathcal{A}_{\theta'}$ contains an activation pattern $A$ with $A^{(1)}_i=A^{(1)}_j=1$, whereas no such
pattern belongs to $\mathcal{A}_\theta$. Therefore
\[
\mathcal{A}_{\theta'}\neq \mathcal{A}_\theta.
\]
Since $\varepsilon>0$ can be chosen arbitrarily
small, $\theta$ is not pattern-stable.
\end{proof}
\subsection{Pattern inclusiveness}\label{subsec:pattern_inclusiveness}
\begin{definition}[Pattern-inclusive parameter]\label{def:patterninclusivenbhd}
    Let $\theta \in \mathcal{P}_{n_0,\ldots,n_m}$, and let
$\mathcal{A}_\theta=\{A_1,\ldots,A_N\}$
be its canonical set of activation patterns. A parameter $\theta$ is \emph{pattern-inclusive} if there exists an open neighborhood $\mathcal{N}_\theta$ of $\theta$ such that, for every $\theta' \in \mathcal{N}_\theta$,
    $$
    \mathcal{A}_\theta \subseteq \mathcal{A}_{\theta'}.
    $$
    We call $\mathcal{N}_\theta$ a \emph{pattern-inclusive neighborhood of $\theta$}.
\end{definition}

\begin{proposition}\label{prop:pattern_inclusive_nbhd_exist}
Let $\theta\in\mathcal{P}_{n_0,n_1,n_2}$ be a parameter of an
unbiased shallow ReLU network
\[
f_\theta(x)=W_2\circ \varsigma \circ W_1(x),
\]
where
\[
W_1=
\begin{pmatrix}
w_1^\top\\
\vdots\\
w_{n_1}^\top
\end{pmatrix},
\qquad w_i\in\mathbb{R}^{n_0}.
\]
 Then $\theta$ is pattern-inclusive if and only if
either
\begin{enumerate}[label=(\roman*)]
    \item $W_1$ has no zero rows, or
    \item $n_1=1$ and $W_1=0$.
\end{enumerate}
In particular, if $n_1\geq 2$, then $\theta$ is
pattern-inclusive if and only if $W_1$ has no zero rows.
\end{proposition}

\begin{proof}
Since the network is shallow, its activation patterns are
determined entirely by the first-layer weight matrix $W_1$.
We consider the following cases.

\emph{Case (i). $W_1$ has no zero rows.}
Let
\[
\mathcal{A}_\theta=\{A_1,\ldots,A_N\}
\]
be the set of activation patterns associated with the
maximal-dimensional cones of the canonical polyhedral complex.

For each $j=1,\ldots,N$, let $C_j$ be a
maximal-dimensional cone with activation pattern $A_j$ on its interior.
Choose
\[
x_j\in\operatorname{int}(C_j),
\qquad \|x_j\|=1.
\]
Since $x_j$ lies in the interior of $C_j$, it lies on none
of the hyperplanes $H_i=\{x:w_i^\top x=0\}$. Thus
\[
w_i^\top x_j\neq 0
\qquad
\text{for every }i=1,\ldots,n_1.
\]
Define
\[
\delta_j
:=
\min_{1\leq i\leq n_1}|w_i^\top x_j|>0
\]
and choose
\[
\varepsilon
:=
\frac12\min_{1\leq j\leq N}\delta_j>0.
\]

Let $\mathcal{N}_\theta$ be an open neighborhood of
$\theta$ such that, for every
$\theta'\in\mathcal{N}_\theta$, the rows
$(w_i')^\top$ of the perturbed first-layer weight matrix
$W_1'$ satisfy
\[
\|w_i'-w_i\|<\varepsilon
\qquad
\text{for every }i=1,\ldots,n_1.
\]
For every $i$ and $j$, the Cauchy--Schwarz inequality gives
\begin{align*}
\bigl|(w_i')^\top x_j-w_i^\top x_j\bigr|
&=\bigl|(w_i'-w_i)^\top x_j\bigr|\\
&\leq\|w_i'-w_i\|\,\|x_j\|\\
&<\varepsilon\\
&\leq\frac{\delta_j}{2}\\
&<|w_i^\top x_j|.
\end{align*}
 It follows that
\[
\operatorname{sgn}\bigl((w_i')^\top x_j\bigr)
=
\operatorname{sgn}(w_i^\top x_j)
\]
for every $i=1,\ldots,n_1$.

Hence $x_j$ has activation pattern $A_j$ with respect to
$\theta'$. Moreover, since
\[
(w_i')^\top x_j\neq 0
\qquad
\text{for every }i,
\]
the point $x_j$ lies in the interior of a
maximal-dimensional cone of $\mathcal C(\theta')$. Therefore
\[
A_j\in \mathcal{A}_{\theta'}.
\]
Since this holds for every $j=1,\ldots,N$, we conclude that
\[
\mathcal{A}_\theta\subseteq \mathcal{A}_{\theta'}
\qquad
\text{for every }\theta'\in\mathcal{N}_\theta.
\]
Thus $\theta$ is pattern-inclusive.

\emph{Case (ii).} If $W_1$ has at least one zero row, we consider the following subcases.

\emph{$W_1$ has at least one zero row and at least one nonzero row.}
Let
\[
I:=\{i\in\{1,\ldots,n_1\}:w_i\neq 0\}.
\]
By assumption, $I$ is nonempty. Choose an index
$k\notin I$, so that $w_k=0$.
Let $C$ be a maximal cell of $\mathcal C(\theta)$, i.e., a closed chamber of the arrangement of the hyperplanes $H_i$, $i\in I$. (Zero rows add no hyperplane.)

For each $i\in I$, let $\eta_i\in\{+1,-1\}$ denote the sign
of $w_i^\top x$ on $\operatorname{int}C$. Then
\[
\operatorname{int}C=
\left\{
x\in\mathbb{R}^{n_0}:
\eta_iw_i^\top x>0
\text{ for every }i\in I
\right\}.
\]
Define
\[
v:=\sum_{i\in I}\eta_iw_i.
\]
For every $x\in \operatorname{int}C$, we have
\[
v^\top x
=
\sum_{i\in I}\eta_iw_i^\top x>0.
\]
In particular, $v\neq 0$.

Let $A_C\in \mathcal{A}_\theta$ denote the activation pattern
associated with $C$. Since $w_k=0$, the $k$-th neuron is
inactive everywhere, so the $k$-th entry of $A_C$ is zero.
Now, for any $\varepsilon>0$, define a perturbed
first-layer weight matrix $W_1'$ by replacing the
$k$-th row of $W_1$ with
\[
(w_k')^\top=\varepsilon v^\top,
\]
while leaving all other rows unchanged. Let $\theta'$
denote the resulting parameter, with all other network
parameters unchanged.

For every $x\in \operatorname{int}C$,
\[
(w_k')^\top x
=
\varepsilon v^\top x>0.
\]
Thus the $k$-th neuron is active throughout $\operatorname{int}C$ under
$\theta'$.

We claim that
\[
A_C\notin \mathcal{A}_{\theta'}.
\]
Suppose, toward a contradiction, that
$A_C\in \mathcal{A}_{\theta'}$. Then there exists a
maximal-dimensional cone $C'$ of $\mathcal C(\theta')$ whose interior has activation pattern $A_C$.
Since each $w_i$, $i\in I$, is nonzero and is also a row of $W_1'$, $\operatorname{int}(C')$ does not meet $H_i$; hence the entries of $A_C$ give
\[
\eta_iw_i^\top x>0
\qquad
\text{for every }i\in I\text{ and every }x\in\operatorname{int}(C').
\]
Therefore
\[
\operatorname{int}(C')\subseteq \operatorname{int}C.
\]
However, the $k$-th neuron is active throughout $\operatorname{int}C$ under
$\theta'$, whereas the $k$-th entry of $A_C$ is zero.
This is a contradiction. Therefore
\[
A_C\notin \mathcal{A}_{\theta'},
\]
and consequently
\[
\mathcal{A}_\theta\nsubseteq \mathcal{A}_{\theta'}.
\]

Since $\varepsilon>0$ can be chosen arbitrarily small,
such a perturbation exists in every open neighborhood
of $\theta$. Thus $\theta$ is not pattern-inclusive.

\emph{$\mathit{W_1=0}$.}
In this case, every hidden neuron is inactive everywhere,
and hence
\[
\mathcal{A}_\theta=\{(0,\ldots,0)\}.
\]

Suppose first that $n_1\geq 2$. Choose any nonzero vector
$v\in\mathbb{R}^{n_0}$. For any $\varepsilon>0$, define
a perturbation $W_1'$ by setting
\[
w_1'=\varepsilon v,
\qquad
w_2'=-\varepsilon v,
\]
and leaving all remaining rows equal to zero.

On the interior of any maximal-dimensional cone of the perturbed
canonical polyhedral complex, we have either
\[
v^\top x>0
\qquad\text{or}\qquad
v^\top x<0.
\]
In the former case, the first neuron is active; in the latter case, the second neuron is active. Therefore the all-inactive pattern $(0,\ldots,0)$ cannot occur on the interior of any maximal-dimensional cone. Thus
\[
(0,\ldots,0)\notin \mathcal{A}_{\theta'},
\]
which implies
\[
\mathcal{A}_\theta\nsubseteq \mathcal{A}_{\theta'}.
\]
Since $\varepsilon>0$ can be arbitrarily small,
$\theta$ is not pattern-inclusive.

Finally, suppose that $n_1=1$. Then
\[
\mathcal{A}_\theta=\{(0)\}.
\]
For any perturbed parameter $\theta'$, let $w_1'$ denote
its first-layer weight vector.

If $w_1'=0$, the unique neuron remains inactive
everywhere, so $(0)\in \mathcal{A}_{\theta'}$.

If $w_1'\neq 0$, then the open half-space
\[
\{x\in\mathbb{R}^{n_0}:(w_1')^\top x<0\}
\]
is nonempty and full-dimensional. The unique neuron
is inactive throughout this region, so again
\[
(0)\in \mathcal{A}_{\theta'}.
\]
Therefore
\[
\mathcal{A}_\theta\subseteq \mathcal{A}_{\theta'}
\]
for every $\theta'$, and hence $\theta$ is
pattern-inclusive.

Combining the three cases gives the desired
characterization.
\end{proof}

\begin{example}\label{ex:noninclusive_deep_network}
As shown in Proposition~\ref{prop:pattern_inclusive_nbhd_exist},
the presence of zero rows in the first-layer weight matrix
plays a crucial role in determining pattern-inclusiveness
for unbiased shallow ReLU networks. However, for deep
networks, a parameter may fail to be pattern-inclusive
even when none of its weight matrices contains a zero row.

To see this, consider a parameter
$\theta\in\mathcal{P}_{2,2,2,1,1}$ of an unbiased
ReLU network
\[
f_\theta
=
W_4\circ\varsigma\circ W_3\circ\varsigma
\circ W_2\circ\varsigma\circ W_1,
\]
where
\[
W_1=
\begin{bmatrix}
1&0\\
0&1
\end{bmatrix},
\quad
W_2=
\begin{bmatrix}
1&1\\
1&1
\end{bmatrix},
\quad
W_3=
\begin{bmatrix}
1&-1
\end{bmatrix},
\quad
W_4=1.
\]

For any $x=(x_1,x_2)^\top\in\mathbb{R}^2$, let
\[
g(x):=\max\{0,x_1\}+\max\{0,x_2\}.
\]
Then
\[
\varsigma\circ W_2\circ \varsigma\circ W_1(x)
=
\begin{bmatrix}
g(x)\\
g(x)
\end{bmatrix}.
\]
Consequently, the preactivation of the neuron in the
third hidden layer is
\[
W_3\circ \varsigma\circ W_2\circ \varsigma\circ W_1(x)
=
g(x)-g(x)=0
\]
for every $x\in\mathbb{R}^2$.
Thus the third-layer neuron is inactive everywhere.

In particular, consider the maximal-dimensional cone
\[
C=\{x\in\mathbb{R}^2:x_1\geq 0,\ x_2\geq 0\}
\]
of $\mathcal C(\theta)$. The deeper layers do not subdivide the quadrants.
The activation pattern on $\operatorname{int}C$ is
\[
A=((1,1),(1,1),(0))\in \mathcal{A}_\theta.
\]

Now, for any $\varepsilon>0$, perturb only the
third-layer weight matrix by setting
\[
W_3'=
\begin{bmatrix}
1+\varepsilon&-1
\end{bmatrix},
\]
while leaving all other weight matrices unchanged.
Let $\theta'$ denote the resulting parameter.

Then
\[
W_3'\circ \varsigma\circ W_2\circ \varsigma\circ W_1(x)
=
\varepsilon g(x).
\]
For every $x\in \operatorname{int}C$, we have
\[
g(x)=x_1+x_2>0,
\]
and hence
\[
\varepsilon g(x)>0.
\]
Therefore, the third-layer neuron is active throughout
$\operatorname{int}C$ under $\theta'$.

Moreover, the interior of any maximal-dimensional cone of $\mathcal C(\theta')$ with
first-layer activation pattern $(1,1)$ is
contained in $\operatorname{int}C$. Since the third-layer neuron
is active throughout $\operatorname{int}C$ under $\theta'$, no
maximal-dimensional cone of the perturbed network
can have activation pattern $A$. Thus
\[
A\notin \mathcal{A}_{\theta'},
\]
which implies
\[
\mathcal{A}_\theta\nsubseteq \mathcal{A}_{\theta'}.
\]

Since $\varepsilon>0$ can be chosen arbitrarily small,
every open neighborhood of $\theta$ contains such
a parameter $\theta'$. Therefore $\theta$ is not
pattern-inclusive.
\end{example}

\subsection{Pattern Stability versus combinatorial stability}\label{subsec:combinatorial_versus_pattern}
\shallowstability*

\begin{proof}
\emph{(i)$\Rightarrow$(ii).} A combinatorial equivalence respecting the ternary labelings induces a bijection between the maximal-dimensional cells that preserves their activation patterns. Thus combinatorial stability implies pattern-stability.

\emph{(ii)$\Rightarrow$(iii)}. Suppose that some $r$ rows are linearly dependent. Choose a minimal linearly dependent subset $I$ and a nonzero relation
\[
\sum_{k\in I}\lambda_k w_k=0.
\]
By minimality, every $\lambda_k$ is nonzero. Replacing the relation by its negative if necessary, assume that some $\lambda_k>0$, and define $b_k=1$ if $\lambda_k>0$ and $b_k=0$ otherwise. If an input $x$ had binary pattern $b$ on $I$, then
$w_k^\top x>0$ when $\lambda_k>0$ and
$w_k^\top x\leq0$ when $\lambda_k<0$. Hence all terms $\lambda_k w_k^\top x$ would be nonnegative and at least one strictly positive, giving the contradiction
\[
0=\left(\sum_{k\in I}\lambda_k w_k\right)^\top x
=\sum_{k\in I}\lambda_k w_k^\top x>0.
\]

Minimality of $I$ also implies that removing one row from $I$ leaves an independent set whose span is proper, since $|I|\leq n_0$. Adding an arbitrarily small nonzero multiple of a vector outside this span to the removed row makes the modified row together with the remaining rows independent. Thus every neighborhood of $\theta$ contains such a parameter $\theta'$, obtained by perturbing only one row of $W_1$ and leaving all other rows and $W_2$ unchanged.

The matrix with rows $(w_k')^\top$, $k\in I$, has full row rank. We can therefore choose $x_0$ satisfying $(w_k')^\top x_0=2b_k-1$ for every $k\in I$. Consequently,
\[
O\coloneqq\{x\in\mathbb R^{n_0}:
(2b_k-1)(w_k')^\top x>0\text{ for every }k\in I\}
\]
is nonempty and open. It meets the interior of a
maximal-dimensional cell of $\mathcal C(\theta')$, whose pattern restricts to $b$ on $I$. Since no input at $\theta$ has this restriction, $\mathcal{A}_{\theta'}\neq \mathcal{A}_\theta$, contradicting pattern-stability.

\emph{(iii)$\Rightarrow$(i).} Assume (iii). All rows of $W_1$ are nonzero and remain so near $\theta$. For each sign vector $X\in\{-1,1\}^{n_1},$ let $X_k$ denote the sign of the $k$-th neuron and define
\[
\begin{aligned}
K_X(\theta')
&=\{x\in\mathbb R^{n_0}:
X_k(w_k')^\top x>0\text{ for every }k=1,\dots,n_1\} \text{ and}
\\
\mathcal T_{\theta'}
&=\{X\in\{-1,1\}^{n_1}:K_X(\theta')\neq\varnothing\}.
\end{aligned}
\]
The elements of $\mathcal{T}_{\theta'}$ are the sign vectors of the maximal-dimensional
cones of the first-layer arrangement. We first prove $\mathcal{T}_{\theta} = \mathcal T_{\theta'}$ near $\theta$ by showing that existing sign vectors persist and that absent sign vectors remain absent.

For each $X\in\mathcal T_\theta$, choose $x\in K_X(\theta)$.
The finitely-many strict inequalities at this fixed input
persist under sufficiently small perturbations, so $x$
still realizes $X$. Since $\mathcal T_\theta$ is finite,
there is an open neighborhood $\mathcal N_0$ on which
\[
\mathcal T_\theta\subseteq\mathcal T_{\theta'}.
\]
If $n_1\leq n_0$, condition (iii) makes $W_1$ surjective,
so $W_1x=X$ is solvable for every $X\in\{-1,1\}^{n_1}$. Thus every sign vector already occurs, and equality $\mathcal T_\theta=\mathcal T_{\theta'}$ holds throughout $\mathcal N_0$. In this case, set $\mathcal N_\theta=\mathcal N_0$.

Now suppose that $n_1>n_0$, fix $X\notin\mathcal T_\theta$,
and put $u_k=X_kw_k$. Since the inequalities
$u_k^\top x>0$ have no simultaneous solution, Gordan's
theorem of the alternative
\cite[Example~2.21]{boyd2004convex}
gives a nonzero nonnegative dependence. Choose one with
support $S$ of smallest cardinality:
\[
\sum_{k\in S}\lambda_k u_k=0,
\qquad \lambda_k>0\quad(k\in S).
\]

The space of linear relations among the vectors $\{u_k:k\in S\}$ is spanned by $\lambda=(\lambda_k)_{k \in S}$. Indeed, if $\mu=(\mu_k)_{k \in S}$ were another relation not obtainable through multiplication by $\lambda$, we could choose its sign so that some $\mu_k>0$ and set
\[
t=\min_{\substack{k\in S\\\mu_k>0}}
\frac{\lambda_k}{\mu_k}.
\]
Then $\lambda-t\mu$ is a nonnegative coefficient
vector with at least one zero entry. It is nonzero
because $\mu$ is not proportional to $\lambda$ and
\[
\sum_{k\in S}(\lambda_k-t\mu_k)u_k=0.
\]
This contradicts the minimality of $S$. Therefore, rank--nullity gives
\[
|S|-1=\dim\operatorname{span}\{u_k:k\in S\}\leq n_0.
\]
On the other hand, (iii) implies that every set of at most
$n_0$ vectors $u_k$ is independent, since multiplying rows
by signs preserves independence. As the vectors indexed
by $S$ are dependent, $|S|\geq n_0+1$. Hence
$|S|=n_0+1$.

Fix $k_0\in S$. The remaining $n_0$ vectors form a basis,
and the positive dependence gives
\[
-u_{k_0}
=\sum_{k\in S\setminus\{k_0\}}a_k u_k,
\qquad a_k=\frac{\lambda_k}{\lambda_{k_0}}>0.
\]
For nearby $\theta'$, put $u_k'=X_kw_k'$ and let
$\mathcal J_{\theta'}$ be the $n_0\times n_0$ matrix with columns $u_k'$, $k\in S\setminus\{k_0\}$,
in a fixed order. This matrix remains invertible, and
the coefficients in
\[
-u_{k_0}'
=\sum_{k\in S\setminus\{k_0\}}a_k(\theta')u_k'
\]
are given by
$a(\theta')=-\mathcal J_{\theta'}^{-1}u_{k_0}'$.
They therefore vary continuously and remain positive
near $\theta$. If $x\in K_X(\theta')$, taking inner
products in this relation gives
\[
0=(u_{k_0}')^\top x
+\sum_{k\in S\setminus\{k_0\}}
a_k(\theta')(u_k')^\top x>0,
\]
a contradiction. Thus each absent sign vector $X$
remains absent on an open neighborhood $\mathcal U_X$
of $\theta$.
Intersecting $\mathcal N_0$ with these finitely many
neighborhoods gives an open neighborhood $\mathcal N_\theta$ on which all existing sign vectors persist and all absent sign vectors remain absent. Hence
\[
\mathcal T_{\theta'}=\mathcal T_\theta
\qquad\text{for every }\theta'\in\mathcal N_\theta.
\]

Now consider the set of all ternary sign vectors
\[
\mathcal L_{\theta'}=
\{\operatorname{sgn}(W_1'x):x\in\mathbb R^{n_0}\}.
\]
This is the covector set of the realizable oriented matroid associated with the first-layer arrangement, while $\mathcal T_{\theta'}$ is its set of topes.
By Mandel's theorem \cite[Theorem~4.2.13]{BLSWZ99}, the covector set is uniquely determined by its set of topes. Hence $\mathcal L_{\theta'}=\mathcal L_\theta$ throughout $\mathcal N_\theta$.
The covectors determine the ternary-labeled faces and their incidence relations, as well as the underlying matroid. Corresponding faces have the same dimensions, since the face with covector $X$ has dimension
\[
n_0-\operatorname{rank}\{w_k':X_k=0\}.
\]

It remains to construct a continuous family of homeomorphisms realizing these correspondences. For each
$\theta'\in\mathcal N_\theta$, let
\[
V_{\theta'}=\operatorname{span}\{w_1',\ldots,w_{n_1}'\},
\qquad
K_{\theta'}=\ker W_1'.
\]
Then $\mathbb R^{n_0}=V_{\theta'}\oplus K_{\theta'}$.
Since the underlying matroid is unchanged, the ranks of all row submatrices are constant. Consequently, these subspaces and the linear spans of corresponding faces vary continuously with $\theta'$.

Every cell $C$ of $\mathcal C(\theta')$ satisfies
\[
C=(C\cap V_{\theta'})\oplus K_{\theta'}.
\]
The cones $C\cap V_{\theta'}$ form a fan of strongly convex cones in $V_{\theta'}$, which we call the essential fan at $\theta'$. We first construct a homeomorphism between the essential fans at $\theta$ and $\theta'$.

Let $\mathcal L$ denote the common covector set
$\mathcal L_\theta=\mathcal L_{\theta'}$.
For each $X\in\mathcal L$, let $F_X(\theta')$ denote
the closed face of the essential fan at $\theta'$
with covector $X$. For each nonzero $X\in\mathcal L$, choose a point in the relative interior
\[
u_X(\theta)\in\operatorname{relint}F_X(\theta),
\]
and let $u_X(\theta')$ be the orthogonal projection of $u_X(\theta)$ onto $\operatorname{span}F_X(\theta')$.
Since these spans vary continuously with $\theta'$,
the points $u_X(\theta')$ also vary continuously.
They satisfy the equalities defining $F_X(\theta')$,
and, for sufficiently small perturbations, the strict
sign inequalities remain satisfied. Hence
\[
u_X(\theta')\in\operatorname{relint}F_X(\theta').
\]
As there are finitely many faces, we may shrink
$\mathcal N_\theta$ so that this holds simultaneously
for every nonzero $X\in\mathcal L$.

Now subdivide the essential fan at $\theta$ into
simplicial cones using the rays generated by the
chosen points $u_X(\theta)$. Each nonzero simplicial
cone corresponds to a chain of nonzero faces
\[
F_{X_1}(\theta)\subsetneq\cdots
\subsetneq F_{X_\ell}(\theta)
\]
and is generated by $u_{X_1}(\theta),\ldots,u_{X_\ell}(\theta)$,
where $\ell$ is the length of the chain.
Since the face relations are preserved, the corresponding points $u_{X_1}(\theta'),\ldots,u_{X_\ell}(\theta')$
generate the simplicial cones of the corresponding
subdivision of the essential fan at $\theta'$.

We can therefore define
\[
H_{\theta'}:V_\theta\longrightarrow V_{\theta'}
\]
on each simplicial cone by
\[
H_{\theta'}\left(
\sum_{i=1}^{\ell}a_i u_{X_i}(\theta)
\right)
=
\sum_{i=1}^{\ell}a_i u_{X_i}(\theta'),
\qquad a_i\geq0,
\]
and set $H_{\theta'}(0)=0$. The generators of each simplicial cone are linearly independent, so the coefficients $a_i$ are unique. The linear maps agree on common faces, and the reverse correspondence gives a continuous inverse. Thus $H_{\theta'}$ is a piecewise-linear homeomorphism preserving the ternary-labeled faces.

On each fixed simplicial cone in $V_\theta$, the
coefficients $a_i$ depend continuously on the input
$v$, and the points $u_{X_i}(\theta')$ depend continuously on $\theta'$. Hence $H_{\theta'}(v)$ is continuous in both $\theta'$ and $v$ on that cone.
Since there are finitely many closed simplicial cones
and the maps agree on their intersections, the map
\[
(\theta',v)\longmapsto H_{\theta'}(v)
\]
is continuous on $\mathcal N_\theta\times V_\theta$.

Finally, since $K_{\theta'}$ varies continuously and
has constant dimension, after shrinking
$\mathcal N_\theta$ once more, the orthogonal projection
\[
Q_{\theta'}=P_{K_{\theta'}}|_{K_\theta}:
K_\theta\longrightarrow K_{\theta'}
\]
is a linear isomorphism.

Using the decomposition
$\mathbb R^{n_0}=V_\theta\oplus K_\theta$, define
\[
\Phi_{\theta'}(v+k)
=
H_{\theta'}(v)+Q_{\theta'}(k),
\qquad v\in V_\theta,\quad k\in K_\theta.
\]
Then $\Phi_{\theta'}$ is a homeomorphism of
$\mathbb R^{n_0}$ mapping each cell of
$\mathcal C(\theta)$ onto its corresponding cell
of $\mathcal C(\theta')$ preserving the ternary
activation labels.
Moreover, the map $(\theta',x)\mapsto\Phi_{\theta'}(x)$ is continuous on $\mathcal N_\theta\times\mathbb R^{n_0}$, and $\Phi_\theta$ is the identity. Therefore, these maps form the continuous family of homeomorphisms required by combinatorial stability. Hence $\theta$ is combinatorially stable.

For the genericity statement, let $W_{1,I}$ denote the submatrix of $W_1$ consisting of the rows indexed by $I$. The complement of $\mathcal G_{n_0,n_1,n_2}$ is the zero set of
\[
p(W_1)=\prod_{\substack{I\subseteq\{1,\ldots,n_1\}\\ |I|=r}}
\det(W_{1,I}W_{1,I}^{\top}).
\]
Each factor is a nonzero polynomial, since $r\leq n_0$.
Thus $p$ is nonzero. Its zero set is closed, has empty interior, and has Lebesgue measure zero, proving the remaining claims.
\end{proof}

\begin{remark}
 The equivalence in \Cref{thm:shallow_pattern_combinatorial_stable} is a local stability statement and does not by itself compare two arbitrary parameters. For unbiased shallow networks, $\mathcal{A}_{\theta_1}=\mathcal{A}_{\theta_2}$ does imply that
$\mathcal C(\theta_1)$ and $\mathcal C(\theta_2)$ are combinatorially equivalent by an equivalence respecting the ternary labelings (\Cref{thm:hierarchy_parameters}(i)). This fails in general. Masden \cite[Theorem~15]{Masden2025} gives two ReLU networks with biases whose
maximal-dimensional cells carry the same activation patterns but whose canonical polyhedral complexes are not combinatorially equivalent. Without biases,
\Cref{ex:pattern_combinatorial} gives two parameters of architecture $(3,4,1,1)$ with equal canonical pattern sets that are not combinatorially equivalent.
\end{remark}

If a parameter $\theta\in\mathcal P_{n_0,\ldots,n_m}$ is combinatorially stable, then $\theta$ is pattern-stable. This follows from the definition:
the ternary-label-preserving bijection between
maximal-dimensional cells preserves their binary
activation patterns. The converse holds for shallow networks (\Cref{thm:shallow_pattern_combinatorial_stable}), but fails in general. The following is a counterexample to the converse for deep networks.
\begin{example}\label{ex:pattern_combinatorial}
    Consider the architecture $(3,4,1,1)$
and the parameter $\theta$ given by
\[
W_1=
\begin{bmatrix}
1&0&0\\
0&1&0\\
1&-2&1\\
0&1&1
\end{bmatrix},
\qquad
W_2=
\begin{bmatrix}
-2&4&1&-1
\end{bmatrix},
\qquad
W_3=1.
\]
For an input $(x,y,z)\in\mathbb R^3$, write
\[
q_1=x,\qquad q_2=y,\qquad
q_3=x-2y+z,\qquad q_4=y+z.
\]
These first-layer preactivations satisfy the linear relation
\begin{equation}\label{eq:equation_of_q_i}
q_3=q_1-3q_2+q_4.
\end{equation}

The preactivation of the neuron in the second hidden layer is
\[
h=-2\ReLU(q_1)+4\ReLU(q_2)
+\ReLU(q_3)-\ReLU(q_4).
\]

We first show that $\theta$ is pattern-stable.
Every three rows of $W_1$ are linearly independent,
so \Cref{thm:shallow_pattern_combinatorial_stable}
implies that the first-layer sign vectors remain
unchanged under sufficiently small perturbations.
Writing these vectors as strings of $\texttt{+}$ and $\texttt{-}$ signs,
the fourteen realized vectors are all those except
$\texttt{+--+}$ and $\texttt{-++-}$.

Direct substitution shows that $h$ takes both positive
and negative values on each of the eight first-layer
cones with sign vectors
\[
\begin{gathered}
\texttt{--++},\quad
\texttt{-+-+},\quad
\texttt{-+++},\quad
\texttt{+-+-},\\
\texttt{+-++},\quad
\texttt{++--},\quad
\texttt{++-+},\quad
\texttt{++++}.
\end{gathered}
\]
Choosing points satisfying the corresponding strict
inequalities shows that both activation statuses
persist on full-dimensional regions under sufficiently
small perturbations.

On the cone with sign vector $\texttt{----}$,
all first-layer neurons are inactive, so $h$ remains
identically zero. On each of the four cones
\[
\texttt{---+},\quad
\texttt{--+-},\quad
\texttt{-+--},\quad
\texttt{+---},
\]
exactly one first-layer neuron is active. Since every
entry of $W_2$ is nonzero, the sign of $h$ on each of
these cones is preserved under sufficiently small
perturbations.

It remains to consider $\texttt{+++-}$. By
\eqref{eq:equation_of_q_i}, on this cone we have
\[
h=-2q_2-q_3+2q_4<0.
\]
Since $q_2,q_3,q_4$ are linearly independent,
the perturbed preactivation on the corresponding
cone can be written as
\[
h_{\theta'}=a(\theta')q_2'
+b(\theta')q_3'+c(\theta')q_4',
\]
where the coefficients vary continuously and equal
$(-2,-1,2)$ at $\theta$. Thus, near $\theta$,
$a(\theta'),b(\theta')<0$ and $c(\theta')>0$.
Since $q_2',q_3'>0$ and $q_4'<0$ on this cone,
$h_{\theta'}<0$ throughout it.

These cases exhaust the fourteen first-layer sign
vectors. Thus all canonical activation patterns persist
and no new ones appear under sufficiently small perturbations. Hence there is an open neighborhood
$\mathcal N_\theta$ such that
\[
\mathcal{A}_{\theta'}=\mathcal{A}_\theta
\qquad\text{for every }\theta'\in\mathcal N_\theta,
\]
so $\theta$ is pattern-stable.

We now show that $\theta$ is not combinatorially stable.
The ray $\rho=\{(0,0,t):t\geq0\}$ is a face of the closed first-layer cone with sign vector $\texttt{++++}$, obtained by imposing
$q_1=q_2=0$. Along $\rho$, we have $(q_1,q_2,q_3,q_4,h)=(0,0,t,t,0)$. Thus the second hidden layer does not subdivide $\rho$, so it is a cell of $\mathcal C(\theta)$ whose relative interior has ternary label $(0,0,1,1,0)$.

For $\varepsilon>0$, perturb only $W_2$ by setting
\[
W_2^\varepsilon=
\begin{bmatrix}
-2&4&1&-1+\varepsilon
\end{bmatrix},
\]
and denote the resulting parameter by $\theta_\varepsilon$.
Then
\[
h_\varepsilon(0,0,t)=\varepsilon t>0
\qquad\text{for } t>0.
\]
Since the first-layer ternary label $(0,0,1,1)$ occurs
only on the relative interior of $\rho$, no cell of
$\mathcal C(\theta_\varepsilon)$ has label $(0,0,1,1,0)$.
Therefore there is no combinatorial equivalence
preserving the ternary labelings. As $\varepsilon$
can be arbitrarily small, $\theta$ is not
combinatorially stable.
\end{example}

\begin{figure}[!htbp]
  \centering
  \captionsetup{font=small,skip=3pt}
  \captionsetup[subfigure]{font=footnotesize,skip=1pt}

  \begin{subfigure}[t]{0.4\linewidth}
    \centering
    \includegraphics[width=\linewidth,height=0.16\textheight,keepaspectratio]
      {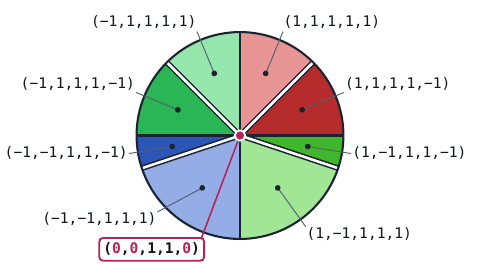}
    \caption{$\theta$: ternary labels.}
    \label{fig:stability-a}
  \end{subfigure}\hspace{0.01\linewidth}%
  \begin{subfigure}[t]{0.4\linewidth}
    \centering
    \includegraphics[width=\linewidth,height=0.16\textheight,keepaspectratio]
      {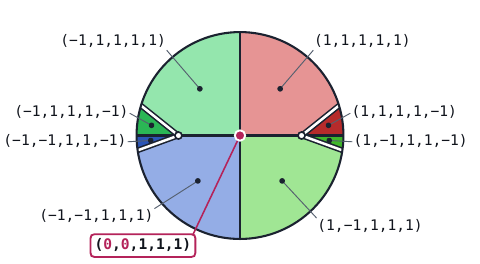}
    \caption{$\theta_\varepsilon$: ternary labels.}
    \label{fig:stability-b}
  \end{subfigure}

  \par\vspace{2pt}

  \begin{subfigure}[t]{0.4\linewidth}
    \centering
    \includegraphics[width=\linewidth,height=0.16\textheight,keepaspectratio]
      {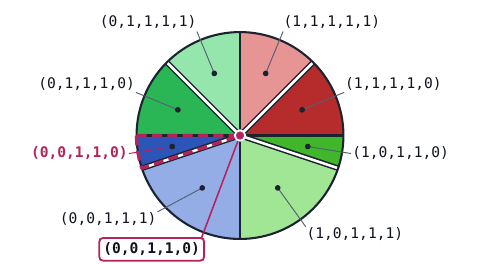}
    \caption{$\theta$: binary patterns.}
    \label{fig:stability-c}
  \end{subfigure}\hspace{0.01\linewidth}%
  \begin{subfigure}[t]{0.4\linewidth}
    \centering
    \includegraphics[width=\linewidth,height=0.16\textheight,keepaspectratio]
      {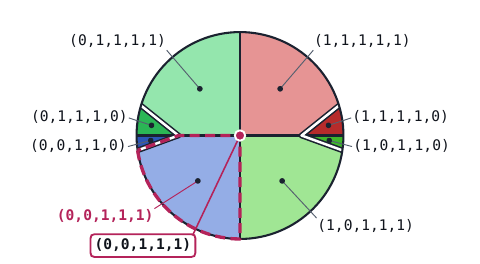}
    \caption{$\theta_\varepsilon$: binary patterns.}
    \label{fig:stability-d}
  \end{subfigure}

  \caption{Slices at $z=1$ for the parameters in
    \Cref{ex:pattern_combinatorial},
    with $\varepsilon=1/4$. Labels list first-layer neurons
    before the second-layer neuron. Matching colors indicate
    identical binary patterns. Dashed outlines mark the region sharing the pattern of the red point $(0,0,1)\in\rho$.}
  \label{fig:stability}
\end{figure}

In the slice $z=1$ shown in \Cref{fig:stability}, the ternary label $(0,0,1,1,0)$ disappears, while the canonical activation patterns remain unchanged.

\newpage 

\section{Generators for shallow networks}\label{sec:generators}

Throughout this section, we consider an unbiased network with architecture $(n_0,n_1,n_2)$.

Write $w_k^\top$ for row $k$ of $W_1$, $c_k$ for
column $k$ of $W_2$, and $V_k=c_kw_k^\top$. For a pattern
$A\in\{0,1\}^{n_1}$, we have
\begin{equation}\label{eq:shallow_rankone_sum}
M_A(\theta)=\sum_{k=1}^{n_1}A_kV_k,
\qquad \operatorname{rank}V_k\leq1.
\end{equation}
We fix a set of distinct patterns $\mathcal A=\{A_1,\ldots,A_N\}$
and use $Y_A$ to denote the $n_2\times n_0$ block of indeterminates indexed by $A$. The set $\mathcal A$ remains fixed while the parameter in $\varphi^{\mathcal A}$ varies over the entire parameter
space. We show that if we take $\mathcal{A}=\mathcal{A}_\theta$, where $\mathcal{A}_\theta$ is the canonical set of activation patterns (\Cref{def:canonical_set_of_activation_patterns}) whose first-layer matrices have rows that are nonzero and pairwise linearly independent, then the ideal
$J^{\mathcal A}$ is generated by linear relations and
$2\times2$ minors.

\subsection{Bending relations}
\label{subsec:shallow_bending}

\begin{lemma}[Bending relations]\label{lem:intersectionnumbersymmetrytogenerators}
Suppose $A,A',B,B'\in\mathcal A$ satisfy
$A'-A=B'-B=e_k$. Then, for every real parameter $\theta$,
\[
M_{A'}(\theta)-M_A(\theta)
=M_{B'}(\theta)-M_B(\theta)=c_kw_k^\top.
\]
Thus the entries of $Y_{A'}-Y_A-Y_{B'}+Y_B$ and the
$2\times2$ minors of $Y_{A'}-Y_A$ belong to $J^{\mathcal A}$.
More generally, the same linear relation holds whenever
$A'-A=B'-B$.
\end{lemma}

\begin{proof}
Subtracting \eqref{eq:shallow_rankone_sum} gives
$M_{A'}-M_A=\sum_k(A'_k-A_k)V_k$.
For $A'-A=e_k$, this is $V_k=c_kw_k^\top$, which has rank
at most one.
\end{proof}

We call these linear equalities \emph{bending relations}.
They and the accompanying minors are special cases of
\cite[eq.~(5), Theorem~34]{alexandr2025constraining}:
\begin{equation}\label{eq:shallow_rank_bound}
\operatorname{rank}\!\left(\sum_{A\in\mathcal A}\lambda_AM_A\right)
\leq
\left|\operatorname{supp}\!\left(\sum_{A\in\mathcal A}\lambda_AA\right)\right|,
\qquad \lambda\in\mathbb Z^N.
\end{equation}
Here $\operatorname{supp}$ is the set of nonzero coordinates.
The two cases above have support size zero and one, respectively.

\subsection{Generators from the wall graph}
\label{subsec:wall_graph_generators}

\begin{definition}[Wall graph]\label{def:wall_graph}
For a nonempty set $\mathcal A\subseteq\{0,1\}^{n_1}$,
the \emph{wall graph} $\Gamma_{\mathcal A}$ is the edge-labeled graph with vertex set $\mathcal A$ constructed by joining two patterns by an edge labeled $k$ if they differ only in coordinate $k$.
\end{definition}

\begin{lemma}\label{lem:cube}
Let $\theta=(W_1,W_2)$ be a shallow parameter whose first-layer
rows are nonzero and pairwise linearly independent.
Then $\Gamma_{\mathcal{A}_\theta}$ is connected and has an edge labeled
$k$ for every $k=1,\ldots,n_1$.
\end{lemma}

\begin{proof}
The maximal cells are the closures of the chambers of the arrangement formed by the distinct hyperplanes 
$$H_k=\{x:w_k^\top x=0\}.$$
For each $k$, choose $z\in H_k\setminus\bigcup_{j\ne k}H_j$; such a point exists because the intersections $H_k\cap H_j$ are proper subspaces of $H_k$. On the two sides of $H_k$ near $z$, only the $k$-th sign changes. This gives an edge labeled $k$.

To prove connectedness, choose $x$ in one chamber and $y$ in another so that the line segment
\[
[x,y]\coloneqq\{(1-t)x+ty:0\leq t\leq1\}
\]
does not meet $H_j\cap H_k$ for any $j\neq k$. This choice is possible: for fixed $x$, every forbidden $y$ lies in one of the proper subspaces $\operatorname{span}(x,H_j\cap H_k)$, and finitely many proper subspaces cannot cover an open chamber. The segment crosses hyperplanes one at a time, so its successive chamber patterns form a path in $\Gamma_{\mathcal{A}_\theta}$. If $n_0=1$, the hypotheses force $n_1=1$, and the graph is simply the edge joining the vertices $(0)$ and $(1)$.
\end{proof}

The next theorem applies to these canonical sets and, more
generally, to every pattern set with the two graph properties
in \Cref{lem:cube}.

\shallowgenerators*

\begin{proof}
\emph{Recovering the neuron matrices.}
The chosen pairs satisfy
\[
\adjdiff{k}\bigl(\slopeembed{\mathcal A}(V_1,\ldots,V_{n_1})\bigr)=V_k.
\]
Thus $\slopeembed{\mathcal A}$ is injective, with inverse given by the differences $\adjdiff{k}$ on its image. Every complex matrix of rank at most one factors as $cw^\top$. Therefore the polynomial map
\[
\kappa:(W_1,W_2)\longmapsto(c_kw_k^\top)_{k=1}^{n_1}
\]
maps the complex weight space onto $\rankonevar^{n_1}$, and $\varphi^{\mathcal A}=\slopeembed{\mathcal A}\circ\kappa$. The real weights are Zariski dense in the complex weight space, so using complex weights does not change $U^{\mathcal A}$. Since an injective linear map is an isomorphism onto a closed linear subspace, $\slopeembed{\mathcal A}(\rankonevar^{n_1})$
is closed. Hence
\begin{equation}\label{eq:shallow_product_image}
U^{\mathcal A}=\slopeembed{\mathcal A}(\rankonevar^{n_1}).
\end{equation}

\emph{The linear relations.}
Let $L$ be the ideal generated by (G1) and (G2).
Along an edge $B'-B=e_k$, relation (G1) gives
$Y_{B'}\equiv Y_B+\adjdiff{k}\pmod L$.
Following a path from $A_*$ to any $A$ and using (G2) therefore gives
\begin{equation}\label{eq:shallow_path_reconstruction}
Y_A\equiv\sum_{k=1}^{n_1}A_k\adjdiff{k}\pmod L.
\end{equation}
Indeed, the signed number of changes in coordinate $k$ along
the path is $A_k-(A_*)_k$.

Let $Z_1,\ldots,Z_{n_1}$ be independent $n_2\times n_0$
matrices of indeterminates. The substitutions
\[
Y_A\longmapsto\sum_k A_kZ_k,
\qquad
Z_k\longmapsto\adjdiff{k}\pmod L
\]
give inverse ring homomorphisms
\begin{equation}\label{eq:shallow_linear_quotient}
\mathbb C[\mathbf Y]/L\cong\mathbb C[Z_1,\ldots,Z_{n_1}].
\end{equation}
The first is well defined because it annihilates (G1) and (G2);
the two compositions are the identity by
\eqref{eq:shallow_path_reconstruction} and the definition of
$\adjdiff{k}$.

\emph{The remaining relations.}
For a generic matrix $Z$, its $2\times2$ minors generate the prime ideal of the matrices of rank at most one \cite[Theorem~2.10 and Remark~2.12]{brunsVetter1988}. 
Since every matrix of rank at most one appears as $cw^\top$, it follows that, for each $k$, the kernel of the homomorphism $\mathbb C[Z_k]\to\mathbb C[c_k,w_k]$ given entrywise by $Z_k\mapsto c_kw_k^\top$ is the ideal  $I_2(Z_k)$ generated by the $2\times2$ minors of $Z_k$. Because the blocks $Z_k$, and likewise the pairs $(c_k,w_k)$, involve disjoint sets of variables, and a tensor product over $\mathbb C$ of injective linear maps is injective, the homomorphism
\[
\kappa^*:\mathbb C[Z_1,\ldots,Z_{n_1}]\longrightarrow\mathbb C[\mathbf W],
\qquad Z_k\longmapsto c_kw_k^\top,
\]
has kernel
\[
I_2(Z_1)+\cdots+I_2(Z_{n_1}).
\]
Since $\mathbb C[\mathbf W]$ is an integral domain and $\kappa$ maps onto $\rankonevar^{n_1}$, this ideal is prime and is the defining ideal of $\rankonevar^{n_1}$.

If $n_0=1$ or $n_2=1$ then $I_2(Z_k)=0$, since then $Z_k$ has no $2\times2$ minors. Under \eqref{eq:shallow_linear_quotient}, the class of a polynomial $F$ maps to $F((\sum_k A_kZ_k)_{A\in\mathcal A})$. By \eqref{eq:shallow_product_image}, $F\in J^{\mathcal A}$ exactly when the polynomial obtained via this substitution vanishes on $\rankonevar^{n_1}$, or, equivalently, belongs to the displayed determinantal ideal. Its inverse image is generated by (G1)--(G3), proving ideal generation and primality.
Finally, $\dim\rankonevar=n_0+n_2-1$
\cite[Proposition~1.1]{brunsVetter1988}, which gives the dimension formula.
\end{proof}

Relation (G1) expresses the bending relations, which require that the difference across every edge $k$ equals $\Delta_k$. When the wall graph is connected, the bending relations determine the blocks $Y_A$ up to the addition of a common matrix. Relation (G2) fixes $Y_{A_*}$ in terms of the edge differences $\Delta_k$ to give a unique linear expression for each $Y_A$. In particular, (G2) is one matrix relation, consisting of $n_2n_0$ scalar equations. It cannot be omitted: setting every $Y_A$ equal to the same nonzero matrix satisfies (G1) and (G3), but violates (G2). When $n_2=1$, the minors in (G3) are vacuous, so $J^{\mathcal A}$ is generated by linear forms.

\Cref{thm:shallow_generators} gives a finite generating set of the kind requested in \cite[Problem~35]{alexandr2025constraining} for these pattern sets, including the canonical sets in \Cref{lem:cube}. Every generator in \Cref{thm:shallow_generators} is obtained from \eqref{eq:shallow_rank_bound}. For (G1), take coefficients $(1,-1,-1,1)$ on $(B',B,A^{k,1},A^{k,0})$; for (G2), expand $Y_{A_*}-\sum_k(A_*)_k\adjdiff{k}$ as an integer combination of the blocks; for (G3), take $(1,-1)$ on $(A^{k,1},A^{k,0})$. If a pattern occurs more than once, combine its coefficients. The resulting combinations of patterns are zero for (G1) and (G2), and $e_k$ for (G3). However, this generating set need not be minimal.

\begin{figure}[htbp]
\centering
\begin{tikzpicture}[scale=1.35,>=Stealth]
  \draw[thick] (-2,0) -- (2,0) node[right] {$H_1$};
  \draw[thick] (-1.6,-1.6) -- (1.6,1.6) node[above right] {$H_2$};
  \draw[thick] (0,-1.8) -- (0,1.8) node[above] {$H_3$};
  \draw[->,thick] (-1.75,0) -- (-1.75,0.4) node[above] {$w_1$};
  \draw[->,thick] (1.4,1.4) -- (1.1,1.7) node[above left] {$w_2$};
  \draw[->,thick] (0,1.55) -- (-0.4,1.55) node[left] {$w_3$};
  \node[align=center] at (1.25,0.5) {$\sigma_1$\\[-2pt] \scriptsize $(1,0,0)$};
  \node[align=center] at (0.5,1.15) {$\sigma_2$\\[-2pt] \scriptsize $(1,1,0)$};
  \node[align=center] at (-0.85,0.8) {$\sigma_3$\\[-2pt] \scriptsize $(1,1,1)$};
  \node[align=center] at (-1.3,-0.5) {$\sigma_4$\\[-2pt] \scriptsize $(0,1,1)$};
  \node[align=center] at (-0.5,-1.25) {$\sigma_5$\\[-2pt] \scriptsize $(0,0,1)$};
  \node[align=center] at (0.95,-0.9) {$\sigma_6$\\[-2pt] \scriptsize $(0,0,0)$};
\end{tikzpicture}
\caption{The six cones in \Cref{ex:threehyperplanesin2D}, for normals
$(0,1),(-1,1),(-1,0)$. Arrows point to the active side of each neuron.}
\label{fig:threehyperplanesin2D}
\end{figure}
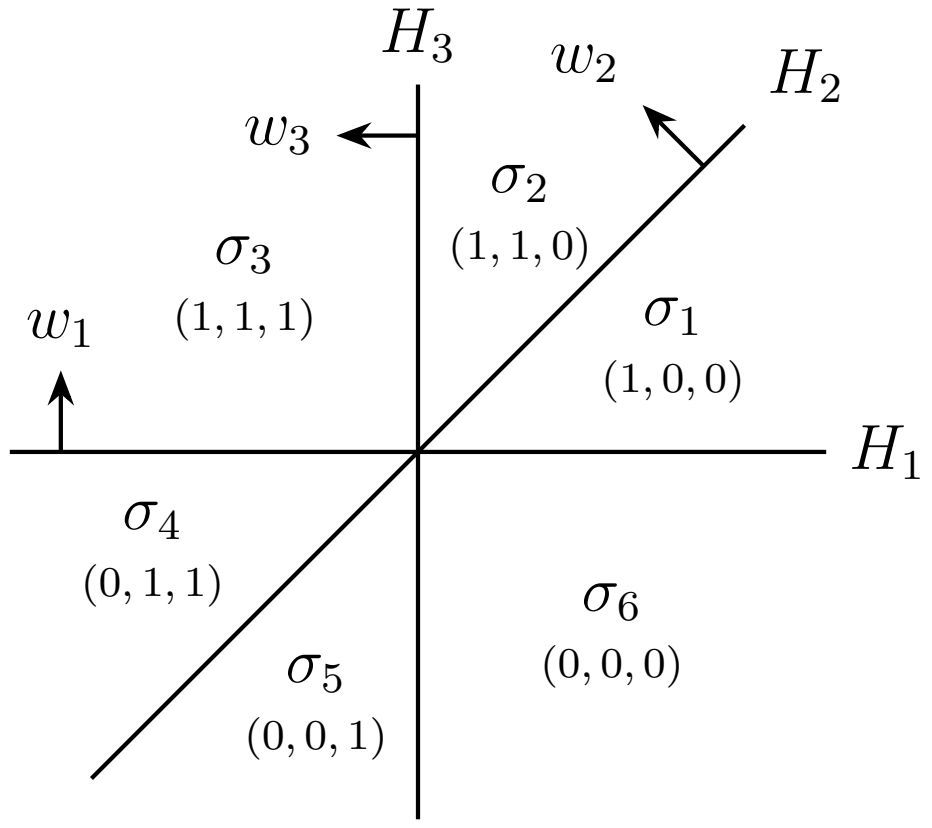

The following example illustrates these three types of generators for a simple planar arrangement. In this case, (G1) gives the bending relations, (G2) fixes the block corresponding to the all-inactive pattern, and (G3) is vacuous because the output dimension is one.
\begin{example}\label{ex:threehyperplanesin2D}
Consider the architecture $(2,3,1)$ and a parameter
$\theta=(W_1,W_2)$ whose first-layer rows are
\[
w_1^\top=(0,1),\qquad
w_2^\top=(-1,1),\qquad
w_3^\top=(-1,0).
\]
Let $\sigma_1,\ldots,\sigma_6$ be the maximal cones shown in
\Cref{fig:threehyperplanesin2D}, listed in counterclockwise order.
Their activation patterns are
\[
\begin{array}{lll}
A_1=(1,0,0),&A_2=(1,1,0),&A_3=(1,1,1),\\
A_4=(0,1,1),&A_5=(0,0,1),&A_6=(0,0,0).
\end{array}
\]

We apply \Cref{thm:shallow_generators}. Choose the reference pairs
\[
(A^{1,0},A^{1,1})=(A_6,A_1),\qquad
(A^{2,0},A^{2,1})=(A_1,A_2),\qquad
(A^{3,0},A^{3,1})=(A_2,A_3),
\]
and take $A_*=A_6$. Thus
\[
\adjdiff{1}=Y_{A_1}-Y_{A_6},\qquad
\adjdiff{2}=Y_{A_2}-Y_{A_1},\qquad
\adjdiff{3}=Y_{A_3}-Y_{A_2}.
\]

The remaining three edges of the wall graph give the bending
relations
\begin{align*}
Y_{A_3}-Y_{A_4}&=\adjdiff{1},
\\
Y_{A_4}-Y_{A_5}&=\adjdiff{2},
\\
Y_{A_5}-Y_{A_6}&=\adjdiff{3}.
\end{align*}
Equivalently, writing
$M_{\sigma_i}=M_{A_i}(\theta)$, these relations are
\begin{align*}
M_{\sigma_6}-M_{\sigma_1}
&=M_{\sigma_4}-M_{\sigma_3},\\
M_{\sigma_1}-M_{\sigma_2}
&=M_{\sigma_5}-M_{\sigma_4},\\
M_{\sigma_2}-M_{\sigma_3}
&=M_{\sigma_6}-M_{\sigma_5}.
\end{align*}
Their sum is identically zero, so any two of them determine the third.

Since $A_*=A_6=(0,0,0)$, the extra relation (G2) is
\[
Y_{A_6}=0.
\]
Each block is a $1\times2$ matrix, so this is a pair of scalar
linear equations. Likewise, each bending relation is a pair of
scalar equations. Any two independent bending relations together
with $Y_{A_6}=0$ therefore give six independent scalar linear
generators.

Because the output dimension is one, each $\adjdiff{k}$ is a
$1\times2$ matrix. Hence the $2\times2$ minors in (G3) are
vacuous. It follows from \Cref{thm:shallow_generators} that these
six linear forms generate $J^{\mathcal{A}_\theta}$ and that
\[
U^{\mathcal{A}_\theta}\cong\rankonevar^3
=(\mathbb C^{1\times2})^3.
\]
Consequently,
\[
\dim U^{\mathcal{A}_{\theta}}
=3(2+1-1)=6.
\]
\end{example}

\section{Pattern equivalence and clock classification} \label{sec:permutations}
\subsection{Permutation actions and ReLU pattern varieties}\label{subsec:pattern_equivalence}
In this section, we introduce the notion of \emph{pattern equivalence} by allowing permutations of the neurons
within each hidden layer. For nondegenerate shallow networks with
two-dimensional input, we classify the resulting equivalence classes by
dihedral orbits of clocks and determine their number, with proofs given in Appendix~\ref{app:relu_clocks}. The relationship between pattern equivalence and equivalence of ReLU pattern varieties is studied in \Cref{sec:hierarchy}.

\begin{definition}\label{def:single_actionvationpattern_permutation}
Let
\[
p=(p_1,\ldots,p_{m-1})
\in S_{n_1}\times\cdots\times S_{n_{m-1}},
\]
where $S_{n_\ell}$ denotes the symmetric group on
$n_\ell$ elements.
For a binary activation pattern
\[
A=(A^{(1)},\ldots,A^{(m-1)}),
\qquad
A^{(\ell)}\in\{0,1\}^{n_\ell},
\]
define $p\cdot A$ by
\[
(p\cdot A)^{(\ell)}_k
=
A^{(\ell)}_{p_\ell^{-1}(k)},
\qquad
1\leq \ell\leq m-1,\quad
1\leq k\leq n_\ell,
\]
so that $p_\ell$ permutes the entries corresponding to the
neurons in the $\ell$-th hidden layer.
This defines a left action of
$S_{n_1}\times\cdots\times S_{n_{m-1}}$
on the set of all binary activation patterns.
\end{definition}

\begin{definition}\label{def:setofactivationpattern_permutation}
For a set $\mathcal A$ of binary activation patterns and
$p\in S_{n_1}\times\cdots\times S_{n_{m-1}}$, define
\[
p\cdot\mathcal A=\{p\cdot A:A\in\mathcal A\}.
\]
The same permutation is applied simultaneously to every pattern in
$\mathcal A$. In particular, this defines $p\cdot \mathcal{A}_{\theta}$ for every
$\theta\in\mathcal P_{n_0,\ldots,n_m}$.
\end{definition}

\begin{definition}[Pattern-equivalent parameters]\label{def:patternequivalentparameter}
Two parameters
$\theta_1,\theta_2\in\mathcal P_{n_0,\ldots,n_m}$
are \emph{pattern-equivalent}, written
$\theta_1\sim\theta_2$, if there exists
\[
p\in S_{n_1}\times\cdots\times S_{n_{m-1}}
\]
such that
\[
p\cdot \mathcal{A}_{\theta_1}=\mathcal{A}_{\theta_2}.
\]
\end{definition}

The actions in
\Cref{def:single_actionvationpattern_permutation,def:setofactivationpattern_permutation}
show that $\sim$ is an equivalence relation.

\begin{remark}
The canonical set $\mathcal{A}_{\theta}$ is unordered. Accordingly, the equalities
\[
\mathcal{A}_{\theta'}=\mathcal{A}_{\theta}
\qquad\text{and}\qquad
p\cdot \mathcal{A}_{\theta_1}=\mathcal{A}_{\theta_2}
\]
in \Cref{def:pattern_stable_neighborhood,def:patternequivalentparameter},
respectively, are equalities of sets and do not depend on any chosen
enumeration of the activation patterns.
\end{remark}

\subsection{\texorpdfstring{Clock classification and enumeration for the architecture $(2,n_1,n_2)$}{Clock classification and enumeration for the architecture (2,n1,n2)}} \label{subsec:clock_pattern_equivalence}

We consider unbiased shallow networks of architecture
$(2,n_1,n_2)$. Let $\mathcal P^{\mathrm{nd}}_{2,n_1,n_2}$ denote the set of parameters $\theta=(W_1,W_2)$ such that every row of $W_1$ is nonzero and every two distinct rows are linearly independent. Equivalently, the rows of $W_1$ are nonzero and no two rows are parallel. For every $\theta\in\mathcal P^{\mathrm{nd}}_{2,n_1,n_2}$, the canonical set of activation patterns $\mathcal{A}_{\theta}$ depends only on $W_1$. Consequently, pattern equivalence is determined entirely by $W_1$, and the number of pattern-equivalence classes is independent of the output dimension $n_2$.

These classes admit a description in terms of \emph{clocks}: cyclic sign sequences of length $2n_1$ in which antipodal entries have opposite signs (Definition~\ref{def:clock}). Write $\clockset{n_1}$ for the set of clocks of size $n_1$ and $D_{2n_1}$ for the dihedral group of order $4n_1$, acting on $\clockset{n_1}$ by rotations and reflections. The construction that associates a clock to a parameter is given in \Cref{app:clock_construction}. This construction involves a choice of one of two opposite rays as the initial ray. The two choices produce clocks that differ by a half-turn and form a \emph{shadow pair}. The dihedral action on clocks is defined in \Cref{app:clock_action}.

\begin{theorem}[Clock characterization of pattern equivalence] \label{thm:clock_pattern_equivalence}

Let
$\theta,\theta'\in\mathcal P^{\mathrm{nd}}_{2,n_1,n_2}$. Then $\theta$ and $\theta'$ are pattern-equivalent if and only if a clock associated with $\theta$ and a clock associated with $\theta'$ belong to the same orbit under $D_{2n_1}$. Since the two clocks in each shadow pair differ by a half-turn, this condition is independent of the choices of associated clocks. Consequently, the pattern-equivalence classes in
$\mathcal P^{\mathrm{nd}}_{2,n_1,n_2}$ are in bijection with
\[
\clockset{n_1}/D_{2n_1}.
\]
\end{theorem}
The proof is given in \Cref{app:clock_classification_proof}.

The resulting orbit-counting problem was studied by
Palmer and Robinson \cite[pp.~208--210]{PalmerRobinson1984}, who
enumerated self-dual two-colorings of necklaces with $2n_1$ beads
where opposite beads have different colors and colorings are identified
under rotations and reflections. Identifying the two colors with $+$
and $-$ gives precisely the clocks and equivalence relation above.
Their formula, written in our notation, yields the following count.

\begin{theorem}[Enumeration of pattern-equivalence classes] \label{thm:clock_counting_formula}
Let $N(n_1)$ be the number of pattern-equivalence classes in $\mathcal P^{\mathrm{nd}}_{2,n_1,n_2}$. Then
\[
N(n_1)=\frac{1}{4n_1}
\bigg(
\sum_{\substack{e\mid n_1\\ e\text{ odd}}}
\phi(e)\,2^{n_1/e}
+n_1\,2^{\lceil n_1/2\rceil}
\bigg),
\]
where $\phi$ denotes Euler's totient function.
\end{theorem}

For completeness, a direct fixed-point calculation in our notation is
given in \Cref{app:clock_counting_proof}; small values appear in
\Cref{tab:clock_counts}. The role of zero and repeated neurons is
discussed in \Cref{app:clock_degenerate}.

\section{Equivalence of pattern varieties and generic pattern stability}\label{sec:hierarchy}
In this section, we compare pattern equivalence with equivalence
of the associated ReLU pattern varieties and with local stability.
We distinguish the information retained by the canonical complex, the canonical set of activation patterns, and its ReLU pattern variety. Our first theorem compares pairs of parameters and separates three notions of equivalence for their varieties. The supporting results explain why these distinctions matter even for shallow networks. We then connect pattern equivalence to local stability and prove, in the second theorem, that pattern stability is generic at every depth.

\subsection{The hierarchy for pairs of parameters}\label{subsec:relations_parameters}
For shallow architectures, write $\mathcal P^{\mathrm{nd}}_{n_0,n_1,n_2}$ for the parameters whose first-layer rows are nonzero and pairwise linearly independent. We call $\theta,\theta'\in\mathcal P_{n_0,\ldots,n_m}$ \emph{combinatorially equivalent} if $\mathcal C(\theta)$ and $\mathcal C(\theta')$ are combinatorially equivalent by a bijection that respects the ternary labelings (\Cref{subsec:combinatorial_stability}). For a shallow network, we identify the single-layer pattern $(A^{(1)})$ with its tuple $A^{(1)}\in\{0,1\}^{n_1}$. For a deep network, we retain the layer structure $A=(A^{(1)},\ldots,A^{(m-1)})$.

When $|\mathcal{A}_{\theta}|=|\mathcal{A}_{\theta'}|=N$, choose enumerations of both sets and identify their ambient spaces with
\[
\mathbb C^{n_m\times Nn_0}
\cong(\mathbb C^{n_m\times n_0})^N.
\]
A \emph{block permutation} only reorders the $N$ matrices $Y_A$; it does not mix their entries. An arbitrary \emph{linear automorphism of the ambient space} may mix all coordinates. Finally, an \emph{isomorphism of affine varieties} means mutually inverse regular maps between the varieties themselves; it need not be specified by a map of their ambient spaces. Thus $U^{\mathcal{A}_{\theta'}}=\Pi(U^{\mathcal{A}_{\theta}})$ asserts equality of embedded sets after applying a specified map $\Pi$, whereas $U^{\mathcal{A}_{\theta'}}\cong U^{\mathcal{A}_{\theta}}$ asserts only abstract isomorphism.

The first link between pattern equivalence and equivalence of varieties comes from permuting the hidden neurons in the polynomial parametrization.

\begin{proposition}
\label{prop:pattern_equivalent_varieties}
Fix an architecture $(n_0,\ldots,n_m)$.
If $\theta_1,\theta_2\in\mathcal P_{n_0,\ldots,n_m}$
are pattern-equivalent, then
\[
|\mathcal{A}_{\theta_1}|=|\mathcal{A}_{\theta_2}|\eqqcolon N.
\]
For any chosen enumerations of these sets, there is
a linear automorphism $\Pi$ of
$\mathbb C^{n_m\times Nn_0}$ that permutes its
$N$ matrix blocks of size $n_m\times n_0$ and satisfies
\[
U^{\mathcal{A}_{\theta_2}}
=
\Pi\bigl(U^{\mathcal{A}_{\theta_1}}\bigr).
\]
In particular,
$U^{\mathcal{A}_{\theta_1}}\cong U^{\mathcal{A}_{\theta_2}}$.
If $p\cdot \mathcal{A}_{\theta_1}=\mathcal{A}_{\theta_2}$ and
$\mathcal{A}_{\theta_1}$ and $\mathcal{A}_{\theta_2}$ are enumerated as
\[
A_1,\ldots,A_N
\qquad\text{and}\qquad
p\cdot A_1,\ldots,p\cdot A_N,
\]
respectively, then $\Pi$ can be chosen to be the identity.
\end{proposition}

\begin{proof}
By Definition~\ref{def:patternequivalentparameter}, there exists
\[
p=(p_1,\ldots,p_{m-1})
\in S_{n_1}\times\cdots\times S_{n_{m-1}}
\]
such that $p\cdot \mathcal{A}_{\theta_1}=\mathcal{A}_{\theta_2}$.
Since the action of $p$ is bijective, the two
pattern sets have the same cardinality $N$.

For $1\leq i\leq m-1$, let $P_i$ be the
permutation matrix satisfying
\[
P_i e_k=e_{p_i(k)},
\]
and set $P_0=I_{n_0}$ and $P_m=I_{n_m}$.
Define the linear automorphism $T_p$ of the
parameter space by
\[
T_p(W_1,\ldots,W_m)
=
(W_1',\ldots,W_m'),
\qquad
W_i'=P_iW_iP_{i-1}^{-1}.
\]
For any binary activation pattern $A$, its diagonal
activation matrices satisfy
\[
D_i(p\cdot A)=P_iD_i(A)P_i^{-1}.
\]
Consequently, for every parameter
$\theta=(W_1,\ldots,W_m)$,
\[
\begin{aligned}
M_{p\cdot A}(T_p(\theta))
&=
W_m'D_{m-1}(p\cdot A)
W_{m-1}'\cdots D_1(p\cdot A)W_1'\\
&=
W_mD_{m-1}(A)
W_{m-1}\cdots D_1(A)W_1\\
&=M_A(\theta).
\end{aligned}
\]

Now enumerate $\mathcal{A}_{\theta_1}=\{A_1,\ldots,A_N\}.$ Since $p\cdot \mathcal{A}_{\theta_1}=\mathcal{A}_{\theta_2}$,
the activation patterns in $\mathcal{A}_{\theta_2}$ are precisely $p\cdot A_1,\ldots,p\cdot A_N$,
possibly in a different order.

It follows from
\[
M_{p\cdot A}(T_p(\theta))=M_A(\theta)
\]
that the polynomial maps
$\varphi^{\mathcal{A}_{\theta_1}}$ and
$\varphi^{\mathcal{A}_{\theta_2}}\circ T_p$
differ only by a permutation of their matrix blocks. Thus there exists a linear automorphism $\Pi$ of the ambient affine space, permuting these blocks, such that
\[
\varphi^{\mathcal{A}_{\theta_2}}\circ T_p
=
\Pi\circ\varphi^{\mathcal{A}_{\theta_1}},
\]
and $\Pi$ is the identity if $\mathcal{A}_{\theta_2}$ is enumerated as $p\cdot A_1,\ldots,p\cdot A_N$.

Since $T_p$ is an automorphism of the parameter space, we have
\[
\varphi^{\mathcal{A}_{\theta_2}}(\mathcal P_{n_0,\ldots,n_m})=\Pi\bigl(\varphi^{\mathcal{A}_{\theta_1}}(\mathcal P_{n_0,\ldots,n_m})\bigr),
\]
and since $\Pi$ is an automorphism of the ambient affine space, taking Zariski closures gives
\[
U^{\mathcal{A}_{\theta_2}}
=
\Pi\bigl(U^{\mathcal{A}_{\theta_1}}\bigr).
\]
\end{proof}

We now place this implication in the full hierarchy.

\equivalencehierarchy*

\begin{proof}
The implication (a)$\Rightarrow$(b) follows because a combinatorial equivalence preserves maximal-dimensional cells and their ternary labels, and a binary activation pattern is obtained from such a label by $1\mapsto1$ and $0,-1\mapsto0$. For (b)$\Rightarrow$(c), use the identity. The implication (c)$\Rightarrow$(d) is \Cref{prop:pattern_equivalent_varieties}. Every block permutation is a linear automorphism, giving (d)$\Rightarrow$(e), and the restriction of $T$ to the variety gives (e)$\Rightarrow$(f).

For the shallow converse (b)$\Rightarrow$(a), first suppose that the first-layer rows are nonzero and pairwise linearly independent. The binary canonical patterns determine the topes by replacing $0$ with $-1$. Mandel's theorem \cite[Theorem~4.2.13]{BLSWZ99} recovers the covectors from the topes. The covectors determine the ternary-labeled face poset and the ranks of the underlying matroid, hence the dimensions of its faces, exactly as in the proof of \Cref{thm:shallow_pattern_combinatorial_stable}.

The same conclusion holds in the degenerate case by reduction to this situation. 

A coordinate is zero in every canonical pattern if and only if its first-layer row is zero: a nonzero row is positive on an open half-space meeting a chamber. 

Among the nonzero rows, positive proportionality is equivalent to equality of the corresponding coordinates in every pattern in $\mathcal{A}_{\theta}$, while negative proportionality is equivalent to those coordinates summing to one in every pattern in $\mathcal{A}_{\theta}$. The forward implications are immediate. For the converses, suppose that two nonzero rows $w_i^\top$ and $w_j^\top$ are not proportional. Then the canonical pattern set $\mathcal{A}_{\theta}$ contains patterns whose $i$-th and $j$-th coordinates are,
respectively,
\[
(1,1),\qquad (1,0),\qquad (0,1),\qquad (0,0).
\]
Therefore these coordinates are neither always equal nor always complementary. 

Thus, the zero rows and all signed parallel classes are determined by $\mathcal{A}_{\theta}$. Delete the zero rows and retain one representative of each parallel class. The resulting topes determine the covectors of the reduced arrangement by Mandel's theorem; reinserting zero coordinates and the prescribed equal or opposite coordinates recovers all original covectors. If all rows vanish, the complex consists only of $\mathbb R^{n_0}$, with its all-zero ternary label. This proves the shallow assertion without a nondegeneracy hypothesis.

For (b)$\not\Rightarrow$(a), use the pattern-stable parameter $\theta\in\mathcal P_{3,4,1,1}$ and its perturbations $\theta_\varepsilon$ from \Cref{ex:pattern_combinatorial}. For all sufficiently small $\varepsilon>0$, pattern stability gives $\mathcal{A}_{\theta_\varepsilon}=\mathcal{A}_{\theta}$, while the proof asserts that $\theta$ and $\theta_{\varepsilon}$ are not combinatorially equivalent.

The remaining assertions are established below: \Cref{ex:equivalent_distinct_patterns} proves (ii), \Cref{prop:block_permutation_nd} and  \Cref{ex:block_classes_scalar} prove (iii), \Cref{cor:converse_fails} proves (iv), and \Cref{ex:abstract_different_ambient} proves (v).
\end{proof}

The distinction between (e) and (f) includes the ambient-space requirement in (e). The example for (f)$\not\Rightarrow$(e) uses different values of $N$; it does not assert failure of ambient linear equivalence when $N$ is fixed. In fact, for nondegenerate shallow parameters with the same $N$, condition (e) always holds by \Cref{prop:nd_single_class}.

\subsection{Shallow varieties: linear structure and counterexamples}\label{subsec:linear_facts}

For an enumerated pattern set $\mathcal A=\{A_1,\ldots,A_N\}$, let $\Lam_{\mathcal A}$ be the $N\times n_1$ matrix with $j$-th row $A_j$. The next three results separate the linear description for one-dimensional output, the rank supplied by nondegeneracy, and the extension to an automorphism of the ambient space.

\begin{proposition}\label{prop:scalar_output_linear}
Consider the shallow architecture $(n_0,n_1,1)$, and let $\mathcal A=\{A_1,\ldots,A_N\}\subseteq\{0,1\}^{n_1}$ be any set of distinct patterns. Identify $\mathbb C^{1\times Nn_0}$ with $\mathbb C^{N\times n_0}$ by arranging the blocks $Y_{A_j}$ as rows, and use $\mathbf Y$ for the resulting matrix of indeterminates. Then
\[
U^{\mathcal A}
=
\{\Lam_{\mathcal A}V:\ V\in\mathbb C^{n_1\times n_0}\}
\]
is a linear subspace of dimension
$n_0\operatorname{rank}\Lam_{\mathcal A}$. If the rows of
\[
K\in
\mathbb C^{(N-\operatorname{rank}\Lam_{\mathcal A})\times N}
\]
form a basis of the left kernel of $\Lam_{\mathcal A}$, then $J^{\mathcal A}$ is generated by the
$(N-\operatorname{rank}\Lam_{\mathcal A})n_0$ entries of $K\mathbf Y$. In particular, $J^{\mathcal A}$ is generated by linear forms.
\end{proposition}

\begin{proof}
For $\theta=(W_1,W_2)\in\mathcal P_{n_0,n_1,1}$, write
the rows of $W_1$ as $w_k^\top$ and
$W_2=(c_1,\ldots,c_{n_1})$. Let
$V(\theta)\in\mathbb R^{n_1\times n_0}$ have rows
$c_kw_k^\top$. Since
\[
M_A(\theta)=\sum_{k:A_k=1}c_kw_k^\top,
\]
we have
\[
\varphi^{\mathcal A}(\theta)
=
\Lam_{\mathcal A}V(\theta).
\]
For any $B\in\mathbb R^{n_1\times n_0}$, choosing
$W_1=B$ and $W_2=(1,\ldots,1)$ gives $V(\theta)=B$.
Hence the image of $\varphi^{\mathcal A}$ is the set of real points of
\[
L:=
\{\Lam_{\mathcal A}V:\ V\in\mathbb C^{n_1\times n_0}\}.
\]
These real points are Zariski dense in $L$, so
$U^{\mathcal A}=L$.

Writing
\[
V=[u_1\mid\cdots\mid u_{n_0}],
\qquad u_i\in\mathbb C^{n_1},
\]
we have
\[
\Lam_{\mathcal A}V
=
[\Lam_{\mathcal A}u_1\mid\cdots\mid
\Lam_{\mathcal A}u_{n_0}].
\]
Thus its columns vary independently over
$\operatorname{im}\Lam_{\mathcal A}$, and therefore
\[
\dim U^\mathcal{A}=\dim L =n_0\operatorname{rank}\Lam_{\mathcal A}.
\]

Finally, since $\ker K=\operatorname{im}\Lam_{\mathcal A}$, a matrix
$Y\in\mathbb C^{N\times n_0}$ lies in $L$ if and only if $KY=0$.
Therefore the ideal of $L$ is generated by the $(N-\operatorname{rank}\Lam_{\mathcal A})n_0$ entries of $K\mathbf Y$.
\end{proof}

For canonical sets, nondegeneracy supplies full column rank.

\begin{lemma}\label{lem:canonical_pattern_matrix_full_rank}
Let $\theta\in\mathcal P^{\mathrm{nd}}_{n_0,n_1,n_2}$
be a nondegenerate shallow parameter and $\mathcal{A}_{\theta}$ the canonical set of activation patterns. Then
\[
\operatorname{rank}\Lam_{\mathcal{A}_{\theta}}=n_1.
\]
\end{lemma}

\begin{proof}
Since the first-layer rows are nonzero and pairwise linearly independent, the hyperplanes $H_k=\{x:w_k^\top x=0\}$ are distinct. For each $k$, choose
\[
x\in H_k\setminus\bigcup_{\ell\neq k}H_\ell.
\]
Such a point exists because the other hyperplanes cut proper linear subspaces of $H_k$, and finitely many proper subspaces cannot cover a real vector space. When $n_0=1$, nondegeneracy forces $n_1=1$, so the union is empty. Crossing $H_k$ near $x$ changes only the $k$-th activation coordinate. The canonical patterns of the two
adjacent maximal cells therefore differ only in their $k$-th coordinate. Hence the row span of $\Lam_{\mathcal{A}_{\theta}}$ contains
$e_k$ for every $k$, so its rank is at least $n_1$. Since
$\Lam_{\mathcal{A}_{\theta}}$ has $n_1$ columns, its rank is exactly $n_1$.
\end{proof}

Combining \Cref{prop:scalar_output_linear} and \Cref{lem:canonical_pattern_matrix_full_rank} gives
\[
\dim U^{\mathcal{A}_{\theta}}=n_0n_1
\qquad\text{for }\theta\in\mathcal P^{\mathrm{nd}}_{n_0,n_1,1}.
\]
For these parameters, the $2\times2$ minors in \Cref{thm:shallow_generators}~(i) are vacuous, so the whole ideal is generated by the linear generators. More generally, the linear relations in \Cref{prop:scalar_output_linear} are the zero-support case of \cite[eq.~(5)]{alexandr2025constraining}: the integer matrix $\Lam_{\mathcal A}$ has a left-kernel basis consisting of integer vectors, and
\[
\sum_j\lambda_j A_j=0
\quad\Longrightarrow\quad
\operatorname{rank}\Bigl(\sum_j\lambda_jM_{A_j}\Bigr)\leq0.
\]
Thus $\sum_j\lambda_jY_{A_j}=0$. The proposition proves that these relations generate the ideal when the output is one-dimensional.

For arbitrary output dimension, \Cref{thm:shallow_generators}~(ii) identifies the varieties of nondegenerate shallow parameters with the same product $\rankonevar^{n_1}$. Equal numbers of blocks allow these identifications to extend to the ambient spaces.

\begin{proposition}\label{prop:nd_single_class}\label{cor:nd_single_class}
Let $\theta,\theta' \in\mathcal P^{\mathrm{nd}}_{n_0,n_1,n_2}$ and suppose that $|\mathcal{A}_{\theta}|=|\mathcal{A}_{\theta'}|\eqqcolon N$. Then there is a linear automorphism
\[
T:\mathbb C^{n_2\times Nn_0}
\longrightarrow
\mathbb C^{n_2\times Nn_0}
\]
such that
\[
T(U^{\mathcal{A}_{\theta}})=U^{\mathcal{A}_{\theta'}}.
\]
In particular, for $n_0=2$, all varieties
$U^{\mathcal{A}_{\theta}}$ with $\theta\in\mathcal P^{\mathrm{nd}}_{2,n_1,n_2}$ are linearly equivalent in their common ambient space.
\end{proposition}

\begin{proof}
By \Cref{thm:shallow_generators}~(ii),
$\slopeembed{\mathcal{A}_{\theta}}$ and $\slopeembed{\mathcal{A}_{\theta'}}$ are injective and
\[
U^{\mathcal{A}_{\theta}}
=
\slopeembed{\mathcal{A}_{\theta}}(\rankonevar^{n_1}),
\qquad
U^{\mathcal{A}_{\theta'}}
=
\slopeembed{\mathcal{A}_{\theta'}}(\rankonevar^{n_1}).
\]
Consequently,
\[
g:=
\slopeembed{\mathcal{A}_{\theta'}}
\circ
\slopeembed{\mathcal{A}_{\theta}}^{-1}
\]
is a linear isomorphism from
\[
E:=\operatorname{im}\slopeembed{\mathcal{A}_{\theta}}
\quad\text{onto}\quad
E':=\operatorname{im}\slopeembed{\mathcal{A}_{\theta'}}
\]
that sends $U^{\mathcal{A}_{\theta}}$ onto $U^{\mathcal{A}_{\theta'}}$.
Choose complements $F,F'$ such that
\[
\mathbb C^{n_2\times Nn_0}
=
E\oplus F
=
E'\oplus F'.
\]
Since $\dim E=\dim E'$, we have $\dim F=\dim F'$; hence
there is a linear isomorphism $h:F\to F'$. Define
\[
T(e+f):=g(e)+h(f),
\qquad e\in E,\quad f\in F.
\]
Then $T$ is a linear automorphism of the ambient space extending $g$, and therefore
\[
T(U^{\mathcal{A}_{\theta}})=U^{\mathcal{A}_{\theta'}}.
\]

For $n_0=2$, every
$\theta\in\mathcal P^{\mathrm{nd}}_{2,n_1,n_2}$ satisfies $|\mathcal{A}_{\theta}|=2n_1$ because $n_1$ distinct lines through the origin cut the plane into $2n_1$ regions, so the first
statement applies.
\end{proof}

We now restrict the ambient maps to block permutations.

\begin{proposition}\label{prop:block_permutation_nd}
Let $n_0,n_2\geq2$ and $\theta,\theta'\in\mathcal P^{\mathrm{nd}}_{n_0,n_1,n_2}$ with $|\mathcal{A}_{\theta}|=|\mathcal{A}_{\theta'}|=N$. There is a linear automorphism $\Pi$ of $\mathbb C^{n_2\times Nn_0}$ permuting the $N$ blocks with $\Pi(U^{\mathcal{A}_{\theta}})=U^{\mathcal{A}_{\theta'}}$ if and only if $\theta\sim\theta'$.
\end{proposition}

\begin{proof}
The ``if'' direction follows from
\Cref{prop:pattern_equivalent_varieties}.
Conversely, suppose that
$\Pi(U^{\mathcal{A}_{\theta}})=U^{\mathcal{A}_{\theta'}}$ for a block permutation $\Pi$, and let $P$ be its $N\times N$ permutation matrix. Write
\[
\Lam=\Lam_{\mathcal{A}_{\theta}},
\qquad
\Lam'=\Lam_{\mathcal{A}_{\theta'}},
\]
and let $\rankonevar\subseteq \mathbb C^{n_2\times n_0}$ be the variety of matrices of rank at most one.
For $B\in\mathbb C^{N\times n_1}$, define the linear map
\[
S_B:(\mathbb C^{n_2\times n_0})^{n_1}\longrightarrow (\mathbb C^{n_2\times n_0})^N,
\qquad
S_B(V_1,\ldots,V_{n_1})
=
\left(\sum_{k=1}^{n_1}B_{jk}V_k\right)_{j=1}^N.
\]
Thus $S_\Lam=\slopeembed{\mathcal{A}_{\theta}}$ and
$S_{\Lam'}=\slopeembed{\mathcal{A}_{\theta'}}$. By
\Cref{thm:shallow_generators}~(ii), our assumption becomes
\[
S_{P\Lam}(\rankonevar^{n_1})
=
S_{\Lam'}(\rankonevar^{n_1}).
\]

Since the matrix units have rank one, $\rankonevar$ spans $\mathbb C^{n_2\times n_0}$, and hence $\rankonevar^{n_1}$ spans $(\mathbb C^{n_2\times n_0})^{n_1}$. Taking linear spans therefore gives
\[
S_{P\Lam}((\mathbb C^{n_2\times n_0})^{n_1})=S_{\Lam'}((\mathbb C^{n_2\times n_0})^{n_1}).
\]
For any fixed matrix-entry position $(a,b)$, the vector of $(a,b)$-entries across the $N$ blocks of $S_B(V)$ ranges over the column space of $B$. Consequently, $P\Lam$ and $\Lam'$ have the same column space. Both have rank $n_1$ by \Cref{lem:canonical_pattern_matrix_full_rank}. Their columns are thus two bases of the same space, so
\[
\Lam'=P\Lam g
\]
for some $g\in\mathrm{GL}_{n_1}(\mathbb C)$.

Define
\[
G_g:(\mathbb C^{n_2\times n_0})^{n_1}\longrightarrow (\mathbb C^{n_2\times n_0})^{n_1},
\qquad
G_g(V_1,\ldots,V_{n_1})
=
\left(\sum_{k=1}^{n_1}g_{jk}V_k\right)_{j=1}^{n_1}.
\]
The identity $\Lam'=P\Lam g$ implies
$S_{\Lam'}=S_{P\Lam}\circ G_g$, and therefore
\[
S_{P\Lam}(\rankonevar^{n_1})
=
S_{P\Lam}\bigl(G_g(\rankonevar^{n_1})\bigr).
\]
By injectivity of $S_{P\Lam}$, we obtain
\[
G_g(\rankonevar^{n_1})=\rankonevar^{n_1}.
\]
Thus every output matrix of $G_g$ must have rank at most one whenever all its input matrices do.

We claim that each row of $g$ has at most one nonzero
entry. Otherwise, for some row $j$, choose $k\neq\ell$ with $g_{jk}\neq0$ and $g_{j\ell}\neq0$, and set
\[
V_k=E_{11},\qquad V_\ell=E_{22},
\qquad V_i=0\quad(i\neq k,\ell),
\]
where $E_{ab}\in\mathbb C^{n_2\times n_0}$ has a single
$1$ in position $(a,b)$ and zeros elsewhere.
Then $V\in\rankonevar^{n_1}$, but the $j$-th component
of $G_g(V)$ is
\[
g_{jk}E_{11}+g_{j\ell}E_{22}.
\]
Since $n_0,n_2\geq2$, its upper-left $2\times2$ minor
is defined and equals $g_{jk}g_{j\ell}\neq0$.
Thus this component has rank two, so $G_g(V)\notin\rankonevar^{n_1}$, contradicting
$G_g(\rankonevar^{n_1})=\rankonevar^{n_1}$.

Since $g$ is invertible, it consequently has exactly one nonzero entry in each row and each column. Hence
$g=QD$, where $Q$ is a permutation matrix and $D$ is an invertible diagonal matrix. Now
\[
\Lam'=P\Lam QD.
\]
Every column of the binary matrix $P\Lam Q$ is nonzero and therefore contains a $1$. At such an entry in column $k$, the corresponding entry of $\Lam'$ equals $D_{kk}$. Since $\Lam'$ is binary and $D_{kk}\neq0$, we have $D_{kk}=1$ for every $k$. Thus $D=I$ and
\[
\Lam'=P\Lam Q.
\]
The matrix $P$ only reorders the patterns, while $Q$
permutes their neuron coordinates. Therefore
$\mathcal{A}_{\theta'}=\{AQ:A\in \mathcal{A}_{\theta}\}$, which is precisely $\theta\sim\theta'$.
\end{proof}

\begin{corollary}[Classification for two-dimensional input]
\label{cor:clock_variety_classes}
Fix $n_1\geq1$ and $n_2\geq2$. Among the varieties $U^{\mathcal{A}_{\theta}}$ with
$\theta\in\mathcal P^{\mathrm{nd}}_{2,n_1,n_2}$, there are exactly
$N(n_1)$ equivalence classes under block permutations.
These varieties nevertheless
form a single equivalence class under linear automorphisms of their
common ambient space.
\end{corollary}

\begin{proof}
By \Cref{prop:block_permutation_nd}, equivalence under block permutations
is equivalent to pattern equivalence. The count follows from
\Cref{thm:clock_counting_formula}, and the ambient linear equivalence follows
from \Cref{prop:nd_single_class}.
\end{proof}

\begin{example}[Pattern equivalence without equality]\label{ex:equivalent_distinct_patterns}
Fix any $n_2\geq1$ and consider $\theta, \theta' \in \mathcal P^{\mathrm{nd}}_{2,3,n_2}$ where the first-layer matrices of $\theta$ and $\theta'$ are respectively
\[
W_1=\begin{pmatrix}0&1\\-1&1\\-1&0\end{pmatrix},
\qquad
W_1'=\begin{pmatrix}0&-1\\-1&1\\-1&0\end{pmatrix}.
\]
The output matrices may be arbitrary. Reading the patterns on the six regions gives
\begin{align*}
\mathcal{A}_{\theta}&=\{(0,0,0),(0,0,1),(0,1,1),(1,0,0),(1,1,0),(1,1,1)\},\\
\mathcal{A}_{\theta'}&=\{(0,0,0),(0,1,0),(0,1,1),(1,0,0),(1,0,1),(1,1,1)\}.
\end{align*}
Swapping the second and third coordinates maps $\mathcal{A}_{\theta}$ onto $\mathcal{A}_{\theta'}$, so $\theta\sim\theta'$. Nevertheless, $(0,0,1)\in \mathcal{A}_{\theta}\setminus \mathcal{A}_{\theta'}$.
The first-layer matrices themselves need not be obtained from one another by this neuron permutation: pattern equivalence requires only the corresponding relation between their canonical pattern sets.
\end{example}

The hypothesis $n_2\geq2$ in \Cref{prop:block_permutation_nd} is essential even for nondegenerate parameters.

\begin{example}[Failure for one-dimensional output]
\label{ex:block_classes_scalar}
When $n_2=1$, every matrix in
$\mathbb C^{n_2\times n_0}$ has rank at most one, so
preservation of $\rankonevar^{n_1}$ imposes no additional
restriction on $g\in\mathrm{GL}_{n_1}(\mathbb C)$.
We construct two nondegenerate parameters whose varieties
agree up to a block permutation although the parameters
are not pattern-equivalent.

Let $\theta,\theta'\in\mathcal P^{\mathrm{nd}}_{2,5,1}$
have first-layer rows
\[
(0,1),\ (-1,2),\ (-1,1),\ (1,0),\ (-1,-1)
\]
and
\[
(0,1),\ (-1,2),\ (1,-1),\ (-1,0),\ (-1,-1),
\]
respectively. Their canonical pattern sets consist of the rows of $B$ and $B'$, respectively, together with the patterns obtained from these rows by replacing every $0$ with $1$ and every $1$ with $0$, where
\[
B=
\begin{pmatrix}
0&0&0&0&1\\
0&0&0&1&0\\
0&0&1&0&1\\
0&1&1&0&1\\
1&1&1&0&0
\end{pmatrix},
\qquad
B'=
\begin{pmatrix}
0&0&0&1&1\\
0&0&1&0&0\\
0&1&0&1&1\\
1&1&0&1&0\\
1&1&0&0&0
\end{pmatrix}.
\]
These patterns can be read off from the first-layer
arrangements, or from the clocks of
\Cref{subsec:clock_pattern_equivalence}.

Both matrices have determinant $1$. The linear map sending the $i$-th row of $B$ to the $i$-th row of $B'$ is therefore given by right multiplication by $g=B^{-1}B'\in\mathrm{GL}_5(\mathbb Z)$.
Explicitly, for a pattern $A=(A_1,A_2,A_3,A_4,A_5)$, 
\[
Ag=(A_2,A_3,A_4,A_5,A_1-A_2+A_5).
\]
Writing $\mathbf1=(1,1,1,1,1)$, we have
$\mathbf1g=\mathbf1$, and hence
\[
(\mathbf1-A)g=\mathbf1-Ag.
\]
Thus this map also matches the complementary patterns, so it maps $\mathcal{A}_{\theta}$ bijectively onto $\mathcal{A}_{\theta'}$.
Consequently,
\[
\Lam_{\mathcal{A}_{\theta'}}=P\Lam_{\mathcal{A}_{\theta}}g
\]
for a permutation matrix $P$. Since $g$ is invertible,
\Cref{prop:scalar_output_linear} gives
\[
U^{\mathcal{A}_{\theta'}}=\Pi(U^{\mathcal{A}_{\theta}}),
\]
where $\Pi$ is the block permutation induced by $P$.

However, $\mathcal{A}_{\theta}$ contains two patterns with exactly one active coordinate, whereas $\mathcal{A}_{\theta'}$ contains only one. This number is invariant under permutations
as in \Cref{def:setofactivationpattern_permutation}, so $\theta\not\sim\theta'$.
\end{example}

\begin{example}[Linear equivalence without a block permutation]\label{cor:converse_fails}
Fix $n_2\geq 2$ and let $\theta, \theta' \in \mathcal P^{\mathrm{nd}}_{2,3,n_2}$ with first-layer matrices
\[
W_1=\begin{pmatrix}0&1\\-1&1\\-1&0\end{pmatrix},
\qquad
W_1'=\begin{pmatrix}0&1\\1&-1\\-1&0\end{pmatrix}.
\]
The output matrices may be arbitrary. Their pattern sets are
\begin{align*}
\mathcal{A}_{\theta}&=\{(0,0,0),(0,0,1),(0,1,1),(1,0,0),(1,1,0),(1,1,1)\},\\
\mathcal{A}_{\theta'}&=\{(0,0,1),(0,1,0),(0,1,1),(1,0,0),(1,0,1),(1,1,0)\}.
\end{align*}
The all-zero pattern belongs only to $\mathcal{A}_{\theta}$, so $\theta \not \sim \theta'$. They represent the two classes counted by $N(3)=2$ in \Cref{thm:clock_counting_formula}. By \Cref{prop:nd_single_class}, a linear automorphism of $\mathbb C^{n_2\times12}$ carries one variety onto the other. 
However, since $n_2\geq2$ and
$\theta\not\sim\theta'$, \Cref{prop:block_permutation_nd}
shows that no block permutation carries one variety onto
the other.
\end{example}

\begin{example}[Abstract isomorphism with different ambient spaces]\label{ex:abstract_different_ambient}
Fix $n_2\geq1$ and choose shallow parameters of architecture $(3,3,n_2)$ with
\[
W_1=\begin{pmatrix}1&0&0\\0&1&0\\1&1&0\end{pmatrix},
\qquad W_1'=I_3.
\]
Both first layers are nondegenerate in the pairwise sense. The first arrangement has six maximal cells, obtained from the six regions in the $(x_1,x_2)$-plane by taking their product with $\mathbb R$; the second has eight regions. Nevertheless, \Cref{thm:shallow_generators}~(ii) gives
\[
U^{\mathcal{A}_{\theta}}\cong\rankonevar^3\cong U^{\mathcal{A}_{\theta'}},
\qquad \rankonevar\subseteq\mathbb C^{n_2\times3}.
\]
Their prescribed ambient spaces are $\mathbb C^{n_2\times18}$ and $\mathbb C^{n_2\times24}$, so they cannot satisfy condition (e) of \Cref{thm:hierarchy_parameters}. This example distinguishes abstract isomorphism from linear equivalence in a common ambient space through the different numbers of blocks.
\end{example}

\subsection{From pattern equivalence to generic stability}\label{subsec:local_properties}

The preceding hierarchy compares two parameters. Pattern stability instead asks whether the exact canonical set remains constant near a single parameter. Each pattern-equivalence class is a union of fibers of the map $\theta\mapsto \mathcal{A}_{\theta}$, so a pattern-stable parameter lies in the interior of its pattern-equivalence class. The following result proves the converse for shallow networks.

\begin{proposition}\label{prop:local_pattern_equivalence}
A shallow parameter $\theta\in\mathcal P_{n_0,n_1,n_2}$ is pattern-stable if and only if it has a neighborhood all of whose elements are pattern-equivalent to $\theta$.
\end{proposition}

In other words, if all parameters in a neighborhood of
$\theta$ have canonical pattern sets equivalent to $\mathcal{A}_{\theta}$, then, after shrinking the neighborhood, all these pattern sets are equal to $\mathcal{A}_{\theta}$.

\begin{proof}
If $\theta$ is pattern-stable, then $\mathcal{A}_{\theta'}=\mathcal{A}_{\theta}$
for all $\theta'$ in some neighborhood of $\theta$,
so every such $\theta'$ is pattern-equivalent to $\theta$.

For the converse, suppose that $\theta$ is not
pattern-stable. We show that arbitrarily close to $\theta$ there are parameters that are not pattern-equivalent to it. For a set $\mathcal A\subseteq\{0,1\}^{n_1}$ of activation patterns, put
\[
\nu(\mathcal A):=\sum_{A\in\mathcal A}\sum_{k=1}^{n_1}A_k.
\]
Thus $\nu(\mathcal A)$ counts the total number of active
coordinates across all patterns in $\mathcal A$. Both $|\mathcal A|$ and $\nu(\mathcal A)$ are invariant
under permutations as in \Cref{def:setofactivationpattern_permutation}.

Suppose first that $W_1$ has no zero row.
By \Cref{prop:pattern_inclusive_nbhd_exist}, there is
a neighborhood $\mathcal N$ of $\theta$ such that $\mathcal{A}_{\theta}\subseteq \mathcal{A}_{\theta'}$ for every $\theta'\in\mathcal N.$
Since $\theta$ is not pattern-stable, there are
$\theta'\in\mathcal N$ arbitrarily close to $\theta$
with $\mathcal{A}_{\theta'}\neq \mathcal{A}_{\theta}$. For these parameters,
the inclusion is strict, and hence
\[
|\mathcal{A}_{\theta'}|>|\mathcal{A}_{\theta}|.
\]

Now suppose that $W_1$ has a zero row $w_k^\top$.
Fix $v\in\mathbb R^{n_0}\setminus\{0\}$ and let
$\theta_\varepsilon$ be obtained by replacing this row
with $\varepsilon v^\top$ for $\varepsilon>0$ and leaving all other parameters
unchanged.
The first-layer arrangement of $\theta_\varepsilon$
is obtained by adding the hyperplane
$\{x:v^\top x=0\}$ to that of $\theta$.
Consequently, each old maximal cell is a union of new
maximal cells, and the new patterns agree with the
corresponding old pattern outside coordinate $k$.
Since distinct maximal cells of a shallow network have
distinct canonical patterns, it follows that
\[
|\mathcal{A}_{\theta_\varepsilon}|\geq |\mathcal{A}_{\theta}|.
\]

If this inequality is strict, the pattern sets cannot
be equivalent. If equality holds, each old maximal cell
corresponds to exactly one new maximal cell.
Under this correspondence, all coordinates other than
$k$ remain unchanged, while coordinate $k$, previously
always zero, becomes one on at least one cell:
indeed, the open half-space $\{x:v^\top x>0\}$ contains
a point outside the finitely many hyperplanes defined
by the nonzero rows of $W_1$. Therefore
\[
\nu(\mathcal{A}_{\theta_\varepsilon})>\nu(\mathcal{A}_{\theta}).
\]
Thus $\theta_\varepsilon$ is not pattern-equivalent
to $\theta$ in either case. Since
$\theta_\varepsilon\to\theta$ as $\varepsilon\to0$,
the converse follows.
\end{proof}

For $n_0=2$, \Cref{thm:shallow_pattern_combinatorial_stable} implies that every nondegenerate parameter is pattern-stable. Hence the $N(n_1)$ pattern-equivalence classes in $\mathcal P^{\mathrm{nd}}_{2,n_1,n_2}$ are nonempty open subsets of the parameter space. In higher input dimensions, pairwise independence alone is insufficient for stability. For example, the rows $e_1,e_2,e_1+e_2$ in $\mathbb R^3$ are pairwise linearly independent, but replacing the last row with $(1,1,\varepsilon)$, $\varepsilon\neq0$, creates the previously absent patterns $(1,1,0)$ and $(0,0,1)$.

Recall also the implications
\begin{equation}\label{eq:hierarchy_local}
\text{combinatorially stable}
\ \Longrightarrow\ \text{pattern-stable}
\ \Longrightarrow\ \text{pattern-inclusive}.
\end{equation}
The first is proved right above \Cref{ex:pattern_combinatorial}, and the second follows from the definitions. The first converse holds for shallow networks by \Cref{thm:shallow_pattern_combinatorial_stable} and fails for architecture $(3,4,1,1)$ by \Cref{ex:pattern_combinatorial}. The second converse fails already for architecture $(2,2,1)$: for the first-layer matrix, rows $(1,0)$ and $(2,0)$ give a pattern-inclusive parameter by \Cref{prop:pattern_inclusive_nbhd_exist}, but replacing the second row with $(2,\varepsilon)$ changes its canonical set from $\{(0,0),(1,1)\}$ to all four patterns for every $\varepsilon\neq0$.

We now prove that pattern stability holds on an open dense set of full Lebesgue measure for every depth. The argument uses the finiteness of the pattern space and the semialgebraic nature of pattern realization.

\begin{theorem}[Generic pattern stability]\label{prop:pattern_stable_generic}\label{thm:pattern_stable_generic}
For every architecture $(n_0,\ldots,n_m)$, the pattern-stable parameters form an open dense subset of $\mathcal P_{n_0,\ldots,n_m}$ whose complement has Lebesgue measure zero.
\end{theorem}

\begin{proof}
Identify the parameter space with $\mathbb R^D$, where $D=\sum_{\ell=1}^m n_\ell n_{\ell-1}$. Let
\[
\mathcal A_{\mathrm{all}}:=\prod_{\ell=1}^{m-1}\{0,1\}^{n_\ell}
\]
be the finite set of all possible activation patterns.

\emph{Step 1: realization on a maximal cell is a semialgebraic condition.}
For a fixed $A=(A^{(1)},\ldots,A^{(m-1)})\in\mathcal A_{\mathrm{all}}$, put
\[
L_1^A(\theta):=W_1,\qquad
L_\ell^A(\theta):=W_\ell D_{\ell-1}(A)W_{\ell-1}\cdots D_1(A)W_1
\quad(2\leq\ell\leq m-1).
\]
These are the preactivation matrices from the preliminaries, with all earlier-layer patterns fixed. By induction through the layers, $\alpha_{\theta}(x)=A$ if and only if
\[
(L_\ell^A(\theta)x)_i>0\quad\text{when }A_i^{(\ell)}=1,
\qquad
(L_\ell^A(\theta)x)_i\leq0\quad\text{when }A_i^{(\ell)}=0
\]
for every hidden neuron $(\ell,i)$. Denote this finite conjunction of polynomial inequalities in $(\theta,x)$ by $\Psi_A(\theta,x)$.

Moreover,
\[
A\in \mathcal{A}_{\theta}\quad\Longleftrightarrow\quad
R_A(\theta)\text{ has nonempty interior}.
\]
Indeed, the interior of a maximal-dimensional cell with pattern $A$ lies in $R_A(\theta)$. Conversely, a nonempty open subset of $R_A(\theta)$ cannot be contained in the finite union of lower-dimensional cells of the canonical complex, so it meets the interior of a maximal-dimensional cell, whose pattern must be $A$.
Thus $A\in \mathcal{A}_{\theta}$ is expressed by the first-order formula
\[
\exists x\in\mathbb R^{n_0}\ \exists r>0\ \forall z\in\mathbb R^{n_0}:
\quad \|z-x\|^2<r^2\ \Longrightarrow\ \Psi_A(\theta,z).
\]
The Tarski--Seidenberg theorem in its first-order form
\cite[Theorem~2.6]{Coste2002Semialgebraic} therefore implies that
$S_A:=\{\theta:A\in \mathcal{A}_{\theta}\}$ is semialgebraic.

\emph{Step 2: the finitely many fibers describe the stable set.}
For each $\mathcal A\subseteq\mathcal A_{\mathrm{all}}$, its fiber is
\[
F_{\mathcal A}:=\{\theta:\mathcal{A}_{\theta}=\mathcal A\}
=\bigcap_{A\in\mathcal A}S_A\ \cap\!
\bigcap_{A\in\mathcal A_{\mathrm{all}}\setminus\mathcal A}(\mathbb R^D\setminus S_A).
\]
These are finitely many semialgebraic sets, and the nonempty ones partition the parameter space. By definition, the pattern-stable set is exactly
\[
\bigcup_{\mathcal A\subseteq\mathcal A_{\mathrm{all}}}\operatorname{int}F_{\mathcal A},
\]
so it is open.

\emph{Step 3: all unstable parameters lie in a polynomial zero set.}
Since the sets $S_A$ are semialgebraic and there are only finitely many of them, there exist finitely many nonzero polynomials $p_1,\ldots,p_q$ such that the signs of these polynomials determine, for every activation pattern $A$, whether $A\in \mathcal{A}_{\theta}$. Collect the nonzero polynomials appearing in these descriptions for all $A\in\mathcal A_{\mathrm{all}}$, and let $P$ be their product, taking $P=1$ if the list is empty.

If $P(\theta)\neq0$, continuity ensures that all these polynomials keep their signs on a sufficiently small neighborhood of $\theta$. Thus, every condition defining the sets $S_A$ remains unchanged there, so $\mathcal{A}_{\theta'}=\mathcal{A}_{\theta}$ for every $\theta'$ in this neighborhood. Hence $\theta$ is pattern-stable, and all unstable parameters lie in
\[
Z:=\{\theta\in\mathbb R^D:P(\theta)=0\}.
\]
The polynomial $P$ is nonzero, so its zero set has Lebesgue measure zero and empty interior. Therefore $\mathbb R^D\setminus Z$ is dense and consists entirely of pattern-stable parameters. Together with Step~2, this proves that the pattern-stable set is open, dense, and of full Lebesgue measure.
\end{proof}

For shallow networks, \Cref{thm:shallow_pattern_combinatorial_stable}
identifies combinatorial stability with pattern stability. Consequently, the combinatorially stable shallow parameters form an open dense set of full Lebesgue measure. At greater depth, however, combinatorial stability is strictly stronger than pattern stability.

The deep setting also leaves open a local question concerning pattern equivalence. Pattern stability implies that a parameter is an interior point of its pattern-equivalence class: throughout a sufficiently small neighborhood, the canonical pattern set remains unchanged and hence, in particular, remains equivalent up to a permutation. For shallow networks, the converse holds. It is not known whether it continues to hold at greater depth.

\begin{problem}\label{prob:local_pattern_equivalence}
Let $\theta\in\mathcal P_{n_0,\ldots,n_m}$ with $m\geq3$, and suppose that every parameter in some neighborhood of $\theta$ is pattern-equivalent to $\theta$. Must $\theta$ be pattern-stable?
\end{problem}

A second question is whether the genericity theorem for pattern stability can be strengthened to combinatorial stability. Grigsby, Lindsey, Meyerhoff, and Wu conjecture full-measure combinatorial stability for networks with biases \cite[\S5.2 and Remark~4.16]{GLMWFunctionalDimension}. Their conjecture does not directly resolve the corresponding question here: our parameter space is unbiased, and the zero-bias subspace has measure zero in the biased parameter space. Moreover, our deep result establishes generic pattern stability, while the two stability notions need not coincide at greater depth.

\begin{problem}\label{prob:deep_genericity}
For $m\geq3$, do the combinatorially stable parameters of $\mathcal P_{n_0,\ldots,n_m}$ form a dense subset, or a subset of full Lebesgue measure?
\end{problem}

\section{Functional dimension and pattern varieties}
\label{sec:functional_dimension}

The variety $U^{\mathcal{A}_{\theta}}$ records the polynomial relations among the linear pieces indexed by the canonical patterns of $\theta$. We compare its dimension with the number of independent first-order variations detected by finitely many evaluations at inputs where every hidden preactivation is nonzero. These evaluations recover exactly the differential of the polynomial parametrization restricted to the patterns accessible to such inputs. Its rank is bounded by the rank of the full canonical pattern parametrization, and then by the dimension of the pattern variety. For shallow networks, all three quantities agree on an open dense set of full Lebesgue measure.

\subsection{Evaluation maps and the sampling convention}
\label{subsec:functional_evaluations}

Throughout this section, parameters are real and all derivatives are taken with respect to the weights. Variety dimensions are over $\mathbb C$. As in the setup, write
\[
D:=\sum_{\ell=1}^{m}n_\ell n_{\ell-1},
\qquad \mathcal P_{n_0,\ldots,n_m}\cong\mathbb R^D.
\]
Our evaluation maps follow the construction of
\cite[Definition~4.6]{GLMWFunctionalDimension}, with their output vector arranged as columns of a matrix. There are three differences in our conventions: our networks have no biases, the output layer is not ReLU activated, and we restrict the sample points to those at which every hidden preactivation is nonzero. In contrast, \cite[Definitions~2.1, 3.9, and~4.12]{GLMWFunctionalDimension} uses biases, a ReLU output layer, and all parametrically smooth sample points. Nonvanishing preactivations are sufficient for parametric smoothness, but are not necessary; see \cite[Corollary~3.6]{GLMWFunctionalDimension}.
Thus the definition below is a restricted unbiased analogue of their functional dimension. We also extend it by zero when there are no admissible sample points, whereas their definition is stated for ordinary parameters, namely those admitting a parametrically smooth point \cite[Definitions~3.12 and~4.12]{GLMWFunctionalDimension}.

\begin{definition}[Unbiased functional dimension on strict activation regions]
\label{def:unbiased_fundim}
Let $\theta\in\mathcal P_{n_0,\ldots,n_m}$.
A point $z\in\mathbb R^{n_0}$ is \emph{$\theta$-generic} if the hidden preactivation 
\[
h^{(\ell)}_{\theta,i}(z)\neq 0 \quad\text{for every }1\leq\ell\leq m-1,
\quad 1\leq i\leq n_\ell.
\]
For a finite ordered set $Z=(z_1,\ldots,z_k)$ of such points, define
\[
E_Z(\theta'):=
\bigl[f_{\theta'}(z_1)\mid\cdots\mid f_{\theta'}(z_k)\bigr]
\in\mathbb R^{n_m\times k}.
\]
Let $JE_Z(\theta)$ be its Jacobian, after arranging its output entries as a vector. The \emph{unbiased functional dimension on strict activation regions} is
\[
\dimfunzero(\theta)\coloneqq\max_Z\operatorname{rank}JE_Z(\theta).
\]
We allow $Z=\varnothing$, with rank zero. In particular,
$\dimfunzero(\theta)=0$ if no $\theta$-generic point exists.
\end{definition}

The next lemma ensures that these Jacobians exist. Their ranks belong to $\{0,\ldots,D\}$, so the maximum in the definition is attained.

\begin{lemma}[Local polynomial evaluation]\label{lem:evaluation_local}
If $z$ is $\theta$-generic and $A=\alpha_\theta(z)$, then $A\in \mathcal{A}_{\theta}$.
There are neighborhoods $\mathcal O$ of $\theta$ and $V$ of $z$ such that
\[
\alpha_{\theta'}(x)=A,
\qquad f_{\theta'}(x)=M_A(\theta')x
\quad\text{for all }(\theta',x)\in\mathcal O\times V.
\]
\end{lemma}

\begin{proof}
Each hidden preactivation $h^{(\ell)}_{\theta,i}(x)$ is jointly continuous in $(\theta,x)$. Since $z$ is $\theta$-generic, all of these finitely many functions are nonzero at $(\theta,z)$, so each keeps its sign on a product neighborhood $\mathcal O\times V$ of $(\theta,z)$. Hence, for all $(\theta',x)\in\mathcal O\times V$, we have $\alpha_{\theta}(x)=A$, which is equivalent to $x\in R_A(\theta')$. Consequently, 
\[
V\subseteq R_A(\theta).
\] Since the nonempty open set $V$ can not be contained in the finite union of lower-dimensional cells of $\mathcal C(\theta)$, it meets $\operatorname{int}C$ for some full-dimensional cell $C$. It follows that $A=A(C)$, and hence $A \in \mathcal{A}_\theta$.
\end{proof}

Write
\[
\mathcal{A}_{\theta}^{\mathrm{str}}
:=\{\alpha_\theta(z):z\text{ is }\theta\text{-generic}\}
\subseteq \mathcal{A}_{\theta}.
\]
These are the canonical patterns that occur at $\theta$-generic points, and hence the only patterns seen by the evaluation maps of \Cref{def:unbiased_fundim}. By \Cref{lem:evaluation_local}, the $\theta$-generic points with a given pattern in $\mathcal{A}_{\theta}^{\mathrm{str}}$ form a nonempty open subset of $\mathbb R^{n_0}$.

\subsection{Comparison with the pattern parametrization}
\label{subsec:functional_comparison}

In $d\varphi^{\mathcal{A}_{\theta}}(\theta)$, the set $\mathcal{A}_{\theta}$ is fixed before
differentiating and the varying parameter ranges over the entire weight
space in the polynomial map $\varphi^{\mathcal{A}_{\theta}}$.
We first record the algebraic rank fact used in the comparison.
For the empty pattern set, we interpret $\varphi^{\varnothing}$ as the
unique map to the zero-dimensional vector space and
$U^{\varnothing}$ as a point.

\begin{lemma}[Generic rank of a pattern parametrization]
\label{lem:pattern_map_generic_rank}
For any fixed finite pattern set $\mathcal A$, the generic rank of
$d\varphi^{\mathcal A}$ equals $\dim U^{\mathcal A}$. In particular,
\[
\operatorname{rank}d\varphi^{\mathcal A}(\theta')
\leq\dim U^{\mathcal A}
\quad\text{for every real parameter }\theta'.
\]
The equality locus is a nonempty real Zariski-open subset of the
parameter space, and its complement has Lebesgue measure zero.
\end{lemma}

\begin{proof}
Write $N=|\mathcal A|$. The map $\varphi^{\mathcal A}$
has $q=Nn_mn_0$ coordinates, namely the entries
of the matrices $M_A$, for all $A\in\mathcal A$. Let $F_1,\ldots,F_q$ denote these coordinate polynomials
in the weight coordinates $t_1,\ldots,t_D$. The coordinate ring of its complex image closure is $\mathbb C[F_1,\ldots,F_q]$. Thus
\[
\operatorname{rank}_{\mathbb C(t_1,\ldots,t_D)}
\left(\frac{\partial F_i}{\partial t_j}\right)_{i,j}
=\operatorname{trdeg}_{\mathbb C}\mathbb C(F_1,\ldots,F_q)
=\dim U^{\mathcal A}.
\]
Put $r=\dim U^{\mathcal A}$. The preceding equality
shows that the Jacobian matrix has rank $r$
over $\mathbb C(t_1,\ldots,t_D)$, so we get
$r\leq\min(q,D)$.

If $r<\min(q,D)$, every $(r+1)\times(r+1)$ minor
vanishes identically. Thus
$\operatorname{rank}d\varphi^{\mathcal A}(\theta')\leq r$.
If $r=\min(q,D)$, the same inequality follows directly
from the size of the Jacobian. In either case,
\[
\operatorname{rank}d\varphi^{\mathcal A}(\theta')
\leq r=\dim U^{\mathcal A}
\qquad\text{for every real parameter }\theta'.
\] If $r>0$, the locus where Jacobian drops rank is the common zero set of the $r\times r$ minors. At least one of these minors is a nonzero polynomial with real coefficients, so it cannot vanish at every real point. The locus where $\operatorname{rank} d\varphi^\mathcal{A}(\theta')=\dim U^\mathcal{A}$ is therefore a nonempty real Zariski-open set. Its complement is contained in the zero set of a nonzero real polynomial, hence has empty interior and Lebesgue measure zero. If $r=0$, equality holds everywhere. 
\end{proof}

\functionaldimension*

\begin{proof}
(i) Fix a finite ordered set $Z=(z_1,\ldots,z_k)$ of
$\theta$-generic points, and write $A_i=\alpha_\theta(z_i)$.
By \Cref{lem:evaluation_local}, there is a common neighborhood
$\mathcal O$ of $\theta$ on which every sample point $z_i$
retains its activation pattern $A_i$. Thus
\[
f_{\theta'}(z_i)=M_{A_i}(\theta')z_i
\qquad
\text{for all }\theta'\in\mathcal O
\text{ and }i=1,\ldots,k.
\]

Define
\[
\begin{aligned}
L_Z:
\prod_{A\in \mathcal{A}_{\theta}^{\mathrm{str}}}
\mathbb R^{n_m\times n_0}
&\longrightarrow \mathbb R^{n_m\times k},\\
(Y_A)_{A\in \mathcal{A}_{\theta}^{\mathrm{str}}}
&\longmapsto
\bigl[Y_{A_1}z_1\mid\cdots\mid Y_{A_k}z_k\bigr].
\end{aligned}
\]
Here the input is a tuple containing one matrix $Y_A$ for
each strict pattern $A$. For each sample point $z_i$, the map selects the matrix $Y_{A_i}$ and evaluates it at $z_i$. Since the points $z_i$ and their labels $A_i$ are fixed, $L_Z$ is linear.

Identifying the block matrix $\varphi^{\mathcal{A}_{\theta}^{\mathrm{str}}}
(\theta')$ with its tuple of matrix blocks, we obtain
\[
E_Z(\theta')
=L_Z\bigl(\varphi^{\mathcal{A}_{\theta}^{\mathrm{str}}}(\theta')\bigr) \qquad\text{for every }\theta'\in\mathcal O.
\]
Differentiating at $\theta$, we get
\[
dE_Z(\theta)
=L_Z\circ d\varphi^{\mathcal{A}_{\theta}^{\mathrm{str}}}(\theta).
\]
Thus,
\[
\operatorname{rank}JE_Z(\theta)
\leq
\operatorname{rank}d\varphi^{\mathcal{A}_{\theta}^{\mathrm{str}}}(\theta).
\]
Conversely, suppose $\mathcal{A}_{\theta}^{\mathrm{str}}$ is nonempty. For each
$A\in \mathcal{A}_{\theta}^{\mathrm{str}}$, choose linearly independent $\theta$-generic points $z_{A,1},\ldots,z_{A,n_0}$ realizing $A$. The matrix $X_A=[z_{A,1}\mid\cdots\mid z_{A,n_0}]$ is invertible.
For the combined sample set $Z$, the outputs belonging to pattern $A$ form the block $Y_AX_A$, from which we recover $Y_A=(Y_AX_A)X_A^{-1}$. Thus $L_Z$ is injective, and
\[
\operatorname{rank}(L_Z\circ d\varphi^{\mathcal{A}_{\theta}^{\mathrm{str}}}(\theta))
=\operatorname{rank}d\varphi^{\mathcal{A}_{\theta}^{\mathrm{str}}}(\theta).
\]
It follows that
\[
\operatorname{rank}JE_Z(\theta)
=\operatorname{rank}d\varphi^{\mathcal{A}_{\theta}^{\mathrm{str}}}(\theta).
\]
Since $\dimfunzero(\theta)$ is the maximum of
$\operatorname{rank}JE_Z(\theta)$ over all admissible
finite sample sets, we obtain
\[
\dimfunzero(\theta)
\geq\operatorname{rank}d\varphi^{\mathcal{A}_{\theta}^{\mathrm{str}}}(\theta).
\]
Together with the opposite inequality proved above,
this gives
\[
\dimfunzero(\theta)
=\operatorname{rank}d\varphi^{\mathcal{A}_{\theta}^{\mathrm{str}}}(\theta).
\]
If $\mathcal{A}_{\theta}^{\mathrm{str}}=\varnothing$, both sides
are zero by definition.

The map $\varphi^{\mathcal{A}_{\theta}^{\mathrm{str}}}$ is obtained from
$\varphi^{\mathcal{A}_{\theta}}$ by deleting blocks, so its differential cannot
have larger rank. The last inequality is
\Cref{lem:pattern_map_generic_rank}.

(ii) The hypothesis gives $\mathcal{A}_{\theta}^{\mathrm{str}}=\mathcal{A}_{\theta}$, so (i)
applies. In a shallow network, a nonzero row $w_k^\top$ vanishes only
on its hyperplane, which does not meet any chamber interior.
Thus all chamber interiors consist of $\theta$-generic points if
every row is nonzero. A zero row prevents any point from being
$\theta$-generic, proving the converse.

(iii) Nondegeneracy gives the wall-graph hypotheses by \Cref{lem:cube}.
Hence \Cref{thm:shallow_generators}~(ii) gives
\[
\varphi^{\mathcal{A}_{\theta}}=\slopeembed{\mathcal{A}_{\theta}}\circ\kappa,
\qquad
\kappa(W_1',W_2')=(c_k'w_k'^\top)_{k=1}^{n_1},
\]
where $\slopeembed{\mathcal{A}_{\theta}}$ is linear and injective. Hence
$\operatorname{rank}d\varphi^{\mathcal{A}_{\theta}}(\theta)
=\operatorname{rank}d\kappa(\theta)$.
The $k$-th component of $\kappa$ uses only the variables $(c_k,w_k)$, and its differential is
\[
(\dot c_k,\dot w_k)
\longmapsto\dot c_k w_k^\top+c_k\dot w_k^\top.
\]
Since $w_k \neq 0$ because of nondegeneracy, if $c_k\neq0$  its kernel is the one-dimensional space spanned by $(c_k,-w_k)$, which is the tangent direction of the rescaling $$(c_k,w_k)\mapsto(\lambda c_k,\lambda^{-1}w_k).$$ 
Indeed, let $(\dot c_k,\dot w_k)$ be in the kernel, so that $\dot c_k w_k^\top+c_k\dot w_k^\top=0$. Multiplying on the right by $w_k$ gives $\dot c_k=t c_k$ with $t=-(\dot w_k\cdot w_k)/|w_k|^2$. Substituting this back gives $c_k(t w_k+\dot w_k)^\top=0$, and hence $\dot w_k=-t w_k$ since $c_k\neq0$. Thus $(\dot c_k,\dot w_k)=t(c_k,-w_k)$. Its rank is therefore $n_2+n_0-1$.

If $c_k=0$, the map is $\dot c_k\mapsto\dot c_k w_k^\top$, which has rank $n_2$. 

By part (ii), we conclude that
\begin{align*}
\dimfunzero(\theta)&=\operatorname{rank}d\varphi^{\mathcal{A}_{\theta}}(\theta) =\operatorname{rank}d\kappa(\theta)\\
&= (n_0+n_2-1)(n_1- \big|\{k:c_k=0\}\big|)+n_2 \cdot\big|\{k:c_k=0\}\big|\\
&=n_1(n_0+n_2-1)-(n_0-1)\cdot \big|\{k:c_k=0\}\big|,
\end{align*}
and $\dim U^{\mathcal{A}_{\theta}}=n_1(n_0+n_2-1)$ follows from \Cref{thm:shallow_generators}~(ii).

(iv) For each $\mathcal A\subseteq\{0,1\}^{n_1}$, let
\[
\Omega_{\mathcal A}
=
\left\{\theta':
\operatorname{rank}d\varphi^{\mathcal A}(\theta')
=\dim U^{\mathcal A}\right\},
\]
with $\Omega_{\varnothing}$ equal to the whole parameter space. By \Cref{lem:pattern_map_generic_rank}, each
$\Omega_{\mathcal A}$ is open, dense, and of full Lebesgue measure. The same holds for
\[
\mathcal H:=\{\theta':W_1'\text{ has no zero row}\},
\]
since its complement is the zero set of the nonzero polynomial
$\prod_{k=1}^{n_1}\sum_{j=1}^{n_0}(W_1')_{kj}^2$.

There are finitely many possible pattern sets, so
\[
\mathcal E
:=\mathcal H\cap
\bigcap_{\mathcal A\subseteq\{0,1\}^{n_1}}\Omega_{\mathcal A}
\]
is open and dense, and its complement has measure zero.
For every $\theta\in\mathcal E$, part~(ii) and
$\theta\in\Omega_{\mathcal{A}_{\theta}}$ give, respectively,
\[
\dimfunzero(\theta)
=\operatorname{rank}d\varphi^{\mathcal{A}_{\theta}}(\theta)
=\dim U^{\mathcal{A}_{\theta}}.
\]
\end{proof}

\Cref{def:unbiased_fundim} samples only at $\theta$-generic points, whereas \cite{GLMWFunctionalDimension} samples at all parametrically smooth points. Following \cite[Definition~3.9]{GLMWFunctionalDimension}, we call $z\in\mathbb R^{n_0}$ \emph{parametrically smooth} for $\theta$ if $(\theta',x)\mapsto f_{\theta'}(x)$ is smooth on a neighborhood of $(\theta,z)$, and we let
\[
\dim^{\mathrm{sm}}_{\mathrm{fun}}(\theta)
=\max_Z\operatorname{rank}JE_Z(\theta),
\]
where $Z$ ranges over finite sets of parametrically smooth points for $\theta$, with $Z=\varnothing$ allowed as in \Cref{def:unbiased_fundim}. The next proposition shows that, at any depth, the two quantities agree under the hypothesis of \Cref{thm:functional_dimension}~(ii). In particular, the formulas in \Cref{thm:functional_dimension}~(iii) and~(iv) do not depend on the sampling convention.

\begin{proposition}[Sampling at all parametrically smooth points] \label{prop:functional_smooth_sampling}
Let $\theta\in\mathcal P_{n_0,\ldots,n_m}$. Then
$\dimfunzero(\theta)\leq\dim^{\mathrm{sm}}_{\mathrm{fun}}(\theta)$, with equality if the interior of every maximal-dimensional cell of $\mathcal C(\theta)$ contains a $\theta$-generic point. For shallow
networks, this condition holds if and only if $W_1$ has no zero row.
\end{proposition}

\begin{proof}
Every $\theta$-generic point is parametrically smooth by
\Cref{lem:evaluation_local}, so
\[
\dimfunzero(\theta)\leq\dim^{\mathrm{sm}}_{\mathrm{fun}}(\theta).
\]

If the interior of every maximal-dimensional cell of $\mathcal{C}(\theta)$ contains a $\theta$-generic point, then we claim that the $\theta$-generic points are dense in $\mathbb R^{n_0}$. To see this, let $C \in \mathcal{C}(\theta)$ with $\dim C=n_0$. We know that for every $1\leq \ell \leq m, 1\leq i \leq n_\ell$, each hidden preactivation $h^{(\ell)}_{\theta,i}$ is strictly positive, strictly negative, or identically zero on $\operatorname{int}(C)$. Since $\operatorname{int}(C)$ contains a $\theta$-generic point, we know that each $h^{(\ell)}_{\theta,i}$ is either strictly positive or negative on $\operatorname{int}(C)$, and hence all points in $\operatorname{int}(C)$ are $\theta$-generic. Since the interiors of the maximal-dimensional cells are dense in $\mathbb R^{n_0}$, we can conclude that $\theta$-generic points are dense in $\mathbb R^{n_0}$ as desired. 

Now let $Z=(z_1,\ldots,z_k)$ be parametrically smooth points for $\theta$. Since $(\theta',x)\mapsto f_{\theta'}(x)$ is smooth near each $(\theta,z_i)$, the matrix $JE_{Z'}(\theta)$ depends continuously on
$Z'=(z_1',\ldots,z_k')$ near $Z$. A nonzero minor remains nonzero under small perturbations, so
$\operatorname{rank}JE_{Z'}(\theta)\geq\operatorname{rank}JE_Z(\theta)$
for all $Z'$ sufficiently close to $Z$. By density, we may choose such a $Z'$ consisting of $\theta$-generic points. Hence
\[
\operatorname{rank}JE_Z(\theta)
\leq\operatorname{rank}JE_{Z'}(\theta)
\leq\dimfunzero(\theta).
\]
Taking the maximum over $Z$ gives
$\dim^{\mathrm{sm}}_{\mathrm{fun}}(\theta)\leq\dimfunzero(\theta)$. As a result, we conclude that
\[
\dimfunzero(\theta) =\dim^\mathrm{sm}_{\mathrm{fun}}(\theta)
\]
if the interior of every maximal-dimensional cell of $\mathcal{C}(\theta)$ contains a $\theta$-generic point.

For shallow networks, the statement follows by \Cref{thm:functional_dimension}~(ii).
\end{proof}

Both inequalities in \Cref{thm:functional_dimension}~(i) can be strict. The first is strict when deleting the blocks indexed by $\mathcal A_\theta\setminus\mathcal A_\theta^{\mathrm{str}}$ lowers the rank of $d\varphi^{\mathcal A_\theta}(\theta)$, as in
\Cref{ex:functional_zero_row}. The second is strict when $d\varphi^{\mathcal A_\theta}$ has less than its generic rank at $\theta$, as in \Cref{ex:functional_zero_column}.

\begin{example}[No $\theta$-generic points]\label{ex:functional_zero_row}
For architecture $(2,2,1)$, take
\[
W_1=\begin{pmatrix}1&0\\0&0\end{pmatrix},
\qquad W_2=(1,1).
\]
The second hidden preactivation vanishes identically, so there are no $\theta$-generic points, and $\dimfunzero(\theta)=0$ by the convention for $Z=\varnothing$. On the other hand,
$\mathcal{A}_{\theta}=\{(1,0),(0,0)\}$ and, in this ordering,
\[
\varphi^{\mathcal{A}_{\theta}}(\theta')
=(c_1'w_1'^\top,0),
\qquad
\operatorname{rank}d\varphi^{\mathcal{A}_{\theta}}(\theta)
=\dim U^{\mathcal{A}_{\theta}}=2.
\]
Thus the first inequality in \Cref{thm:functional_dimension}~(i) is
strict, and the second is an equality.
\end{example}

\begin{example}[A zero output column]\label{ex:functional_zero_column}
For architecture $(2,2,1)$, take $W_1=I_2$ and $W_2=(1,0)$. This parameter is nondegenerate and $W_2$ has one zero column, so
\Cref{thm:functional_dimension}~(ii) and~(iii) give
\[
\dimfunzero(\theta)
=\operatorname{rank}d\varphi^{\mathcal{A}_{\theta}}(\theta)=3
<4=\dim U^{\mathcal{A}_{\theta}}.
\]
Thus the first inequality in \Cref{thm:functional_dimension}(i) is an
equality here, and only the second is strict.
\end{example}

\appendix

\section{Clocks (Appendix by Patrick Frem, Eduardo Madero, and Weiya Wang)}\label{app:relu_clocks}

This appendix gives the clock construction and proves
\Cref{thm:clock_pattern_equivalence,thm:clock_counting_formula}.
Throughout, assume $n_1\geq1$. The indices of rays, regions, vertices, and patterns are understood modulo $2n_1$.

We begin with the geometric description of the activation patterns, translate the oriented lines into clocks, and then identify the clock symmetries that correspond to a simultaneous permutation of the neurons. The final step counts the resulting dihedral orbits. Only the first-layer arrangement of a shallow unbiased network is used; the second-layer matrix does not affect the construction.

\subsection{Clock construction and labeled patterns}
\label{app:clock_construction}

Let $\theta\in\mathcal{P}^{\mathrm{nd}}_{2,n_1,n_2}$ be an unbiased
parameter, and write its first-layer weight matrix as
\[
W_1=
\begin{pmatrix}
w_1^\top\\
\vdots\\
w_{n_1}^\top
\end{pmatrix},
\qquad w_i\in\mathbb{R}^2.
\]

Each nonzero neuron $w_i$ determines an oriented line
through the origin,
\[
H_i=\{x\in\mathbb{R}^2:w_i^\top x=0\},
\]
where the orientation is specified by the active
half-plane
\[
H_i^+=\{x\in\mathbb{R}^2:w_i^\top x>0\}.
\]

When $\theta \in \mathcal{P}^{\mathrm{nd}}_{2,n_1,n_2}$, the $n_1$ distinct lines divide $\mathbb{R}^2$ into exactly $2n_1$ maximal-dimensional cones. Each cone is associated with a unique activation pattern (\S\ref{sec:prelim}; see also \Cref{lem:app_adjacency}), and hence
\[
|\mathcal{A}_{\theta}|=2n_1.
\]

\paragraph{A running example.}
Consider a shallow ReLU neural network with architecture $(2,3,1)$. Let $W_1\colon\mathbb R^2\to\mathbb R^3$ be the linear map defined by
\[
W_1=\begin{pmatrix}1&0\\0&1\\1&-1\end{pmatrix}.
\]
For an input $(x,y)\in\mathbb R^2$, the linear transformation yields the three preactivations $x$, $y$, and $x-y$. Applying the ReLU activation coordinate-wise gives
\[
\varsigma\left(W_1\begin{pmatrix}x\\y\end{pmatrix}\right)
=\begin{pmatrix}\max(0,x)\\\max(0,y)\\\max(0,x-y)\end{pmatrix}.
\]
The three hidden neurons therefore induce the lines
\[
H_1=\{x=0\},\qquad H_2=\{y=0\},\qquad H_3=\{x-y=0\},
\]
which partition the input plane into six regions. Each region is uniquely characterized by the signs of these preactivations. Because the signs are constant on each region, the network is linear there; the closures of these regions, together with their faces, form the canonical polyhedral complex $\mathcal C(\theta)$.

\begin{figure}[htbp]
\centering
\includegraphics[width=.48\linewidth]{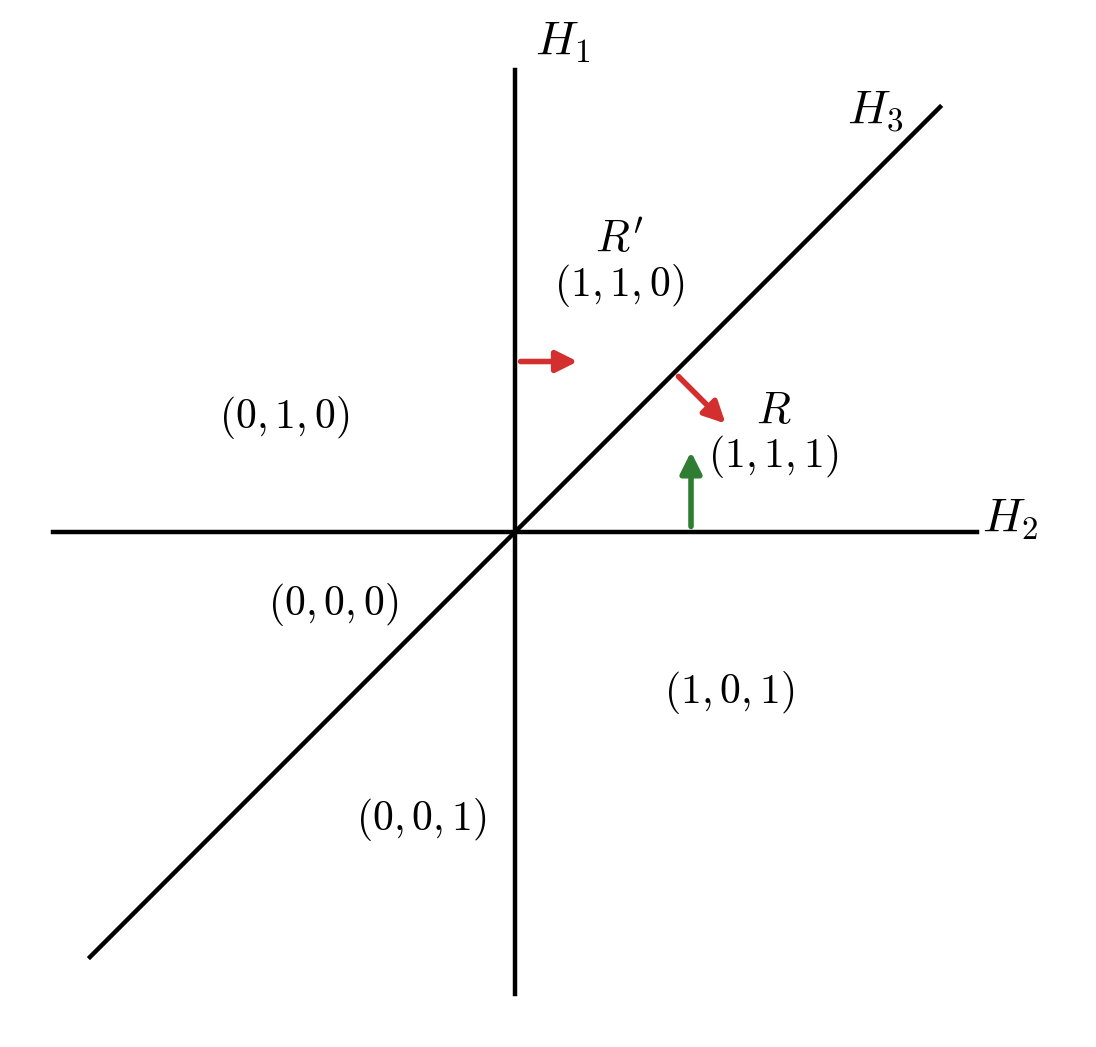}
\caption{The first-layer arrangement of the running example, with its six activation patterns. The labels $R$ and $R'$ refer to the two regions described in the text.}
\label{fig:app_clock_running}\label{fig:relu_fan_example}
\end{figure}

For example, consider $R=\{x>y>0\}$. Within $R$ we have $x>0$, $y>0$, and $x-y>0$, so all three neurons are active and the corresponding activation pattern is $(1,1,1)$. Similarly, in $R'=\{y>x>0\}$ we have $x>0$, $y>0$, and $x-y<0$, giving the activation pattern $(1,1,0)$. In \Cref{fig:app_clock_running}, the arrows indicate the positive side of each hyperplane, that is, the side on which the corresponding hidden neuron is active. Crossing a ray changes the corresponding entry of the activation pattern from $0$ to $1$, or from $1$ to $0$. Opposite regions have complementary patterns because $w_k^\top(-x)=-w_k^\top x$. These observations will be made precise
in \Cref{lem:app_crossing,lem:app_adjacency}.

The combinatorial information of this oriented
hyperplane arrangement can be encoded by a regular
$2n_1$-gon, called a \emph{clock}.

\begin{definition}[Clock]
\label{def:clock}
A \emph{clock} $\clocksign$ of size $n_1$ is a sign function
\[
\clocksign:\mathbb{Z}/2n_1\mathbb{Z}
\longrightarrow\{+,-\}
\]
satisfying
\[
\clocksign(i+n_1)=-\clocksign(i)
\]
for every $i\in\mathbb{Z}/2n_1\mathbb{Z}$.

We visualize its cyclically indexed domain as $2n_1$ vertices, labeled
\[
0,1,\ldots,2n_1-1,
\]
and, when $n_1\geq2$, arranged as a regular $2n_1$-gon.
We denote by $\clockset{n_1}$ the collection of all clocks of size $n_1$.
\end{definition}

A clock is completely determined by the signs
assigned to its first $n_1$ vertices. Consequently,
\[
|\clockset{n_1}|=2^{n_1}.
\]

Given $\theta\in\mathcal P^{\mathrm{nd}}_{2,n_1,n_2}$, we construct an associated clock by choosing one ray of $H_1$ as vertex $0$ and labeling the remaining $2n_1-1$ rays
counterclockwise. Thus the rays $r_0,\ldots,r_{2n_1-1}$ occur in counterclockwise order and $r_{i+n_1}=-r_i$. Let $R_i$ be the open region between $r_i$ and $r_{i+1}$. For each $i\in\mathbb Z/2n_1\mathbb Z$,
let
\[
h_i\in\{1,\ldots,n_1\}
\]
be the unique neuron index such that
\[
r_i\subseteq H_{h_i}.
\]
Then $h_{i+n_1}=h_i$, and $i\mapsto h_i$ restricts to a bijection from $\{0,\ldots,n_1-1\}$ onto
$\{1,\ldots,n_1\}$. In particular, $h_0=1$.

Each vertex is assigned a sign according to the
orientation of its corresponding hyperplane:
the sign $\clocksign_\theta(i)$ is $+$ if neuron $h_i$ is active on the counterclockwise side of the ray $r_i$, i.e.\ on $R_i$, and $-$ otherwise. Since $R_{i+n_1}=-R_i$, we have $\clocksign_\theta(i+n_1)= -\clocksign_\theta(i)$.
Equivalently, $\clocksign_\theta(i)=+$ if and only if, crossing $r_i$ from $R_{i-1}$ to $R_i$, coordinate $h_i$ changes from $0$ to $1$.

Thus, the vertices of the clock correspond to
the rays of the activation hyperplanes, while
the maximal-dimensional cones correspond to
the regions between consecutive rays.

\begin{figure}[htbp]
\centering
\includegraphics[width=.92\linewidth]{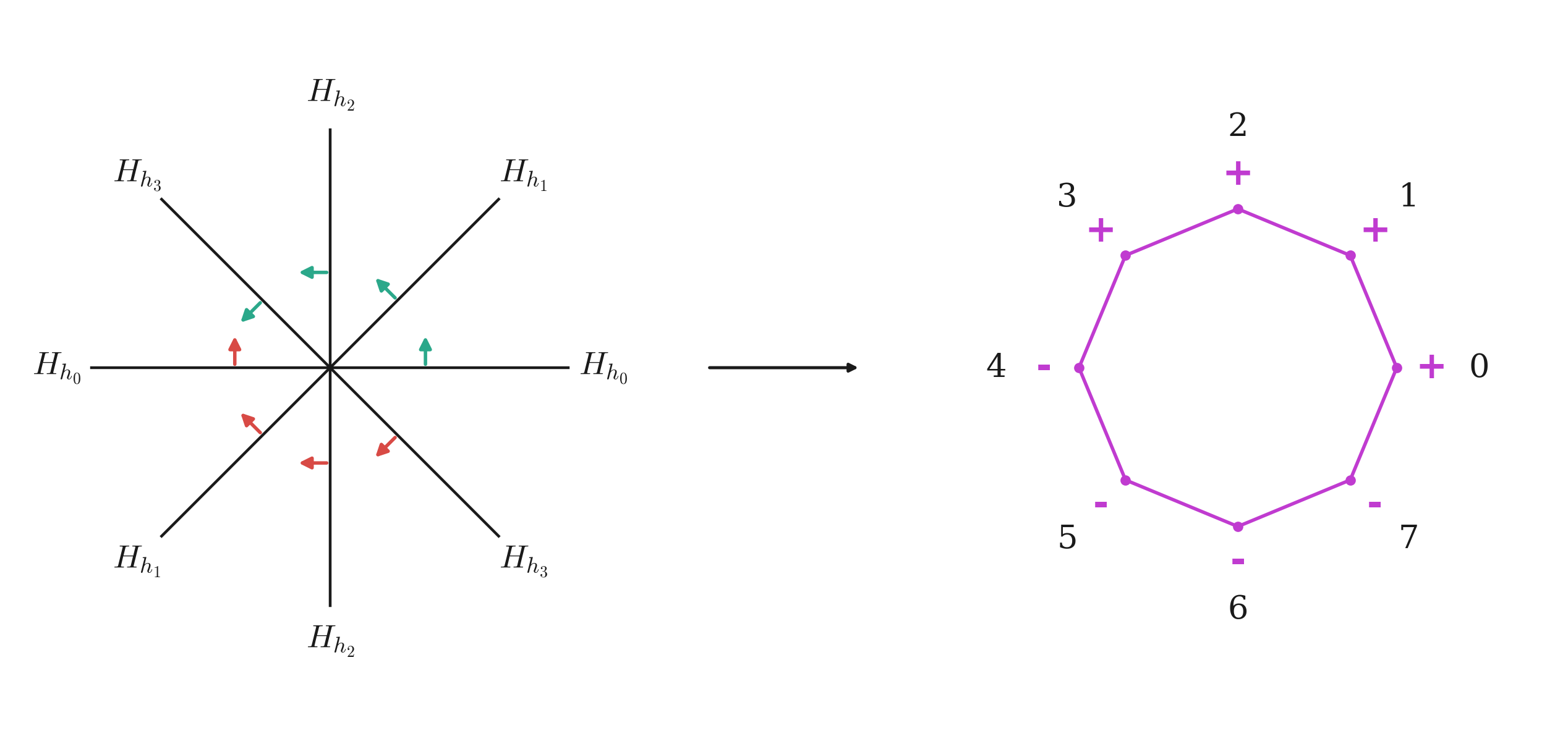}
\caption{Translation of an oriented arrangement with $n_1=4$ into a clock. The line containing ray $r_i$ is $H_{h_i}$. The angles are replaced by equal spacing, while cyclic order and crossing signs are retained.}
\label{fig:app_clock_construction}\label{fig:clock}
\end{figure}

The activation patterns can be read off from a clock together with the neuron labels attached to its antipodal vertex pairs. Formally, a \emph{ray labeling} of size $n_1$ is a sequence
\[
h=(h_i)_{i\in\mathbb Z/2n_1\mathbb Z}
\]
such that $h_{i+n_1}=h_i$ and
$(h_0,\ldots,h_{n_1-1})$ is a permutation of
$(1,\ldots,n_1)$. For each $k\in\{1,\ldots,n_1\}$, let
$a_k\in\mathbb Z/2n_1\mathbb Z$ be the unique index satisfying
\[
h_{a_k}=k
\qquad\text{and}\qquad
\clocksign(a_k)=+.
\]
Define the \emph{labeled pattern set}
\[
P(\clocksign,h)
:=
\{A_\ell:\ell\in\mathbb Z/2n_1\mathbb Z\},
\]
where
\[
(A_\ell)_k
:=
\begin{cases}
1,
& \ell-a_k\in\{0,1,\ldots,n_1-1\}
  \pmod{2n_1},\\
0,
& \text{otherwise}.
\end{cases}
\]
Neuron $k$ is active exactly on the $n_1$ regions
$R_{a_k},\ldots,R_{a_k+n_1-1}$ counterclockwise from its positive ray. Thus, $A_\ell$ is the activation pattern of $R_\ell$, and
\[
\mathcal{A}_{\theta}=P(\clocksign_\theta,h)
\]
by \Cref{lem:app_adjacency}.

Notice, however, that the clock associated with
$\theta$ is not uniquely determined. Indeed,
the construction requires choosing one of the
two rays of $H_1$ as vertex $0$. The two possible
choices yield clocks that differ by a rotation
of $n_1$ vertices. Thus, if one choice gives the sign
function $\clocksign_\theta$, then the other gives
$-\clocksign_\theta$. The ray labeling is unchanged because
$h_{i+n_1}=h_i$.

We say that two clocks $\clocksign$ and $\clocksign'$ are
\emph{shadow-equivalent}, written
$\clocksign\sim_{\mathrm{shadow}}\clocksign'$, if
$\clocksign'=\clocksign$ or $\clocksign'=-\clocksign$.
Equivalently, $\clocksign'=\clocksign$ or
$\clocksign'=r^{n_1}\cdot\clocksign$. Each equivalence class
$[\clocksign]_{\mathrm{shadow}}=\{\clocksign,-\clocksign\}$
is called a \emph{shadow pair}.

\begin{figure}[htbp]
\centering
\includegraphics[width=.68\linewidth]{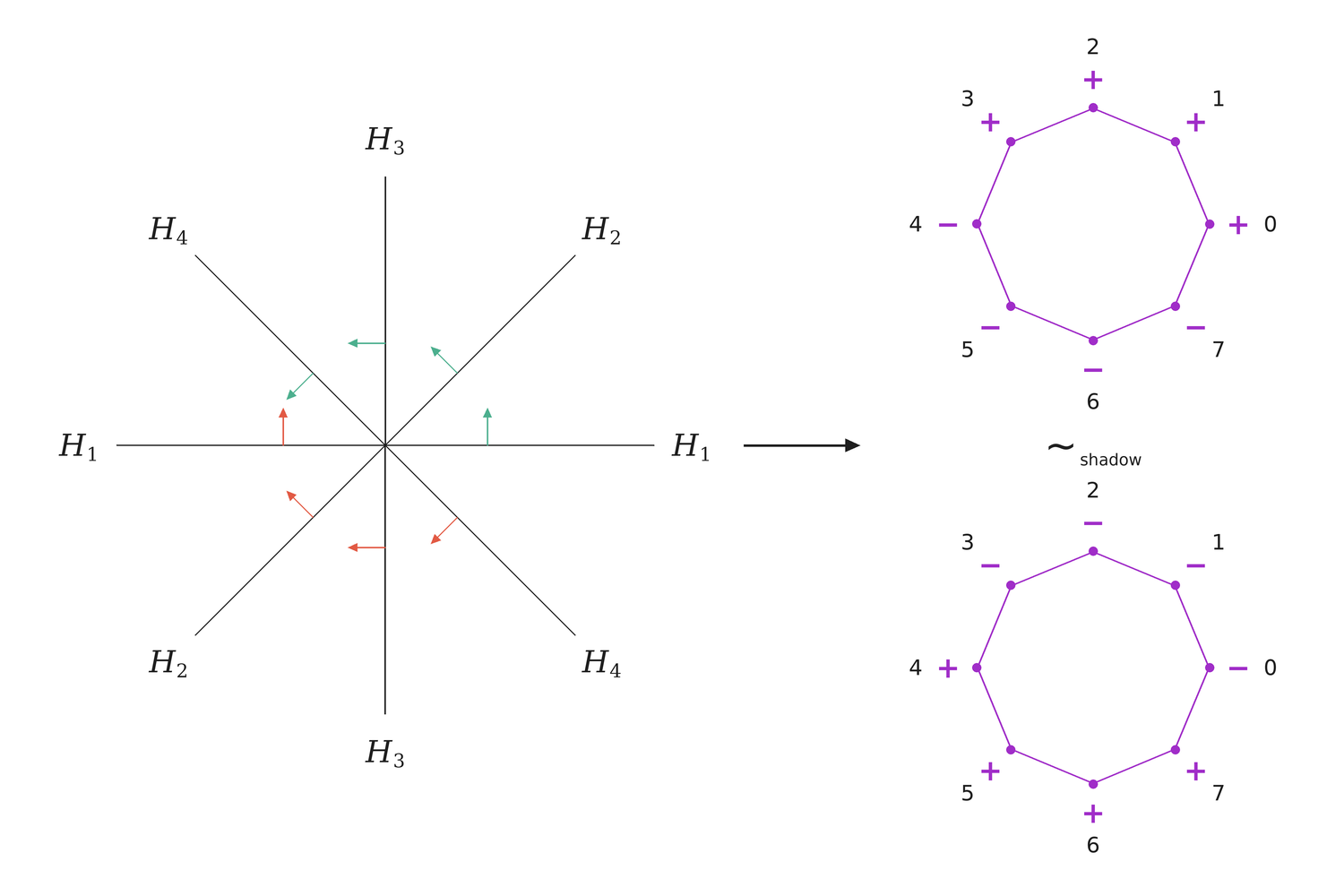}
\caption{The two choices of the initial ray of $H_1$ produce a shadow pair. Their sign functions differ by a half-turn: $r^{n_1}\cdot\clocksign=-\clocksign$.}
\label{fig:app_clock_shadow}
\end{figure}

Since the choice of the initial ray is the only choice in the
construction, the assignment
\[
\alpha:
\mathcal{P}_{2,n_1,n_2}^{\mathrm{nd}}
\longrightarrow
\clockset{n_1}/\sim_{\mathrm{shadow}},
\qquad
\theta\longmapsto[\clocksign_\theta]_{\mathrm{shadow}},
\]
is well defined, where $\mathcal{P}_{2,n_1,n_2}^{\mathrm{nd}}$
denotes the set of unbiased parameters whose
first-layer neurons are nonzero and pairwise
linearly independent. The map $\alpha$ is surjective (\Cref{lem:app_dihedral}(c)).

Moreover, a shadow-equivalence class of clocks
determines a unique pattern-equivalence class
of activation-pattern sets: relabeling the neurons by
$\pi\in S_{n_1}$ replaces each $h_i$ by $\pi(h_i)$ and applies
the same coordinate permutation to every pattern, while
$P(-\clocksign,h)=P(\clocksign,h)$
(\Cref{lem:app_dihedral}(a)).

\begin{figure}[htbp]
\centering
\includegraphics[width=.92\linewidth]{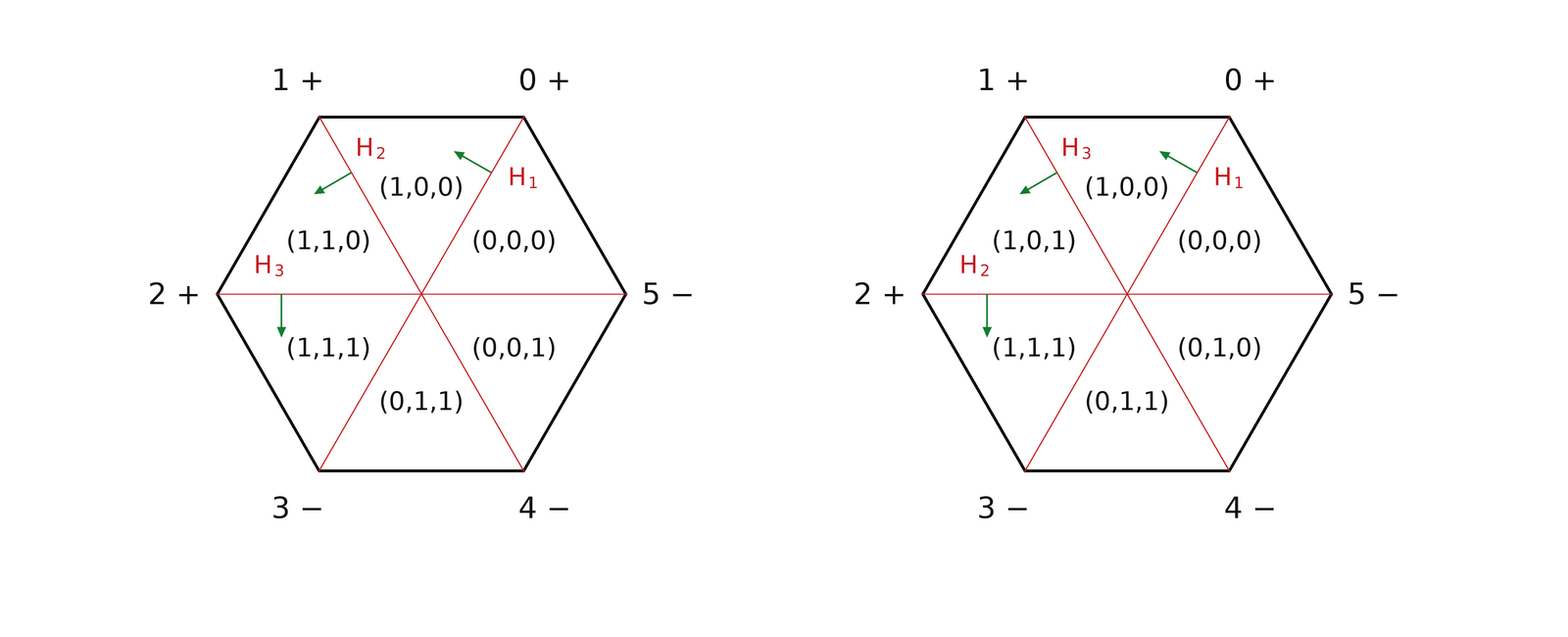}
\caption{The same clock with two ray labelings. On the left, $(h_0,h_1,h_2)=(1,2,3)$; on the right, $(h'_0,h'_1,h'_2)=(1,3,2)$. The pattern sets differ by the single simultaneous coordinate permutation $(2\ 3)$.}
\label{fig:app_clock_labeling}
\end{figure}

\paragraph{Why the labeling matters.}
A clock records which side of each line is active, but not the neuron
names attached to those lines. Each active half-plane consists of
$n_1$ consecutive regions together with their separating rays.
Consequently, the clock determines the activation-pattern set up to
the ray labeling $h$. Changing that labeling permutes the same
coordinates in every pattern.
For example, the two labeled clocks in
\Cref{fig:app_clock_labeling} have identical signs, but interchanging
the labels $2$ and $3$ changes their labeled pattern sets.
Thus a clock without a specified labeling determines a
pattern-equivalence class, rather than a particular labeled set.

To retain the parameter-to-clock description, let
\[
\mathcal S_{\mathrm{nd}}
=
\{\mathcal{A}_{\theta}:
\theta\in\mathcal P^{\mathrm{nd}}_{2,n_1,n_2}\}.
\]
For $\mathcal A,\mathcal A'\in\mathcal S_{\mathrm{nd}}$, write
\[
\mathcal A\sim\mathcal A'
\]
if there exists $\pi\in S_{n_1}$ such that
\[
\pi\cdot\mathcal A=\mathcal A',
\]
where the action on sets of activation patterns is defined in \Cref{def:setofactivationpattern_permutation}. This is the orbit equivalence relation induced by the action of $S_{n_1}$ on $\mathcal S_{\mathrm{nd}}$. In particular, $\theta\sim\theta'$ in the sense of
\Cref{def:patternequivalentparameter} if and only if
$\mathcal{A}_{\theta}\sim \mathcal{A}_{\theta'}$.

\begin{figure}[htbp]
\centering
\includegraphics[width=.65\linewidth]{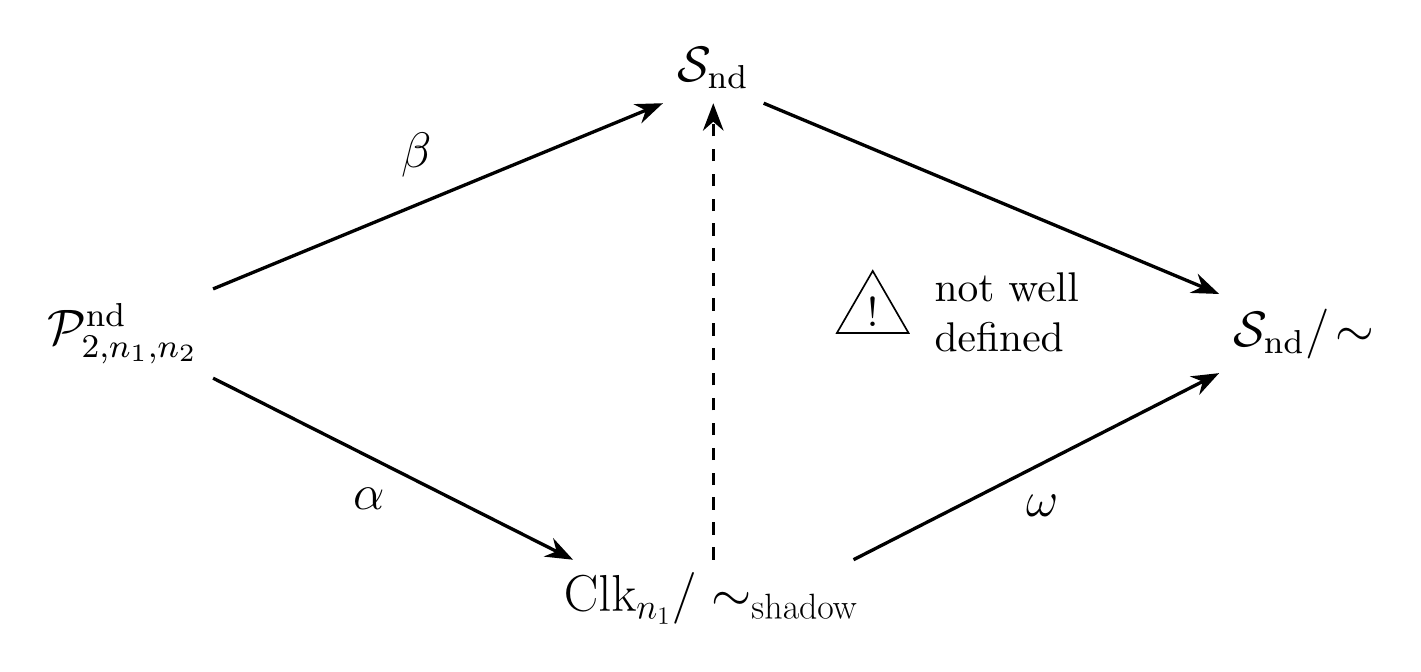}
\caption{The parameter-to-clock and parameter-to-pattern maps. The dotted arrow indicates that a labeled pattern set cannot be recovered without choosing neuron labels. The quotient uses one simultaneous permutation for the entire set.}
\label{fig:app_clock_maps}
\end{figure}

The maps represented in \Cref{fig:app_clock_maps} are
\[
\begin{aligned}
\beta &:
\mathcal P^{\mathrm{nd}}_{2,n_1,n_2}
\longrightarrow
\mathcal S_{\mathrm{nd}},
&&
\theta\longmapsto \mathcal{A}_{\theta},
\\
\omega &:
\clockset{n_1}/\sim_{\mathrm{shadow}}
\longrightarrow
\mathcal S_{\mathrm{nd}}/\sim,
&&
[\clocksign]_{\mathrm{shadow}}
\longmapsto
[P(\clocksign,h)]_{\sim},
\end{aligned}
\]
where $h$ is any ray labeling of size $n_1$ with $h_0=1$. The second map is independent of $h$ and of the representative of the shadow pair, by \Cref{lem:app_dihedral}(a). If $q$ is the quotient map on $\mathcal S_{\mathrm{nd}}$, then
\[
q\circ\beta=\omega\circ\alpha.
\]

\subsection{The dihedral action}\label{app:clock_action}

The natural symmetries of clocks are rotations
and reflections of the $2n_1$-gon, described
by the dihedral group $D_{2n_1}$ of order $4n_1$, generated by the rotation $r(i)=i+1$ and the reflections $f_c(i)=c-i$ of $\mathbb Z/2n_1\mathbb Z$ for $n_1 \geq 2$. It acts on sign functions by $(g\cdot\clocksign)(i):=\clocksign(g^{-1}(i))$. This action preserves the antipodal condition because every $g\in D_{2n_1}$ preserves antipodal pairs. Indeed, for a rotation $r^k$,
\[
r^k(i+n_1)=r^k(i)+n_1,
\]
while for a reflection $f_c(i)=c-i$,
\[
f_c(i+n_1)=c-i-n_1
\equiv c-i+n_1=f_c(i)+n_1
\pmod{2n_1}.
\] 
When $n_1=1$, the canonical activation-pattern set is
$\{0,1\}$ for every nondegenerate parameter, while the two size-one clocks belong to the same dihedral orbit. Thus \Cref{thm:clock_pattern_equivalence} holds in this case, and \Cref{thm:clock_counting_formula} gives $N(1)=1$.
Hence, throughout the remainder of this appendix, we assume $n_1\geq2$.

The shadow operation is precisely the rotation
by $n_1$ vertices, $\clocksign\mapsto r^{n_1}\cdot\clocksign=-\clocksign$, and $r^{n_1}$ is central in $D_{2n_1}$. Therefore, the dihedral
action descends to the shadow-equivalence
classes, and we have
\[
\left(
\clockset{n_1}/\sim_{\mathrm{shadow}}
\right)/D_{2n_1}
\cong
\clockset{n_1}/D_{2n_1}.
\]

\begin{figure}[htbp]
\centering
\includegraphics[width=.88\linewidth]{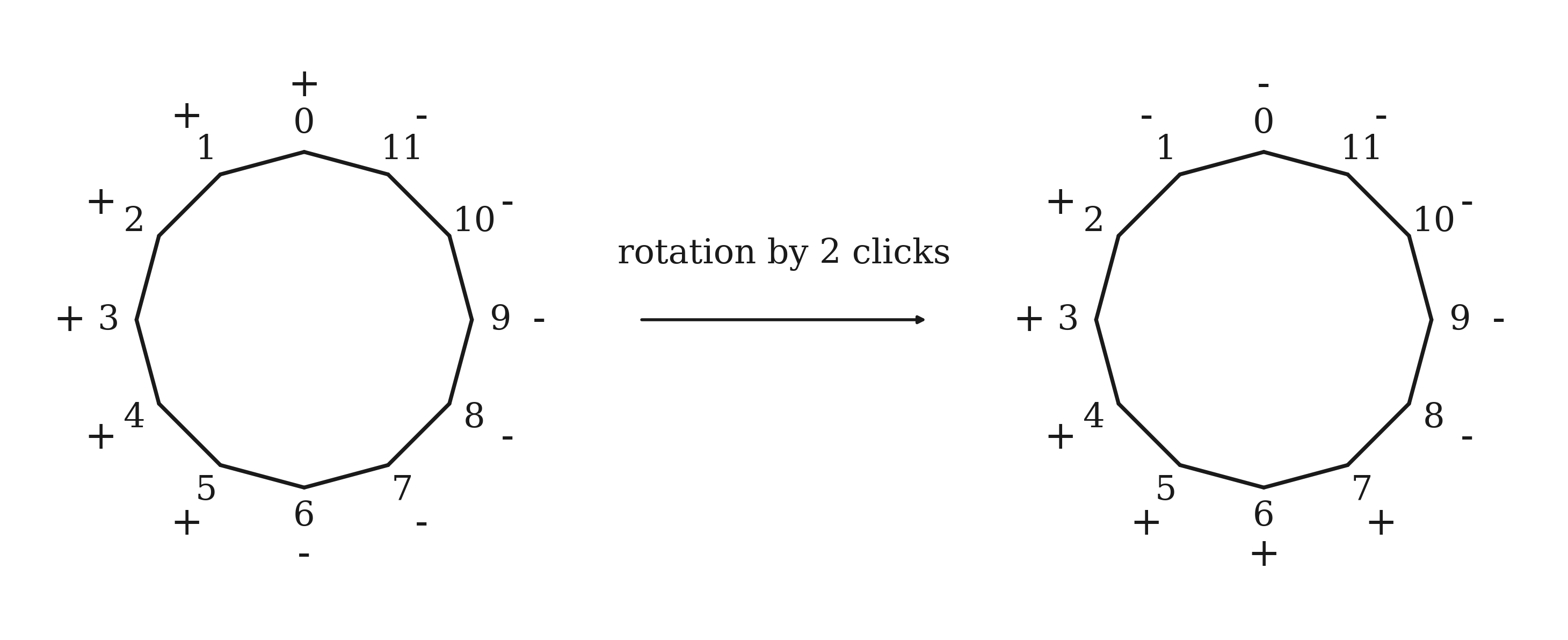}
\caption{A counterclockwise rotation by two clicks for $n_1=6$: the right-hand clock is $r^2\cdot\clocksign$.}
\label{fig:app_clock_rotation}
\end{figure}

\begin{figure}[htbp]
\centering
\includegraphics[width=.88\linewidth]{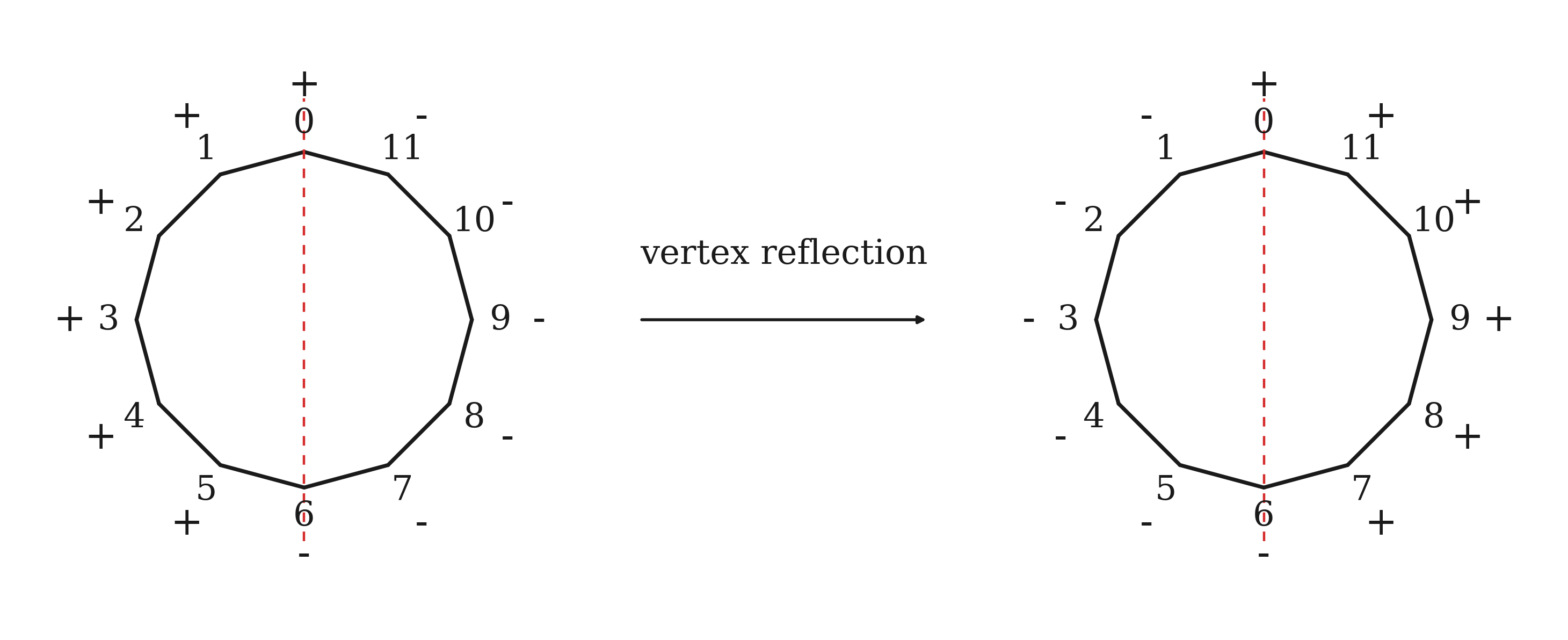}
\caption{Vertex reflection for $n_1=6$ across the axis through vertices $0$ and $6$, namely $f_0(i)=-i$.}
\label{fig:app_clock_vertex}
\end{figure}

\begin{figure}[htbp]
\centering
\includegraphics[width=.88\linewidth]{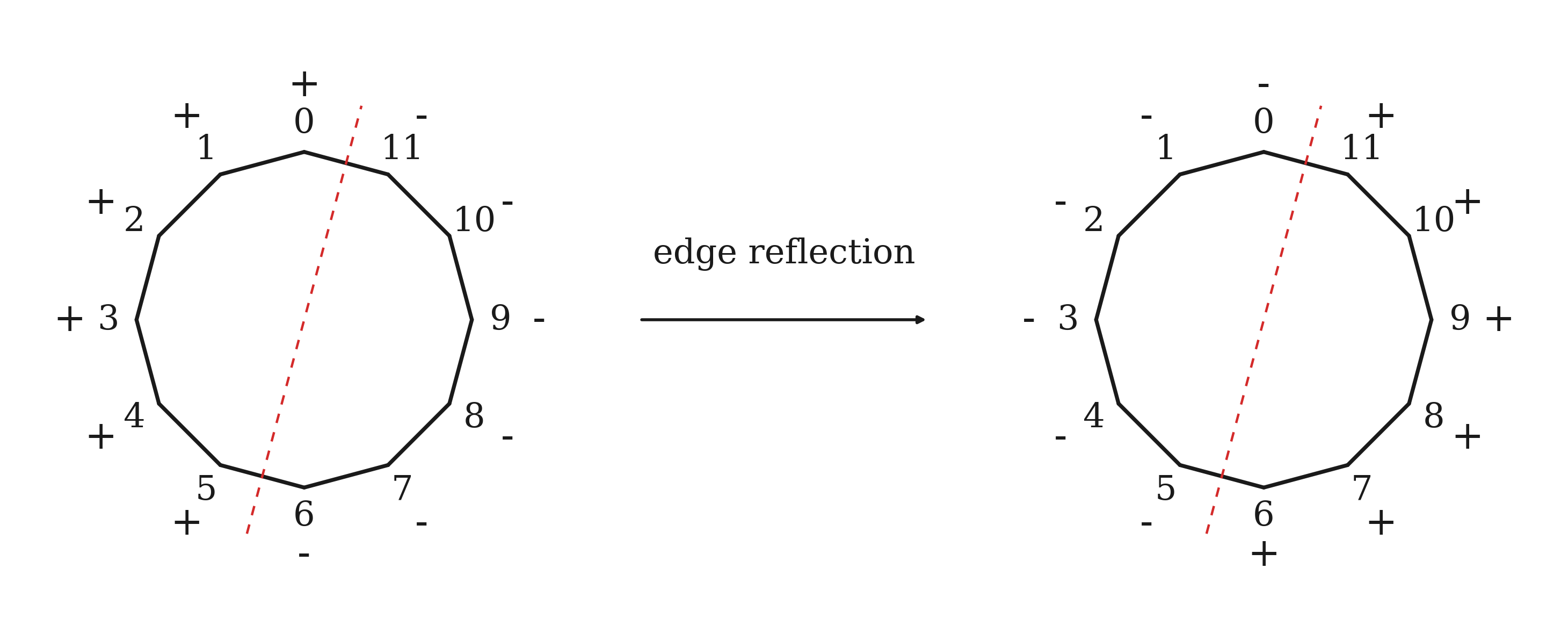}
\caption{Edge reflection for $n_1=6$ across the axis through the midpoints of edges $\{11,0\}$ and $\{5,6\}$, namely $f_{11}(i)=11-i$.}
\label{fig:app_clock_edge}
\end{figure}

\paragraph{Rotations and the two types of reflection.}
The first $n_1$ signs
$(\clocksign(0),\ldots,\clocksign(n_1-1))$ determine the whole clock.
A rotation by $k$ clicks transports the sign at $i$ to $i+k$, so
\[
(r^k\cdot\clocksign)(i)=\clocksign(i-k).
\]
A vertex reflection has the form $f_{2k}(i)=2k-i$ and its axis
passes through vertices $k$ and $k+n_1$. An edge reflection has
the form $f_{2k+1}(i)=2k+1-i$ and its axis passes through the
midpoints of opposite edges. Both act by
\[
(f_c\cdot\clocksign)(i)=\clocksign(c-i).
\]
These operations are illustrated in
\Cref{fig:app_clock_rotation,fig:app_clock_vertex,fig:app_clock_edge}.

\newpage 

\subsection{Crossings, adjacency, and symmetries}
\label{app:clock_lemmas}

For $i,j\in\mathbb Z/2n_1\mathbb Z$, let $d(i,j)$ denote their
cyclic distance. Given a clock $\clocksign$ and a ray labeling
$h=(h_i)$, write
\[
P(\clocksign,h)
=
\{A_\ell:\ell\in\mathbb Z/2n_1\mathbb Z\}
\]
for the associated labeled pattern set defined in
\Cref{app:clock_construction}.

\begin{lemma}[Crossing a ray]\label{lem:app_crossing}
For every $i\in\mathbb Z/2n_1\mathbb Z$, the patterns
$A_{i-1}$ and $A_i$ differ exactly in coordinate $h_i$. This coordinate changes from $0$ to $1$ if $\clocksign(i)=+$ and from $1$ to $0$ if
$\clocksign(i)=-$. Moreover,
\[
A_{i+n_1}=\mathbf 1-A_i.
\]
\end{lemma}

\begin{proof}
Fix $k\in\{1,\ldots,n_1\}$. By definition,
\[
\{\ell\in\mathbb Z/2n_1\mathbb Z:(A_\ell)_k=1\}
=
\{a_k,\ldots,a_k+n_1-1\},
\]
which is a cyclic interval of length $n_1$. Therefore, as one passes
from $A_{i-1}$ to $A_i$, the $k$-th coordinate changes if and only
if $i\in\{a_k,a_k+n_1\}$, or equivalently, if and only if
$h_i=k$. When $h_i=k$, the two possibilities are
\[
i=a_k
\quad\Longleftrightarrow\quad
\clocksign(i)=+,
\]
in which case the $k$-th coordinate changes from $0$ to $1$, and
\[
i=a_k+n_1
\quad\Longleftrightarrow\quad
\clocksign(i)=-,
\]
in which case it changes from $1$ to $0$. Since $k$ was arbitrary, this proves the first assertion.

For every $i$, the cyclic interval
$\{a_k,\ldots,a_k+n_1-1\}$ contains exactly one of $i$ and
$i+n_1$. Hence
\[
(A_{i+n_1})_k=1-(A_i)_k.
\]
Since this holds for every $k$, we obtain
\[
A_{i+n_1}=\mathbf 1-A_i.
\]
\end{proof}

\begin{lemma}[Adjacency]\label{lem:app_adjacency}
\begin{enumerate}[label=(\alph*)]
\item For every clock $\clocksign$, every ray labeling $h=(h_i)$
of size $n_1$, and every $i,j\in\mathbb Z/2n_1\mathbb Z$, the patterns $A_i$ and $A_j$ differ in exactly $d(i,j)$ coordinates.
Consequently, the $2n_1$ patterns $A_i$ are distinct, and $A_i$ and $A_j$ differ in exactly one coordinate if and only if $j=i\pm1$.

\item Let
$\theta\in\mathcal P^{\mathrm{nd}}_{2,n_1,n_2}$, with rays $r_i$,
regions $R_i$, ray labeling $h=(h_i)$, and sign function
$\clocksign_\theta$ as in \Cref{app:clock_construction}. Then the
activation pattern of the maximal cell $\overline{R_i}$ of
$\mathcal C(\theta)$ is $A_i$. Consequently,
\[
\mathcal{A}_{\theta}=P(\clocksign_\theta,h)
\qquad\text{and}\qquad
|\mathcal{A}_{\theta}|=2n_1.
\]
\end{enumerate}
\end{lemma}

\begin{proof}
(a) If $i=j$, the assertion is immediate. Otherwise, after possibly
interchanging $i$ and $j$ and choosing the shorter direction around
the cycle, we may write $j=i+d$ with $d=d(i,j)\leq n_1$. Passing from $A_i$ to $A_j$ through
$A_{i+1},\ldots,A_{i+d-1}$ changes, by
\Cref{lem:app_crossing}, the coordinates
$h_{i+1},\ldots,h_{i+d}$, one at a time. These $d$ coordinates are distinct, because $h_\ell=h_{\ell'}$ only if $\ell-\ell'\in\{0,n_1\}$, while two of the indices $i+1,\ldots,i+d$ differ by at most $d-1<n_1$. Thus, each of these $d$ coordinates changes exactly once, while every other coordinate remains unchanged. Consequently, $A_i$ and $A_j$
differ in exactly $d=d(i,j)$ coordinates. If $i\neq j$, then $d(i,j)\geq1$, so $A_i\neq A_j$. Hence the $2n_1$ patterns $A_i$ are distinct. Moreover, $A_i$ and $A_j$ differ in exactly one coordinate if and only if $d(i,j)=1$, which is equivalent to $j=i\pm1$.

(b) The $n_1$ lines $H_1,\ldots,H_{n_1}$ are distinct and each contains two opposite rays, so there are $2n_1$ rays. For each $j$, every line other than $H_{h_j}$ meets the open half-plane on the counterclockwise side of $r_j$ in exactly one ray. Hence exactly $n_1-1$ rays lie strictly between $r_j$ and $-r_j$ in counterclockwise order, so $-r_j=r_{j+n_1}$. This half-plane is the union of the regions $R_j,\ldots,R_{j+n_1-1}$ and the rays between them. Let $k=h_j$. Neuron $k$ is active on one of the two open half-planes bounded by $H_k$, and by the definition of $\clocksign_\theta$ it is active on $R_j$ if and only if $\clocksign_\theta(j)=+$. Thus it is active exactly on the regions $R_{a_k},\ldots,R_{a_k+n_1-1}$, where
$a_k\in\{j,j+n_1\}$ is the index with $\clocksign_\theta(a_k)=+$. This is the definition of $(A_i)_k$. Since this holds for every $k$, the activation pattern on $R_i$ is $A_i$. The maximal cells of $\mathcal C(\theta)$ are the closures of the regions, and the pattern of a maximal cell is its pattern on the interior. Therefore $\mathcal{A}_{\theta}=P(\clocksign_\theta,h)$, and by
(a) these $2n_1$ patterns are distinct.
\end{proof}

The dihedral group $D_{2n_1}$ acts on clocks by
$(g\cdot\clocksign)(i)=\clocksign(g^{-1}(i))$
(\Cref{app:clock_construction,app:clock_action}), where
$r(i)=i+1$ and $f_c(i)=c-i$. Since every $g\in D_{2n_1}$ preserves antipodal pairs, it also acts on ray labelings by
\[
(g\cdot h)_i:=h_{g^{-1}(i)}.
\]
For $\pi\in S_{n_1}$, define the relabeled sequence $\pi h$ by
\[
(\pi h)_i:=\pi(h_i).
\]
We write $\pi\cdot A$ for the pattern satisfying
$(\pi\cdot A)_{\pi(k)}=A_k$, as in
\Cref{def:single_actionvationpattern_permutation}.

\begin{lemma}[Symmetries]\label{lem:app_dihedral}
Let $\clocksign,\clocksign'$ be clocks and let
$h=(h_i),h'=(h_i')$ be ray labelings of size $n_1$.
\begin{enumerate}[label=(\alph*)]
\item For $\pi\in S_{n_1}$ and $g\in D_{2n_1}$,
\[
P(\clocksign,\pi h)=\pi\cdot P(\clocksign,h),
\qquad
P(g\cdot\clocksign,g\cdot h)=P(\clocksign,h).
\]
In particular, $P(-\clocksign,h)=P(\clocksign,h)$.
\item There is $\pi\in S_{n_1}$ with
\[
\pi\cdot P(\clocksign,h)=P(\clocksign',h')
\]
if and only if $\clocksign'=g\cdot\clocksign$ for some
$g\in D_{2n_1}$.
\item For every clock $\clocksign$ and ray labeling
$h=(h_i)$ with $h_0=1$, choose distinct rays
\[
r_0,r_1,\ldots,r_{n_1-1}
\]
contained in a common open half-plane and listed in
counterclockwise order, and set
\[
r_{j+n_1}=-r_j.
\]
For each $j$, choose a nonzero vector $u_j\in r_j$, and let $J$ denote counterclockwise rotation by $\pi/2$. Identifying the signs with $\pm1$, define the first-layer rows by
\[
w_{h_j}=\clocksign(j)Ju_j,
\qquad j=0,\ldots,n_1-1.
\]
Then
\[
w_{h_j}^{\perp}=\operatorname{span}(u_j),
\]
and the sign $\clocksign(j)$ makes neuron $h_j$ active on the prescribed side of $r_j$. Together with any $W_2$, this gives a parameter in $\mathcal P^{\mathrm{nd}}_{2,n_1,n_2}$ having
clock $\clocksign$ and ray labeling $h$. Consequently, $\alpha$ is surjective.
\end{enumerate}
\end{lemma}

\begin{proof} 
\begin{enumerate}[label=(\alph*),leftmargin=*]
    \item Let $a_k'$ denote the positive ray indexing corresponding to the $k$th neuron, i.e., $(\pi h)_{a_k'}=k$. Hence, $(\pi h)_{a_k'}=\pi (h_{a_k'})=k$, which implies $h_{a_k'}=\pi^{-1}(k)$, and thus $a_k'=a_{\pi^{-1}(k)}$. Let $A_i'\in P(\clocksign,\pi h)$ be the pattern corresponding to region $R_i$, and let $A_i\in P(\clocksign,h)$ be the original pattern corresponding to $R_i$. Then
\[
(A_i')_k=(A_i)_{\pi^{-1}(k)}
=(\pi\cdot A_i)_k
\]
for every $k\in\{1,\ldots,n_1\}$. Hence
$A_i'=\pi\cdot A_i$, and therefore $P(\clocksign,\pi h)
=\pi\cdot P(\clocksign,h).$

For the second identity, set
$\clocksign'=g\cdot\clocksign$ and $h'=g\cdot h$, and let
$a_k',A_i'$ denote the corresponding data. Then
$h_i'=k$ and $\clocksign'(i)=+$ if and only if
$h_{g^{-1}(i)}=k$ and $\clocksign(g^{-1}(i))=+$, so
$a_k'=g(a_k)$. If $g=r^{j_0}$, then
$i-a_k'=(i-j_0)-a_k$, and therefore $A_i'=A_{i-j_0}$.
If $g=f_c$, then
\[
i-a_k'\in\{0,\ldots,n_1-1\}
\]
if and only if
\[
(c-i-1)-a_k\in\{n_1,\ldots,2n_1-1\}.
\]
Thus $A_i'=\mathbf 1-A_{c-i-1}=A_{c-i-1+n_1}$ by
\Cref{lem:app_crossing}. In either case the two pattern sets coincide.
\item Suppose first that $\clocksign'=g\cdot\clocksign$. Because
$g\cdot h$ and $h'$ are both ray labelings, there is a unique
$\pi\in S_{n_1}$ such that
\[
h'=\pi(g\cdot h),
\qquad\text{that is,}\qquad
h_i'=\pi\bigl(h_{g^{-1}(i)}\bigr).
\]
Part (a) then gives
\[
P(\clocksign',h')
=\pi\cdot P(\clocksign',g\cdot h)
=\pi\cdot P(\clocksign,h).
\]

Conversely, suppose
$\pi\cdot P(\clocksign,h)=P(\clocksign',h')$.
Write $A_i$ and $A'_i$ for the corresponding labeled patterns.
We use two observations. First, by \Cref{lem:app_crossing},
$A_{i-1}$ and $A_i$ differ in exactly one coordinate, namely
$h_{i}$; $\clocksign(i)=+$ exactly when that coordinate goes
from $0$ in $A_{i-1}$ to $1$ in $A_i$, and $\clocksign(i)=-$
when it goes from $1$ to $0$. Second, an entry permutation only
relabels coordinates, so it sends patterns differing in exactly one
coordinate to patterns differing in exactly one coordinate, and a
change from $0$ to $1$ to a change from $0$ to $1$.

By \Cref{lem:app_adjacency}(a), the patterns are distinct, so $\pi$
induces a bijection $\psi$ of $\mathbb Z/2n_1\mathbb Z$ with
$\pi\cdot A_i=A'_{\psi(i)}$. By the second observation and
\Cref{lem:app_adjacency}(a), $\psi$ preserves adjacency in both
directions. Since $n_1\geq2$, it is an automorphism of the
$2n_1$-cycle. The automorphisms of this cycle form the dihedral group,
so $\psi$ is either a rotation $i\mapsto i+k$ or a reflection
$i\mapsto k-i$ for some $k\in\mathbb Z/2n_1\mathbb Z$.

It remains to check that the orientations move with the patterns.
If $\psi(i)=i+k$, the transition from $A_{i-1}$ to $A_i$
corresponds, after applying $\pi$, to the transition from
$A'_{i+k-1}$ to $A'_{i+k}$, in the same direction. By the two
observations, $\clocksign'(i+k)=\clocksign(i)$ for all $i$,
so $\clocksign'=r^k\cdot\clocksign$.

If $\psi(i)=k-i$, the same transition corresponds to the transition
from $A'_{k-i+1}$ to $A'_{k-i}$, in the opposite direction.
Thus a change from $0$ to $1$ becomes a change from $1$ to $0$
when read in increasing cyclic order, and
\[
\clocksign'(k+1-i)=-\clocksign(i)=\clocksign(i+n_1).
\]
Setting $c=k+n_1+1$, this says
$\clocksign'(j)=\clocksign(c-j)$ for all $j$, so
$\clocksign'=f_c\cdot\clocksign$.
In both cases the clocks differ by a dihedral action.
\item For each $j=0,\ldots,n_1-1$,
\[
H_{h_j}
=
w_{h_j}^{\perp}
=
(Ju_j)^{\perp}
=
\operatorname{span}(u_j).
\]
Since the rays $r_0,\ldots,r_{n_1-1}$ lie in a common
open half-plane and are distinct, these $n_1$ lines are distinct. In particular, the rows of $W_1$ are nonzero and no two are scalar multiples of one another.

For sufficiently small $\varepsilon>0$, the vector
\[
x_\varepsilon=u_j+\varepsilon Ju_j
\]
lies in the region $R_j$ immediately counterclockwise from $r_j$. Using $u_j^\top Ju_j=0$, we obtain
\[
\begin{aligned}
w_{h_j}^{\top}x_\varepsilon
&=\clocksign(j)(Ju_j)^\top
   \bigl(u_j+\varepsilon Ju_j\bigr)\\
&=\clocksign(j)\varepsilon\lVert u_j\rVert^2.
\end{aligned}
\]
Therefore neuron $h_j$ is active on $R_j$ if and only if $\clocksign(j)=+$. The antipodal relations
\[
r_{j+n_1}=-r_j,
\qquad
h_{j+n_1}=h_j,
\qquad
\clocksign(j+n_1)=-\clocksign(j)
\]
give the prescribed labels and signs on the remaining rays. Thus the constructed parameter has clock $\clocksign$ and ray labeling $h$. Since $h_0=1$ and $r_0=\mathbb R_{>0}u_0\subseteq H_1$, it also satisfies the normalization used in the clock construction. Consequently, $\alpha$ is surjective.
\end{enumerate}
\end{proof}

\begin{figure}[htbp]
\centering
\includegraphics[width=.84\linewidth]{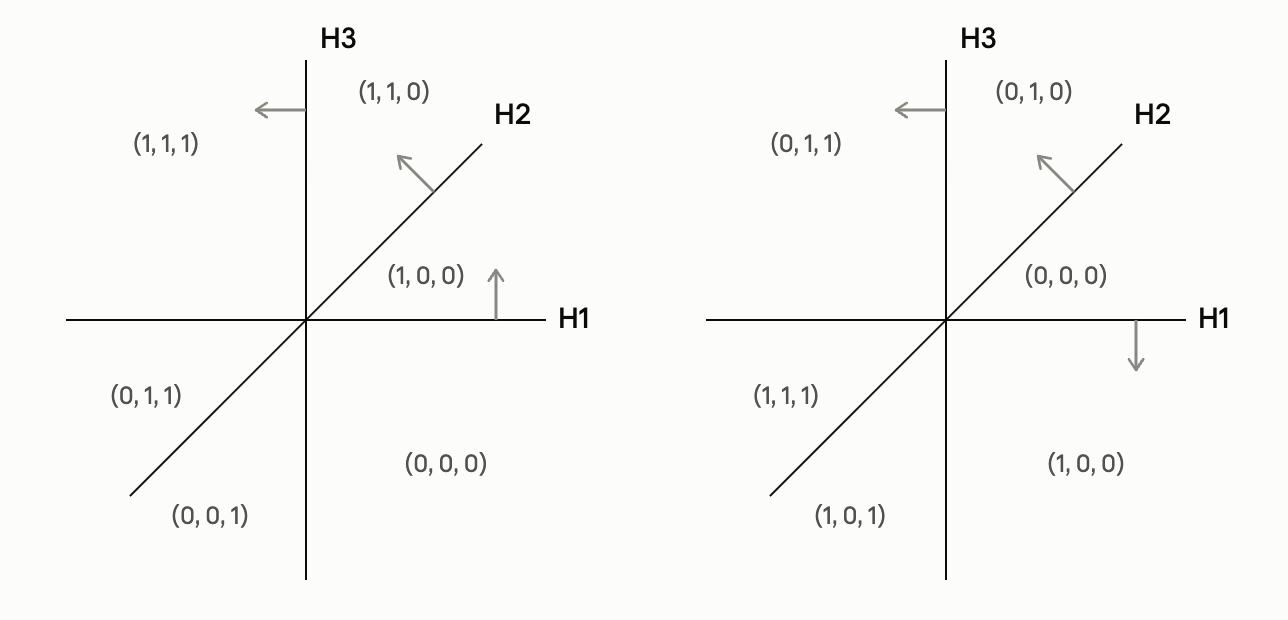}
\caption{The same three lines under two orientations of $H_1$. Their activation-pattern sets agree after exchanging coordinates $2$ and $3$, although the two weight matrices are not row permutations of one another.}
\label{fig:app_clock_orientation}\label{fig:h1-orientation}
\end{figure}

\paragraph{Pattern equivalence without permuting the weight rows.}
A permutation of neurons induces the corresponding permutation of entries in every activation pattern. The converse fails: two weight matrices whose
pattern sets differ by an entry permutation need not be related by a permutation of neurons. For the two arrangements in \Cref{fig:app_clock_orientation}, one may take
\[
W_1=\begin{pmatrix}0&1\\-1&1\\-1&0\end{pmatrix},
\qquad
W'_1=\begin{pmatrix}0&-1\\-1&1\\-1&0\end{pmatrix}.
\]
Only the first neuron reverses orientation. Nevertheless, the sets
\[
\begin{aligned}
\mathcal{A}_{\theta}&=\{(1,1,0),(1,1,1),(0,1,1),(0,0,1),(0,0,0),(1,0,0)\},\\
\mathcal{A}_{\theta'}&=\{(0,1,0),(0,1,1),(1,1,1),(1,0,1),(1,0,0),(0,0,0)\}
\end{aligned}
\]
are related by the simultaneous permutation $(2\ 3)$; here a binary
word abbreviates its corresponding three-entry tuple. This does not
mean that reversing an arbitrary neuron always preserves pattern
equivalence; it is a property of this example.

\newpage 

\begin{figure}[htbp]
\centering
\includegraphics[width=.66\linewidth]{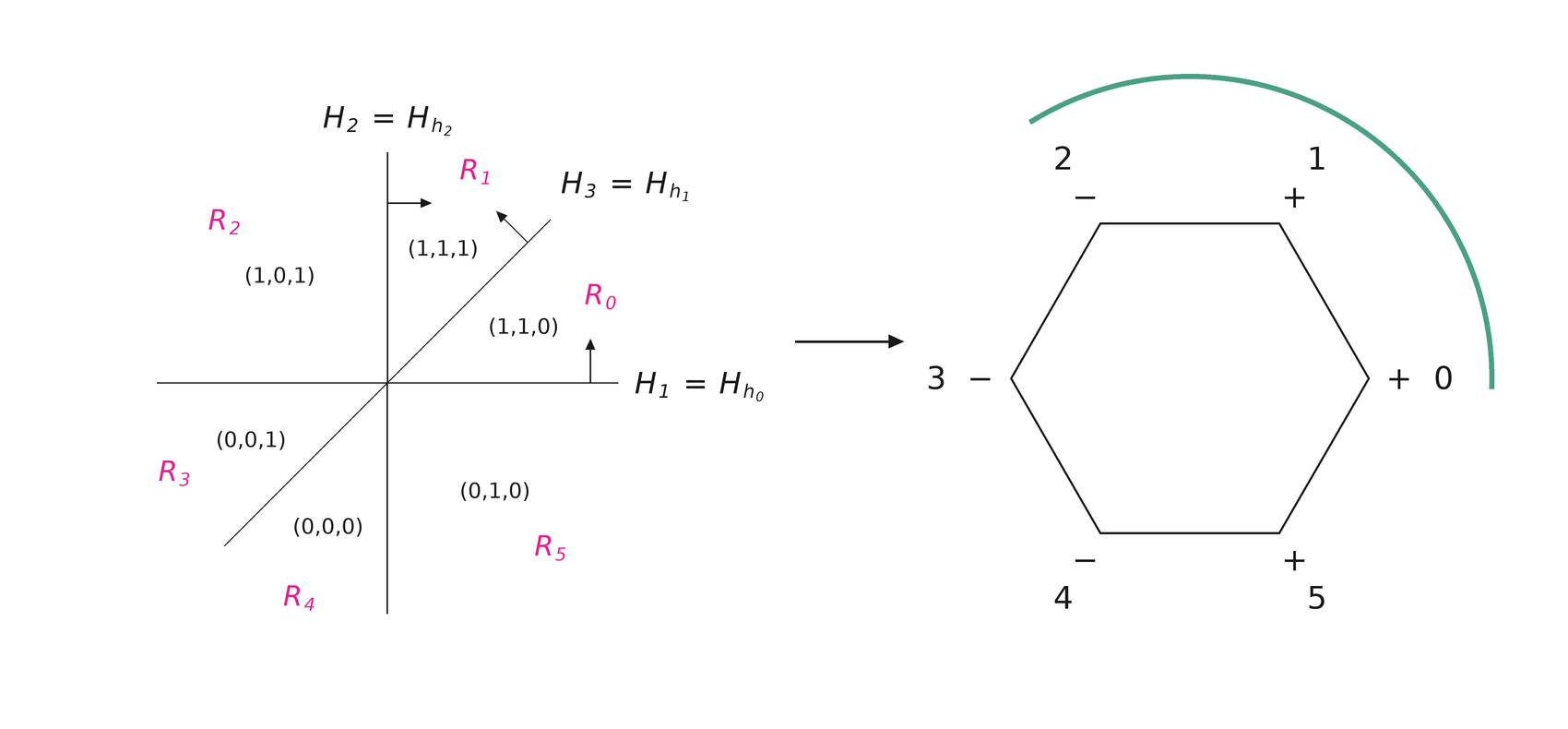}
\par\medskip
\includegraphics[width=.66\linewidth]{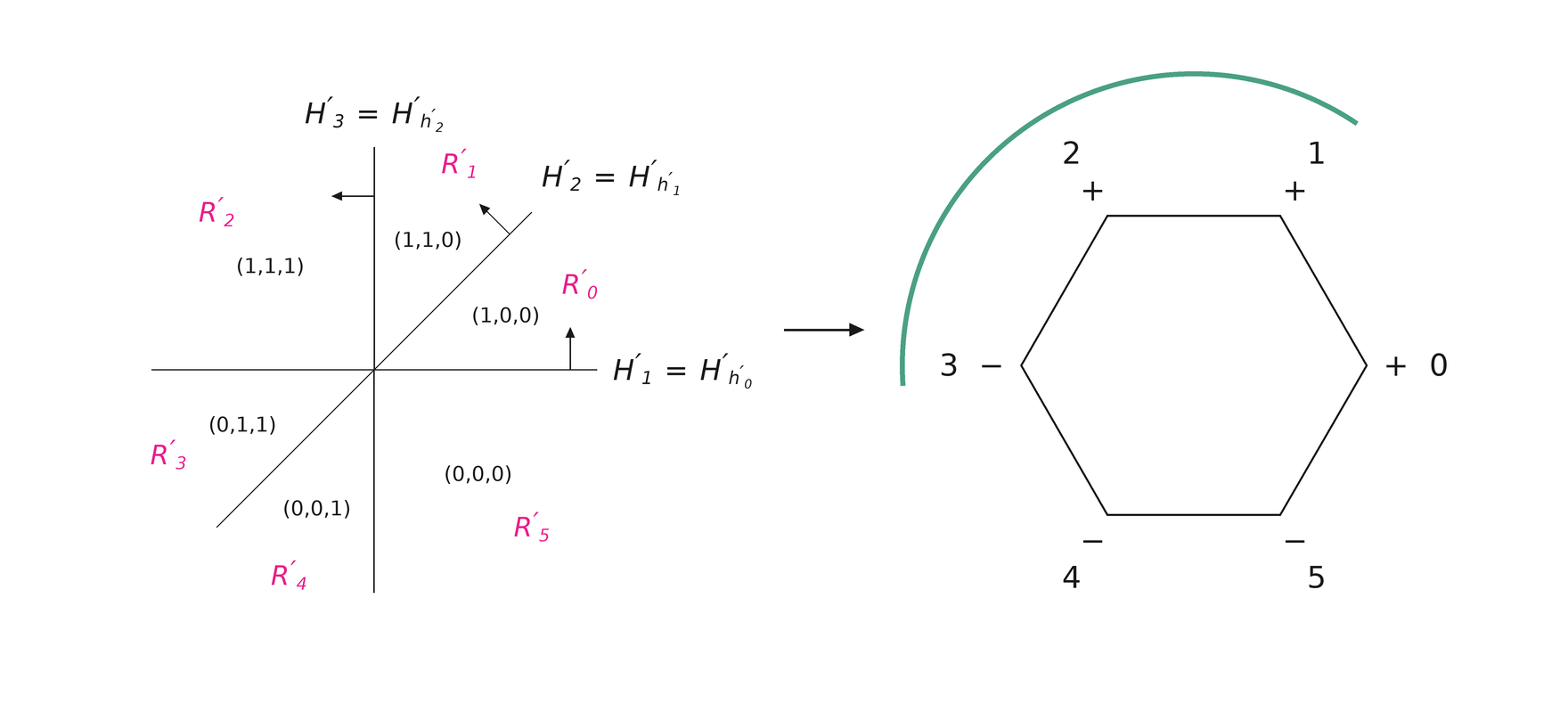}
\caption{Two oriented arrangements whose clocks differ by one counterclockwise click. The line containing $r_i$ is $H_{h_i}$, and the line containing $r'_i$ is $H'_{h'_i}$. The coordinate permutation is $(1\ 2)$, and the corresponding regions are $R_i$ and $R'_{i+1}$.}
\label{fig:app_clock_rotation_example}
\end{figure}

\paragraph{A one-click rotation.}
The two arrangements in \Cref{fig:app_clock_rotation_example} have
\[
\begin{aligned}
(\clocksign(0),\ldots,\clocksign(5))&=(+,+,-,-,-,+),\\
(\clocksign'(0),\ldots,\clocksign'(5))&=(+,+,+,-,-,-).
\end{aligned}
\]
Thus $\clocksign'=r\cdot\clocksign$. The ray labelings are
\[
(h_0,h_1,h_2)=(1,3,2),
\qquad
(h'_0,h'_1,h'_2)=(1,2,3).
\]
The permutation $\pi=(1\ 2)$ satisfies
$h'=\pi(r\cdot h)$, and hence
\[
A'_{i+1}=\pi\cdot A_i\qquad (i\in\mathbb Z/6\mathbb Z).
\]
In particular, $A_0=(1,1,0)$ is sent to $A'_1=(1,1,0)$, while $A_2=(1,0,1)$ is sent to $A'_3=(0,1,1)$. This exhibits the same coordinate permutation on every region, as required by pattern equivalence.

\newpage

\subsection{Proof of the clock characterization}
\label{app:clock_classification_proof}

\begin{proof}[Proof of \Cref{thm:clock_pattern_equivalence}]
Choose associated clocks $\clocksign,\clocksign'$ with ray labelings
$h=(h_i)$ and $h'=(h_i')$ as in
\Cref{app:clock_construction}.
By \Cref{lem:app_adjacency}(b),
\[
\mathcal{A}_{\theta}=P(\clocksign,h),\qquad
\mathcal{A}_{\theta'}=P(\clocksign',h').
\]
By \Cref{def:patternequivalentparameter}, pattern equivalence means
that these sets agree after a simultaneous permutation of their
coordinates. By \Cref{lem:app_dihedral}(b), this holds precisely
when $\clocksign$ and $\clocksign'$ lie in the same dihedral orbit.
Replacing either clock by its shadow is a rotation by $n_1$ vertices
(\Cref{app:clock_action}), so it does not change that orbit.
Finally, every clock is realized by a parameter, by
\Cref{lem:app_dihedral}(c). Thus the induced map from
pattern-equivalence classes to $\clockset{n_1}/D_{2n_1}$ is both
injective and surjective.
\end{proof}

\subsection{Fixed clocks and the counting formula}
\label{app:clock_counting_proof}

We have reduced the problem of counting pattern-equivalence classes to counting clock aisles. To do this, we apply Burnside's lemma to the dihedral action on clocks. We first determine which clocks are fixed by each rotation and reflection, and then apply the orbit-counting identity.

For a group element $g$, write $\operatorname{Fix}(g)$ for the set of clocks fixed by $g$.

\begin{proposition}[Fixed clocks]\label{prop:app_fixed_points}
Let $n_1\geq1$.
\begin{enumerate}[label=(\alph*)]
\item For $k\in\mathbb Z/2n_1\mathbb Z$ and $d:=\gcd(2n_1,k)$, the rotation $r^k$ fixes $2^{d/2}$ clocks if $d\nmid n_1$ and no clock if $d\mid n_1$.
\item The reflection $f_c$ fixes $2^{\lceil n_1/2\rceil}$ clocks if $c+n_1$ is odd and no clock if $c+n_1$ is even. Thus exactly $n_1$ of the $2n_1$ reflections have fixed clocks: for $n_1$ odd those with $c$ even, whose axis passes through two opposite vertices, and for $n_1$ even those with $c$ odd, whose axis passes through the midpoints of two opposite edges.
\end{enumerate}
\end{proposition}

\begin{proof}
For $n_1=1$, the two clocks are $(+,-)$ and $(-,+)$.
Both are fixed by $r^0$ and $f_0$, and neither is fixed by $r$ or $f_1$, giving the stated counts. Assume henceforth that $n_1\geq2$. We compute the fixed clocks for each type of group element. A clock fixed by $g\in D_{2n_1}$ must be constant on each orbit of $\langle g\rangle$ acting on the vertices, and must satisfy
\[
\clocksign(i+n_1)=-\clocksign(i)
\qquad\text{for every vertex }i.
\]
The antipodal map $\tau(i)=i+n_1$ is an involution that commutes with $g$, so it permutes the $\langle g\rangle$-orbits. If $\tau$ fixes some orbit setwise, that orbit contains an antipodal pair $\{i,i+n_1\}$, forcing $\clocksign(i)=-\clocksign(i)$, which is impossible. Otherwise $\tau$ pairs the orbits; each pair contributes one free binary choice. Thus, if there are $q_g$ vertex orbits and none is fixed setwise by $\tau$, there are $2^{q_g/2}$ fixed clocks.

(a) Put $d=\gcd(2n_1,k)$. The rotation $r^k$ decomposes the vertices into $d$ disjoint cycles of length $2n_1/d$. A cycle contains an antipodal pair if and only if there is an integer
$m$ with
\[
i+mk\equiv i+n_1\pmod{2n_1},
\quad\text{that is,}\quad mk\equiv n_1\pmod{2n_1}.
\]
This congruence is solvable exactly when $d\mid n_1$.
In that case every cycle contains an antipodal pair and no clock is fixed. If $d\nmid n_1$, the antipodal map pairs the $d$ cycles into $d/2$ pairs; choosing a sign on one cycle of each pair determines its partner, giving $2^{d/2}$ fixed clocks.

(b) A reflection is an involution whose vertex orbits are its fixed vertices and its transposed pairs. Its axis passes either through two opposite vertices, giving $n_1+1$ vertex orbits, or through the midpoints of two opposite edges, giving $n_1$ vertex orbits. Writing the reflection as $f_c(i)=c-i$, it maps a vertex to its own antipode exactly when
\[
c-i\equiv i+n_1\pmod{2n_1},
\quad\text{that is,}\quad 2i\equiv c-n_1\pmod{2n_1}.
\]
This is solvable precisely when $c$ and $n_1$ have the same parity, equivalently when $c+n_1$ is even. For those $n_1$ reflections, some vertex orbit contains an antipodal pair and there are no fixed clocks.

For the other $n_1$ reflections, the antipodal map pairs all the vertex orbits. If $n_1$ is odd, then $c$ is even, the axis passes through opposite vertices, and the $n_1+1$ vertex orbits give $2^{(n_1+1)/2}$ fixed clocks. If $n_1$ is even, then $c$ is odd, the axis passes through opposite edge midpoints, and the $n_1$ vertex orbits give $2^{n_1/2}$ fixed clocks. Both cases give
$2^{\lceil n_1/2\rceil}$, as claimed.
\end{proof}

\begin{proof}[Proof of \Cref{thm:clock_counting_formula}]
By \Cref{thm:clock_pattern_equivalence}, the number of pattern-equivalence classes equals the number of orbits of $D_{2n_1}$ acting on $\clockset{n_1}$. Therefore, by the orbit-counting identity,
\[
N(n_1)
=
\frac{1}{|D_{2n_1}|}
\sum_{g\in D_{2n_1}}
|\operatorname{Fix}(g)|
=
\frac{1}{4n_1}
\sum_{g\in D_{2n_1}}
|\operatorname{Fix}(g)|.
\]

For a rotation $r^k$, $k\in\mathbb Z/2n_1\mathbb Z$, let $d:=\gcd(2n_1,k)$. By \Cref{prop:app_fixed_points}(a), the number of fixed clocks is
\[
|\operatorname{Fix}(r^k)|
=
\begin{cases}
2^{d/2}, & d\nmid n_1,\\
0, & d\mid n_1.
\end{cases}
\]
The number of $k\in\mathbb Z/2n_1\mathbb Z$ with $\gcd(2n_1,k)=d$ is $\phi(2n_1/d)$. Put $e:=2n_1/d$. Then $d\nmid n_1$ if and only if $e$ is odd, in which case $e\mid n_1$ and $d/2=n_1/e$. Summing over all rotations therefore gives
\[
\sum_{\substack{e\mid n_1\\e\text{ odd}}}
\phi(e)\,2^{n_1/e}.
\]

By \Cref{prop:app_fixed_points}(b), among the $2n_1$ reflections, exactly
$n_1$ fix $2^{\lceil n_1/2\rceil}$
clocks each, while the remaining $n_1$
fix no clocks: for $n_1$ odd the fixing reflections are the $f_c$ with $c$ even, whose axis passes through two opposite vertices, and for $n_1$ even they are the $f_c$ with $c$ odd, whose axis passes through the midpoints of two opposite edges.

Combining the contributions from rotations
and reflections, we obtain
\[
N(n_1)=
\frac{1}{4n_1}
\left(
\sum_{\substack{e\mid n_1\\ e\text{ odd}}}
\phi(e)\,2^{n_1/e}
+
n_1\,2^{\lceil n_1/2\rceil}
\right).
\]
\end{proof}

\begin{table}[ht]
\centering
\caption{Values of the formula in \Cref{thm:clock_counting_formula}.}
\label{tab:clock_counts}
\begin{tabular}{c|cccccccc}
$n_1$ & 1 & 2 & 3 & 4 & 5 & 6 & 7 & 8\\
\hline
$N(n_1)$ & 1 & 1 & 2 & 2 & 4 & 5 & 9 & 12
\end{tabular}
\end{table}

\subsection{Zero and repeated neurons}
\label{app:clock_degenerate}

The clock construction can also describe certain degenerate parameters. Suppose that the nonzero rows of $W_1$ determine $r$ distinct directions, where positive scalar multiples are identified, and that no two nonzero rows are negative scalar multiples of one another. Removing the zero rows and retaining one representative of each direction gives a reduced arrangement of $r$ distinct oriented lines. If $r\geq1$, this arrangement has an associated clock of size $r$.

The original activation-pattern set is recovered from the reduced labeled pattern set by repeating each coordinate according to the number of neurons in that direction and inserting a zero coordinate for each zero row. Thus $N(r)$ counts the possible reduced pattern-equivalence classes, but does not by itself count the classes before reduction. One must also record the multiplicities attached to the lines. Different multiplicity assignments can become equivalent under symmetries of the reduced clock, so their contribution depends on the particular clock.

\begin{example}
\label{ex:clock_degenerate_symmetry}
Consider $n_1=4$, with no zero rows and exactly three distinct nonzero directions. The multiplicities must be $2,1,1$: one of the three directions occurs twice.

There are $N(3)=2$ reduced clock classes, represented by
\[
\clocksign_1=(+,+,+,-,-,-),
\qquad
\clocksign_2=(+,-,+,-,+,-),
\]
where the entries are listed in cyclic order.
Label the antipodal pairs by
$h_{0}=1$, $h_{1}=2$, and $h_{2}=3$.

For $\clocksign_1$, reflection by $f_2(i)=2-i$
exchanges the first and third lines and fixes the second. Thus repeating the first or third coordinate gives equivalent pattern sets. Repeating the second coordinate gives a different class: the numbers of ones in the six patterns, listed in increasing order, are
\[
\begin{array}{c|c}
\text{Repeated coordinate} & \text{Numbers of ones}\\
\hline
1\text{ or }3 & 0,1,2,2,3,4\\
2 & 0,1,1,3,3,4.
\end{array}
\]
Since coordinate permutations preserve these numbers, the two classes are distinct.

For $\clocksign_2$, rotation by two vertices preserves the clock and cyclically permutes its three antipodal pairs. Consequently, all three choices of the repeated coordinate give equivalent pattern sets.

The first reduced clock class therefore contributes two classes, while the second contributes one. These contributions cannot overlap, because deleting repeated coordinates recovers the reduced pattern-equivalence class. Hence there are exactly three classes under these restrictions, illustrating why counting multiplicities does not amount to multiplying $N(r)$ by a common factor.
\end{example}

The exclusion of oppositely directed neurons ensures that each line has a single active side and hence a single crossing sign at each ray. Allowing both orientations on the same line would require additional data beyond the clock used here. If $r=0$, all rows are zero and the only activation pattern is the all-zero vector.

\bibliographystyle{plain}
\bibliography{refs}
\bigskip
\textsc{University of Wisconsin, Madison; Department of Mathematics; 320 Van Vleck Hall; Madison, WI 53706}

\textit{Email address:} \texttt{blfinkel@wisc.edu}
\bigskip

\textsc{University of California, Riverside; Department of Mathematics; 262B Skye Hall; Riverside, CA 92521}

\textit{Email address:} \texttt{yaoyingf@ucr.edu}
\bigskip

\textsc{University of Wisconsin, Madison; Department of Mathematics; 517 Van Vleck Hall; Madison, WI 53706}

\textit{Email address:} \texttt{cwu367@wisc.edu}

\end{document}